\documentclass[11pt]{article}
\usepackage{graphicx}
\usepackage{amssymb,latexsym}
\usepackage{epstopdf}
\newtheorem{theorem}{Theorem}[section]
\newtheorem{prop}[theorem]{Proposition}
\newtheorem{proposition}[theorem]{Proposition}
\newtheorem{remark}[theorem]{Remark}
\newtheorem{lemma}[theorem]{Lemma}
\newtheorem{claim}[theorem]{Claim}
\newtheorem{defi}[theorem]{Definition}

\newtheorem{cor}[theorem]{Corollary}
\newtheorem{conj}[theorem]{Conjecture}

\newcommand{\forme}[1]{}
\newcommand{\pf}{\noindent{\em Proof: }}
\newcommand{\pfmainthm}{\noindent{\em Proof of Theorem \ref{b-i-small-thm}: }}
\newcommand{\HH}{{\texttt h}}
\newcommand{\epf}{\hfill\hbox{\rule{3pt}{6pt}}}

\newcommand{\g}{\mathcal{G}}
\newcommand{\T}{\mathcal{T}}
\newcommand{\TT}{{\texttt t}}
\newcommand{\Gap}{{\texttt {Gap}}}
\newcommand{\Len}{{\texttt {Len}}}
\newcommand{\I}{\mathcal{I}}
\newcommand{\J}{\mathcal{J}}
\newcommand{\N}{{\mathbb N}}
\begin{document}

\title{There are only finitely many distance-regular graphs
of fixed valency greater than two}
\author{
{\bf S.~Bang}
\\Department of Mathematics, Pusan National University \\Geumjeong Gu, Busan 609-735 Republic of Korea\\
e-mail: sjbang3@pusan.ac.kr\\
{\bf A.~Dubickas}
\\Department of Mathematics and Informatics, Vilnius University\\
Naugarduko 24 LT-03225 Vilnius, Lithuania\\
e-mail: arturas.dubickas@mif.vu.lt\\
{\bf J.~H.~Koolen\footnote{Corresponding author}}\\ Pohang Mathematics Institute and Department of Mathematics, POSTECH\\
Hyoja-dong, Namgu, Pohang 790-784 Korea\\
e-mail: koolen@postech.ac.kr \\
{\bf V.~Moulton}
\\ School of Computing Sciences, University of East Anglia,\\ Norwich, NR4 7TJ, UK.
\\ e-mail: vincent.moulton@cmp.uea.ac.uk }
\date{\today}
\maketitle

\begin{abstract}
In this paper we prove the Bannai-Ito conjecture, namely that
there are only finitely many distance-regular graphs of fixed
valency greater than two.
\end{abstract}
\newpage
\footnotesize
\tableofcontents
\normalsize
\newpage

\section{Introduction}\label{intro}

A finite, connected graph $\Gamma$ with vertex set $V(\Gamma)$ and
path-length distance $d$ is said to be {\em distance-regular} if,
for any vertices $x,y \in V(\Gamma)$ and any integers $1 \le i,j
\leq \max \{d(z,w)\,:\,z,w \in V(\Gamma)\}$, the number of
vertices at distance $i$ from $x$ and distance $j$ from $y$
depends only on $i, j$ and $d(x,y)$, independent of the choice of
$x$ and $y$. Many distance-regular graphs arise from classical
objects, such as the Hamming graphs, the Johnson graphs, the
Grassmann graphs, the bilinear forms graphs, and the dual polar
graphs amongst others. In particular, distance-regular graphs give
a framework to study these classical objects from a combinatorial
point of view. In addition, distance-regular graphs and
association schemes give an algebraic-combinatorial framework to
study, for example, codes and designs \cite{bcn,delsarte}.

In their 1984 book, E.~Bannai and T.~Ito conjectured that there
are only finitely many distance-regular graphs of fixed valency
greater than two (cf. \cite[p.237]{banito}). In this paper we
prove that their conjecture holds:

\begin{theorem}\label{b-i-thm}
There are only finitely many distance-regular graphs of fixed
valency greater than two.
\end{theorem}

{\bf History}{\ \\} A distance-transitive graph is a connected
graph $\Gamma$ such that for every four  (not necessarily
distinct) vertices $x,y,u,v$ in $V(\Gamma)$ with $d(x,y) =
d(u,v)$, there exists an automorphism $\tau$ of $\Gamma$ such that
$\tau(x) = u$ and $\tau(y) =v$ both hold. It is straight-forward
to see that distance-transitive graphs are distance-regular
graphs. In \cite{cameron,dtg}, P.~J.~Cameron, C.~E.~Praeger,
J.~Saxl and G.~M.~Seitz proved that there are only finitely many
finite distance-transitive graphs of fixed valency greater than
two. They did this by applying Sims' conjecture \cite{sim} for
finite permutation groups (i.e. that there exists an integral
function $f$ such that $\left|G_{x}\right| \leq f(d_{G_x})$ holds,
where, for $G$ a primitive permutation group acting on a finite
set $\Omega$, $G_{x}$ denotes the stabilizer of $x$, $x \in
\Omega$, and $d_{G_x}$ denotes the length of any $G_x$-orbit in
$\Omega \setminus \{x\}$), which they also showed to hold by using
the classification of the finite simple groups (in \cite{dtg} they
gave a proof without many details, and in \cite{cameron} Cameron
worked out a detailed proof with an explicit diameter bound).

Note that for small diameter there are many distance-regular
graphs which are not distance-transitive. On the other hand there
are only five families of distance-regular but not distance-transitive graphs known with
unbounded diameter, namely the Doob graphs \cite{doob} (see also
\cite[p.262]{bcn}), the quadratic forms graphs \cite{egawa} (see
also  \cite[p.290]{bcn}),  the Hemmeter graphs  \cite{hem} and the
Ustimenko graphs \cite{ust} (for both, see also \cite[p.279]{bcn})
and the twisted Grassmann graphs \cite{D-K}. Any member of the
first four families is vertex-transitive, whereas the twisted
Grassmann graphs have exactly two orbits under the full
automorphism group \cite{D-K}.

The first class of  distance-regular graphs for which the
Bannai-Ito conjecture was shown is the class of regular
generalized $n$-gons. Feit and Higman  \cite{feit-higman} (cf.
\cite[Theorem 6.5.1]{bcn}) showed that a regular generalized
$n$-gon has either valency 2 or $n \in \{3, 4,6,8,12\}$. In
addition, R.~M.~Damerell, and E.~Bannai and T.~Ito have
independently shown that there are only finitely many Moore graphs
with valency at least three \cite{banito-73,damerell}.

In the series of papers
\cite{ban-ito-I,ban-ito-88,ban-ito-II,ban-ito-89}, E.~Bannai and
T.~Ito showed that their conjecture holds for valencies $k = 3,
4$, as well as for the special class of bipartite distance-regular
graphs. In \cite{val5} and \cite{val10}, J.~H.~Koolen and
V.~Moulton also showed that the conjecture holds for
distance-regular graphs of fixed valency $k = 5, 6$ or $7$, and
for triangle-free distance-regular graphs of fixed valency $k = 8,
9$ or $10$. More recently, in \cite{rnp}, together with S.~Bang,
they showed that the Bannai-Ito conjecture holds for regular near
polygons and geodetic distance-regular graphs.

The proof for the Bannai-Ito conjecture
that we present in this paper builds upon many of the concepts and
ideas developed in \cite{two-thm,rnp,val5,val10}.

{\bf Structure of the paper}{\ \\}
In Section~\ref{pre}, we
present some definitions and previous results
concerning distance-regular graphs and associated
sequences and related structures,
and in Section~\ref{graphical seq} we
present some properties of certain
generalizations of these sequences.
In Section~\ref{central}, we state without
proof the key result of the paper (Theorem~\ref{b-i-small-thm})
and used this to prove Theorem~\ref{b-i-thm}.
We also present an outline proof of
Theorem~\ref{b-i-small-thm}, before proving it
in Sections~\ref{no-thy-section} to \ref{distri}.
In Section~\ref{10}, we will present an
application of Theorem~\ref{b-i-thm} to
distance-regular graphs of order $(s,t)$.
We conclude in Section~\ref{11} by discussing
some possible future directions.

\section{Preliminaries}\label{pre}

In this section, we review some of the well-known theory of
Christoffel numbers for orthogonal polynomials, interlacing and
distance-regular graphs that will be used in this paper. We refer
the reader to \cite{banito}, \cite{bcn} and \cite{szego} for more
details.

\subsection{Christoffel Numbers}\label{cn}

Let $L_1$ be the arbitrary tridiagonal matrix
defined by
\begin{equation}\label{tri-mtx}
 L_1:= \left\lgroup
 \begin{tabular}{llllllll}
 $\alpha_0$ & $\beta_0$\\
 $\gamma_1$ & $\alpha_1$ & $\beta_1$\\
 & $\gamma_2$ & $\alpha_2$ & $\beta_2$\\
 & & . & . & .\\
 & & & $\gamma_i$ & $\alpha_i$ & $\beta_i$\\
 & & & & . & . & .\\
 & & & & & $\gamma_{n-1}$ & $\alpha_{n-1}$ & $\beta_{n-1}$\\
& & & & &\makebox{\hspace{.324cm}} &
 $\gamma_{n}$ & $\alpha_{n}$
 \end{tabular}
 \right\rgroup,
\end{equation}
where $\alpha_i\geq 0$, $\beta_{i-1},\gamma_i>0$ are real numbers
with $\alpha_0=\gamma_0=\beta_n=0$ and $\gamma_1=1$, and
$\alpha_i+\beta_i+\gamma_i=\beta_0$ holding for all $1\leq i\leq
n$. Let $v_i(x)~~(0\leq i\leq n+1)$ be the polynomials defined
recursively by the equations
\begin{eqnarray}
v_0(x)&:=&1,~v_1(x):=x,\label{v_i-0,1}\\
xv_i(x)&=&
\beta_{i-1}v_{i-1}(x)+\alpha_iv_i(x)+\gamma_{i+1}v_{i+1}(x)~~~(1\leq
i\leq n-1)~,\label{v_i}\\
v_{n+1}(x)&=&(x-\alpha_n)v_n(x)-\beta_{n-1}v_{n-1}(x)\label{v_n+1},
\end{eqnarray}
and $F_i(x)~~(0\leq i\leq n)$ be the monic polynomials defined by
setting $F_0(x):=1$, $F_1(x):=x+1$ and
\[F_i(x):=\gamma_2\cdots
\gamma_i(v_0(x)+v_1(x)+\cdots +v_i(x))~~(2\leq i\leq n).\]
Note that for each $2\leq i\leq n$, the polynomial $F_i(x)$ satisfies the recurrence relation
\begin{equation}\label{F_i-rec}
F_i(x)=(x-\beta_0 +\beta_{i-1}+
\gamma_i)F_{i-1}(x)-\beta_{i-1}\gamma_{i-1}F_{i-2}(x).
\end{equation}
Moreover, by (\ref{v_i-0,1})--(\ref{F_i-rec}), for each $0\le i \le
n$, the polynomials $v_i(x)$ and $F_i(x)$ have degree $i$ and have
exactly $i$ distinct real roots in the closed interval
$[-\beta_0,\beta_0]$ (cf. \cite[Theorem 3.3.1]{szego}). Note that the polynomial
$(x-\beta_0)F_n(x)$ is the minimal polynomial of
the matrix $L_1$.

Now, let $\kappa:=\beta_0$ and define
\begin{eqnarray}
\kappa_i&:=& v_i(\kappa)~~~(0\leq i\leq n ), \mbox{ and }\label{def-k_i}\\
u_i(x)&:=& \frac{v_i(x)}{\kappa_i}~~~(0\leq i\leq
n).\label{vi-ui}
\end{eqnarray}
Put $\kappa:=\kappa_1$. Then the polynomials $u_i(x)~(0\leq i\leq
n)$ satisfy
\begin{eqnarray}
&& u_i(\kappa)=1~(0\leq i\leq n) ;\label{u_i-initial}\\
&& u_0(x)=1,~u_1(x)=\frac{x}{\kappa},
~ xu_i(x)=\gamma_iu_{i-1}(x)+\alpha_iu_i(x)+\beta_iu_{i+1}(x)~(1\leq i\leq n ).
\label{u_i-three-terms}
\end{eqnarray}
The sequence $\left(u_i(x)\right) _{i=0}^{n}$ is called
the {\em standard sequence} of $L_1$, and
if $\theta$ is an eigenvalue of $L_1$, then
the column vector
$(u_0(\theta),u_1(\theta),\ldots,u_{n}(\theta))^T$ is a right
eigenvector of $L_1$ associated to $\theta$, by
(\ref{u_i-three-terms}).

Note also that it follows by (\ref{vi-ui})
that, for each eigenvalue $\theta$ of the matrix $L_1$, the equation
\begin{equation}\label{cnumber}
\sum_{i=0}^{n}\frac{v_i^2(\theta)}{\kappa_i}
=\sum_{i=0}^{n}\kappa_iu_i^2(\theta)~
\end{equation}
holds.

Now, let $\beta_0=\theta_0 >\theta_1>\theta_2>\cdots >\theta_n$ be
the eigenvalues of $L_1$ and, for $i=0,1,\ldots, n$, define
\begin{equation}\label{crn}
m_i:= \left(\frac{ \sum_{j=0}^n \kappa_j } { \sum_{j=0}^n
\frac{v_j^2(\theta_i)}{\kappa_j}} \right)
\end{equation}
as well as the symmetric bilinear form $(\cdot, \cdot)$ on the
polynomial ring ${\mathbb R}[x]$ by
$$
(f,g) := \sum_{i=0}^n m_i f(\theta_i)g(\theta_i).
$$
Then, $(v_i, v_i) \neq 0$ holds for all $0 \le i \le n$, and $(v_i, v_j) = (v_i, v_i) \delta_{i,j}$ holds for all
$0 \le i,j \le n$, where
$\delta_{i,j}$ is the Kronecker delta function on $\N _0 \times \N
_0$, where $\N _0$ is the set of non-negative integers. In particular, it follows that $(v_i)_{i=0}^n$ is a sequence of
orthogonal polynomials with respect to $( \cdot, \cdot)$. Note
that within the theory of orthogonal polynomials, the numbers
$m_i$ are referred to as the {\em Christoffel numbers} of the sequence
$(v_i)_{i=0}^n$ (\cite[Theorem 3.4.1]{szego},
\cite[p.201]{banito}). Analogously, we call the number $m_i$ as
defined in (\ref{crn}), the {\em Christoffel number} of $L_1$
associated with $\theta_i$.

\subsection{Interlacing}\label{interlacing}

We now recall two results stated in \cite{two-thm} that provide us
with some interrelationships between the eigenvalues of the matrix
$L_1$ as defined in (\ref{tri-mtx}). The first generalizes the
well-known Interlacing Theorem~\cite[Theorem 3.3.1]{bcn}, from
which it immediately follows.

\begin{lemma}(\cite[Lemma 3.1]{two-thm}){\ \\}  \label{inter}
Suppose that $A$ is a real $n \times n$ matrix for which there
exists a non-singular matrix $Q$ such that the matrix $Q^{-1} A Q$
is real and symmetric. If $\eta_1\leq \dots \leq \eta_n$ are the
eigenvalues of $A$ and $\theta_1\leq\dots\leq\theta_{n-1}$ are the
eigenvalues of the matrix obtained by removing the $i$th row and
$i$th column of $A$, with $i\in \{1,\ldots,n\}$, then
$$
 \eta_1 \le \theta_1 \le \eta_2 \le \dots \le
 \eta_{n-1} \le \theta_{n-1} \le \eta_n.
$$
\end{lemma}

In particular, since $\beta_i \gamma_{i+1} > 0~(0\leq i\leq n-1)$
and $L_1$ is tridiagonal, it follows that $L_1$ satisfies the
conditions on $A$ given in Lemma \ref{inter}, and therefore the
eigenvalues of $L_1$ must satisfy the inequalities given in this
lemma.

The second result guarantees the existence of
eigenvalues of $L_1$ lying within certain limits.


\begin{lemma}(\cite[Theorem 3.2]{two-thm}){\ \\} \label{two-thm-3.2}
Let $\alpha_i,\beta_i,\gamma_i~(0\leq i\leq n)$ be non-negative
integers satisfying $\alpha_0=\gamma_0=\beta_n=0$,
$\beta_{i-1},\gamma_i>0$, $\alpha_i+\beta_i+\gamma_i=\beta_0$,
$\beta_{i-1}\geq \beta_{i}$ and $\gamma_{i}\geq \gamma_{i-1}$ for
all $1 \le i \le n$, and let $L_1$ be the tridiagonal matrix as
defined in (\ref{tri-mtx}). For each $1\leq i\leq n-1$, let
$\ell(i):=|\{j \,:\,(\gamma_j,\alpha_j,\beta_j)
=(\gamma_i,\alpha_i,\beta_i), 1\leq j\leq n-1\}|$.
Then the following statements hold.\\
 (i) If $\ell(i) \ge 2$ then there is an eigenvalue
 $\theta$ of $L_1$ with
 $$
 \alpha_i + 2 \sqrt{\beta_i \gamma_i} \cos\left( \frac{2 \pi }{\ell(i)+1}
 \right) \leq \theta < k.
 $$\\
 (ii) If $\ell(i) \ge 3$ then there is an eigenvalue
 $\theta$ of $L_1$ with
 $$
 \alpha_i + 2 \sqrt{\beta_i \gamma_i} \cos\left( \frac{j \pi }{\ell(i)+1}
 \right) \le \theta
 \le \alpha_i + 2 \sqrt{\beta_i \gamma_i} \cos\left( \frac{(j-2)
 \pi}{\ell(i)+1} \right),
 $$
 for all $j=3,\ldots, \ell(i)$.
\end{lemma}

\subsection{Distance-Regular Graphs}\label{drg}

We now review some basic definitions and results concerning
distance-regular graphs.

For $\Gamma$ a finite, connected graph, denote by $d(x,y)$ the
path-length distance
between any two vertices $x,y$ in the vertex set $V(\Gamma)$ of
$\Gamma$ (i.e. the length of a shortest path), and by
$D=D_{\Gamma}$ the diameter of $\Gamma$ (i.e. the maximum distance
between any two vertices of $\Gamma$). For any $y \in V(\Gamma)$,
let $\Gamma_i(y)$ be the set of vertices in $\Gamma$ at distance
precisely $i$ from $y$, where $i \in \N_0$ is a non-negative
integer not exceeding $D$. In addition, define $\Gamma_{-1}(y) = \Gamma_{D+1}(y) :=
\emptyset$.

Following \cite[p.126]{bcn}, a
finite, connected graph $\Gamma$ is called a {\em distance-regular
graph} if there are integers $b_i,c_i$, $i=0,1,\ldots,D$, such
that, for any two vertices $x,y \in V(\Gamma)$ at distance $i =
d(x,y)$, there are precisely $c_i$ neighbors of $y$ in
$\Gamma_{i-1}(x)$ and $b_i$ neighbors of $y$ in $\Gamma_{i+1}(x)$.
In particular, $\Gamma$ is regular with valency $k := b_0$. The
numbers $c_i, b_i$ and
$$
a_i := k - b_i -c_i \, \, \, \, \, \, (0 \le i \le  D)
$$
(i.e. the number of neighbors of $y$ in $\Gamma_i(x)$ for
$d(x,y)=i$) are called the {\em intersection numbers} of $\Gamma$.
Note that $b_D = c_0 =a_0 := 0$ and $c_1 = 1$. In addition, we
define $k_i:=|\Gamma_i(y)|$ for any vertex $y\in V(\Gamma)$,
$i=0,1,\ldots,D$. This definition for distance-regular
graphs is easily seen to be equivalent to the one given in the introduction.

For $\Gamma$ a distance-regular graph as above, we define
\begin{equation}\label{induced tridiagonal seq}
\T_{\Gamma}
:= \Big((c_i,a_i,b_i) \Big) _{i=1}^{D}
\end{equation}
and we let
\begin{equation}\label{induced graphical seq}
\g_{\Gamma} := \Big((\gamma_i,\alpha_i,\beta_i)\Big)_{i=1}^{g+1}
\end{equation}
denote the (necessarily unique)
maximal length subsequence
of $\T_{\Gamma}$ for which the $i$\,th term of
$\g_{\Gamma}$ is not equal to
the $(i+1)$\,th term of $\g_{\Gamma}$ for
all $1 \le i \le D-1$.
In addition, we define the numbers
\begin{eqnarray}
\HH = \HH_{\Gamma} &:=& |\{j \,:\, (c_j,a_j,b_j)
=(c_1,a_1,b_1),\,1\le j \le D-1\}|\,, \mbox{ and } \label{head}\\
\TT=\TT_{\Gamma}&:=&|\{j \,:\, (c_j,a_j,b_j)=(b_1,a_1,c_1),\,\HH < j \le D-1\}|\label{tail}
\end{eqnarray}
which are called the {\em head} and the {\em tail}
of $\Gamma$, respectively.
Note that by
\cite[Lemma~2.1]{two-thm}, it follows that tail $\TT$ satisfies the following :
\begin{equation}\label{drg-h>t}
\TT\leq \HH \mbox{ and, } \mbox{if~} \TT\ge 1  \mbox{~then~}
(c_{D-\TT},a_{D-\TT},b_{D-\TT})=\cdots=(c_{D-1},a_{D-1},b_{D-1})=(b_1,a_1,1).
\end{equation}

\subsubsection{Intersection Numbers}

For the rest of Section \ref{pre}, we suppose that $\Gamma$ is a
distance-regular graph with valency $k\geq 3$, diameter $D\geq 2$,
intersection numbers $a_i,b_i,c_i$, $0 \le i \le D$ and
$\g_{\Gamma}=\Big((\gamma_i,\alpha_i,\beta_i)\Big)_{i=1}^{g+1}$.

In \cite[Proposition 4.1.6]{bcn} and \cite[Lemma 2.1 (ii)]{two-thm},
it is shown that the following inequalities always hold :
\begin{equation}\label{inequals}
k=b_0 > b_1 \ge b_2 \ge \cdots \ge b_{D-1}>b_D =0 \mbox{ and } 1 = c_1 \le c_2 \le \cdots \le c_D \le k,
\end{equation}
\begin{equation}\label{lambda}
a_i\ge a_1+1-\min\{b_i,c_i\} ~~(1 \le i \le D-1).
\end{equation}
In particular, it follows that for every term
$(\gamma_i,\alpha_i,\beta_i)$ in $\g_{\Gamma}$, $\beta_{i}\geq \beta_{i+1}$ and
$\gamma_i\leq \gamma_{i+1}$ hold. For each $1 \le i \le g$, define
\begin{eqnarray}
 s(i) = s_{\Gamma}(i)&:=& \min \{j\,:\, (c_j,a_j,b_j)=
(\gamma_i,\alpha_i,\beta_i)\,,\,1\le j\le D-1\},
\label{f(i)}\\
\ell(i)=\ell_{\Gamma}(i)&:=& |\{j \,:\, (c_j,a_j,b_j)=
(\gamma_i,\alpha_i,\beta_i),\,1\le j \le D-1\}|,\label{l(i)}
\end{eqnarray}
and define $s(g+1)=D$. Note that $s(1)=1$ $\ell(1) = \HH_{\Gamma}$,
$\ell(g+1)=1$, and that $s(i+1)-s(i)=\ell(i)$ holds for all $1 \le i
\le g$.

\subsubsection{Diameter Bounds}

The following result is originally due to
A.~A.~Ivanov \cite{ivanov} (cf. \cite[Theorem 5.9.8]{bcn}). Note that
$\N$ denotes the set of positive integers.

\begin{theorem}(A.~A.~Ivanov's Diameter Bound)\label{ivanov} {\ \\}
Let $k \geq 3$ be an integer. Then there is a
function $F: \N  \rightarrow \N $ so that, for all distance-regular
graphs $\Gamma$ with valency $k$, diameter $D_{\Gamma}$,
and head  $\mbox{\em \HH}_{\Gamma}$, the inequality
\[D_{\Gamma} \leq F(k) ~\mbox{\em \HH}_{\Gamma} \] holds.
\end{theorem}

Note that it was also shown in \cite{ivanov} (cf. \cite[Theorem 5.9.8]{bcn}) that
one can in fact take $F(k) =4^k$ in the last theorem.\\

Now, in order to show that there are only finitely many of
distance-regular graphs $\Gamma$ with fixed valency $k\geq 3$, it
suffices to show that the diameter $D_{\Gamma}$ of any such graph
is bounded above by some function $f:\N \to \N$ depending only on
$k$, since $|V(\Gamma)|\leq1+\sum_{i=1}^{D_{\Gamma}}k(k-1)^{i-1}$.
Thus, in view of Theorem~\ref{ivanov}, it also suffices to show
that the head $\HH_{\Gamma}$ is bounded above by some function $g$
in $k$. In particular, the following result also holds (as we can
take $g(k)$ to be a constant function).

\begin{cor}\label{ivanov-conclusion}
Suppose that $k\geq 3$ and $C\geq 1$ are positive integers.
Then there are only finitely many distance-regular
graphs $\Gamma$ with valency $k$ and head $\mbox{\em \HH}_{\Gamma}\leq C$.
\end{cor}

\subsubsection{Eigenvalues of Distance-Regular Graphs}

The tridiagonal matrix $L_1 = L_1(\Gamma)$ associated to $\Gamma$
is defined by
$$
L_1:= \left\lgroup
\begin{tabular}{llllllll}
$0$ & $b_0$\\
$c_1$ & $a_1$ & $b_1$\\
& $c_2$ & $a_2$ & $b_2$\\
& & . & . & .\\
& & & $c_i$ & $a_i$ & $b_i$\\
& & & & . & . & .\\
& & & & & $c_{D-1}$ & $a_{D-1}$ & $b_{D-1}$\\
& & & & &\makebox{\hspace{.324cm}} & $c_D$ & $a_D$
\end{tabular}
\right\rgroup,
$$
and $\theta \in {\mathbb R}$ is an eigenvalue of
$\Gamma$ if $\theta$ is an eigenvalue of $L_1(\Gamma)$
(\cite[p.129]{bcn}). Note that any distance-regular graph $\Gamma$
with diameter $D=D_{\Gamma}$ has exactly $D+1$ distinct
eigenvalues (\cite[p.128]{bcn}). Moreover, if $\theta$ is
an eigenvalue of $\Gamma$, then
$(u_0, u_1, \ldots,u_D)^T$ is called the standard sequence
of $\Gamma$ associated with $\theta$, which
is a right eigenvector of $L_1(\Gamma)$ associated
with eigenvalue $\theta$, and the multiplicity $m(\theta)$ of
$\theta$ is given by
\begin{equation}\label{biggs}
m(\theta)=\frac{|V(\Gamma)|}{\sum_{i=0}^{D}k_iu_i^2(\theta)}.
\end{equation}
This equation is known as {\em Biggs' formula} (\cite[Theorem 21.4]{biggs}).
Note that in view of
Equations~(\ref{cnumber}) and (\ref{crn}) it follows by this last formula
that the multiplicity of eigenvalue $\theta_i$ of $\Gamma$
is equal to the Christoffel number $m_i$ of $L_1(\Gamma)$.

\forme{ It is shown in \cite{ban-ito-II, ban-ito-89} that the conjecture holds for $k = 3, 4$ by proving that there are only finitely many triangle-free
distance-regular graphs for fixed $k$ and $D-\HH-\TT$. In \cite{suzuki}, this was generalized by Suzuki to include
distance-regular graphs that are not necessarily triangle-free, but that satisfy $k\geq (a_1+1)(a_1+2)$. In Theorem \ref{new-two-thm-1}, we proved that this latter condition is unnecessary. \\
The proof of the Bannai-Ito conjecture for $k = 5, 6,7$ in \cite{val10} employs a generalization of Theorem \ref{new-two-thm-1}, whose proof uses a combinatorial argument that is based on the Terwilliger Tree bound \cite{terw-tree}.
In \cite{val10}, an even more general result is shown to hold for triangle-free distance-regular graphs, whose proof is based on algebraic arguments.
Here, Theorem \ref{new-two-thm} extend this result to arbitrary distance-regular graphs.}

\section{Graphical Sequences}\label{graphical seq}

In this section, we define graphical sequences and tridiagonal
sequences. Note that these are similar (but not identical) to the
ones presented in \cite{rnp}. The definition for these sequences
is motivated by the sequences ${\mathcal G}_{\Gamma}$ and
$\T_{\Gamma}$ associated to $\Gamma$ a distance-regular graph that were
presented in the last section.

For integers $\kappa \ge 3$ and $\lambda \ge 0$ with $\lambda
\le\kappa-2$, define
\[
V_{\kappa,\lambda}:= \{(\gamma,\alpha,\beta)\in \mathbb{N}_0^3
\,:\, \beta, \gamma \ge 1,~\gamma+\alpha+\beta= \kappa \mbox{~and~}
\alpha \ge \max\{\lambda+1-\beta,\lambda+1-\gamma\} \}.
\]

\begin{defi}\label{GS}
With $\kappa$, $\lambda$ and $V_{\kappa,\lambda}$ as just defined above,
a sequence $\mathcal{G}=\Big(
(\gamma_i,\alpha_i,\beta_i)\Big)_{i=1}^{g+1}$ of distinct terms in
$\mathbb{N}_0^3$ is called a {\em $(\kappa,\lambda)$-graphical
sequence} if it satisfies the following conditions:
\begin{itemize}
\item[{\em (G0)}] $(\gamma_i,\alpha_i,\beta_i)\in V_{\kappa,\lambda}~(1\leq i\leq g)$,
\item[{\em (G1)}]
$(\gamma_1,\alpha_1,\beta_1)=(1,\lambda,\kappa-\lambda-1)$,
\item[{\em (G2)}] $\beta_i\ge \beta_{i+1}~(1 \le  i \le g-1)$ and $\gamma_i\le
\gamma_{i+1}~(1 \le  i \le g)$,
\item[{\em (G3)}] $\beta_{g+1}=0$ and $\gamma_{g+1}+\alpha_{g+1}=\kappa$.
\end{itemize}
\end{defi}

Let $\g=\Big( (\gamma_i,\alpha_i,\beta_i) \Big)_{i=1}^{g+1}$
be a $(\kappa,\lambda)$-graphical sequence and let
$\ell:\{1,\ldots, g+1\}\rightarrow \mathbb{N}$ be a function with
$\ell(g+1)=1$. For each $1\leq i\leq g+1$, define $s_{\ell}(i) =s(i)$ by
\begin{eqnarray}
&& s(1):=1,\nonumber\\
&& s(i):=1+\sum_{j=1}^{i-1}\ell(j)~~(2\leq i\leq g+1).\label{s(i)}
\end{eqnarray}

\begin{defi}\label{TDS}
With $\g$, $\ell$ and $s$ as just defined above,
the sequence of triples
$\T=\mathcal{T}(\g,\ell):=\Big((c_m,a_m,b_m)\Big)_{m=1}^{s(g+1)}$
given by putting, for each $1\leq i\leq g+1$,
$$
(c_{s(i)+j},a_{s(i)+j},b_{s(i)+j})=
(\gamma_i,\alpha_i,\beta_i)~~(0\leq j\leq \ell(i)-1)
$$
is called the {\em $(\kappa,\lambda)$-tridiagonal sequence
(associated with $\g$ and $\ell$)}.
\end{defi}

Given $\T=\mathcal{T}(\g,\ell)$ as
in this last definition, we define
the {\em head} $\HH=\HH_{\mathcal{T}}$, the
{\em tail} $\TT=\TT_{\T}$ and the {\em diameter}
$D=D_{\mathcal{T}}$ of $\T$ to be
\begin{eqnarray}
\HH_{\mathcal{T}}&:=&\ell(1),\label{def-h}\\
\TT_{\mathcal{T}}&:=& \left|\{\HH_{\T}<i\leq s(g+1)\,:\,(c_i,a_i,b_i)=(\kappa-\lambda-1,\lambda,1)\}\right|,
\mbox{ and } \label{def-t}\\
D_{\mathcal{T}}&:=&s(g+1), \label{D(T)}
\end{eqnarray}
respectively.
Note that $\HH$ and $\TT$ satisfy
\[
\HH \ge \TT \mbox{ and, if } \TT\ge 1 \mbox{~then~}
(c_{D-\TT},a_{D-\TT},b_{D-\TT})=\cdots=(c_{D-1},a_{D-1},b_{D-1})=(b_1,a_1,1)
\] (see (\ref{drg-h>t})).

Note that if $\Gamma$ is a distance-regular graph,
with diameter $D_{\Gamma}$ and
$\T_{\Gamma} = \Big( (c_m,a_m,b_m) \Big)_{m=1}^{D_{\Gamma}}$,
then, referring to (\ref{induced graphical seq}) and (\ref{l(i)}),
it follows that the sequence $\g_{\Gamma}$ is a
$(b_0,a_1)$-graphical sequence and that
$\T_{\Gamma}$ is the $(b_0,a_1)$-tridiagonal sequence
$\T(\g_{\Gamma},\ell_\Gamma)$.

Now, given a
$(\kappa,\lambda)$-graphical sequence
$\g$, function $\ell$ and the
$(\kappa,\lambda)$-tridiagonal sequence
$\T=\T(\g,\ell)=\Big((c_m,a_m,b_m)\Big)_{m=1}^{D}$ as
in Definitions~\ref{GS} and \ref{TDS},
we define the tridiagonal matrix $L_1(\T)$ associated
to $\T$ by
$$
 L_1(\T):= \left\lgroup
 \begin{tabular}{llllllll}
 $0$ & $\kappa$\\
 $c_1$ & $a_1$ & $b_1$\\
 & $c_2$ & $a_2$ & $b_2$\\
 & & . & . & .\\
 & & & $c_i$ & $a_i$ & $b_i$\\
 & & & & . & . & .\\
 & & & & & $c_{D-1}$ & $a_{D-1}$ & $b_{D-1}$\\
& & & & &\makebox{\hspace{.324cm}} &
 $c_{D}$ & $a_{D}$
 \end{tabular}
 \right\rgroup.
$$
It follows by the results in Section~\ref{cn}, that the tridiagonal
matrix $L_1(\T)$ has exactly $D+1$ distinct eigenvalues,
$\kappa=\theta_0
>\theta_1>\theta_2>\cdots >\theta_D$, say, which we call {\em the
eigenvalues of $\T$} and denote by
\begin{equation}\label{def-Eg}
\mathcal{E}_{\T}:=\{\theta_i \,:\, 0\leq i\leq D\}.
\end{equation}
Note that applying formulae (\ref{v_i}) and (\ref{def-k_i}) to the
matrix $L_1(\T)$ we obtain, for each $1\leq i\leq g+1$,
\begin{equation} \label{kappa_i}
\kappa_{s(i)+j}=\frac{\kappa}{\beta_i}\left(\frac{\beta_i}{\gamma_i}\right)^{j}
\prod_{t=1}^{i-1}\left( \frac{\beta_t}{\gamma_t} \right)
^{\ell(t)}=\frac{\kappa\, b_1\cdots b_{s(i)+j-1} }{c_1c_2\cdots
c_{s(i)+j}}~~~(j=0,\ldots, \ell(i)-1).
\end{equation}
We define the {\em Christoffel numbers} of
$\T$ to be the Christoffel
numbers associated with $L_1(\T)$ (cf. Section~\ref{cn}).

Now, in case $L_1 = L_1(\T_{\Gamma})$ for a distance-regular graph
$\Gamma$ then, for any $\theta, \theta'$ distinct algebraic
conjugate eigenvalues of $\Gamma$, the multiplicities of $\theta$
and $\theta'$ are equal (\cite[Proposition III.1.5]{banito}).
Hence so are the corresponding Christoffel numbers,
which implies that $\sum_{i=0}^{D}\kappa_i u_i^2(\theta)=\sum_{i=0}^{D}\kappa_i u_i^2(\theta')$ holds.\\
Motivated by this fact, we will be interested in
$(\kappa,\lambda)$-tridiagonal sequences $\T$ that satisfy the
following key property:
\begin{itemize}
\item[(AC)] Any two eigenvalues of $\T$ which are
algebraically conjugate (over $\mathbb{Q}$)
have the same Christoffel numbers.
\end{itemize}

We conclude this section with a useful result concerning graphical
sequences. Suppose that
$\g=\Big((\gamma_i,\alpha_i,\beta_i)\Big)_{i=1}^{g+1}$
is a $(\kappa,\lambda)$-graphical sequence for some integers $\kappa
\geq 3$ and $0\leq \lambda\leq \kappa-2$, as in Definition~\ref{GS}.
For each $1 \le i \le g$
we define the {\em i\,th right} and {\em i\,th left guide point} by
\begin{equation}\label{R-L}
\frak{R}_i=\frak{R}_i(\mathcal{G}):=\alpha_i+2\sqrt{\beta_i\gamma_i}\mbox{~~and~~}
\frak{L}_i=\frak{L}_i(\mathcal{G}):=\alpha_i-2\sqrt{\beta_i\gamma_i}~~
(1\leq i\leq g)
\end{equation}
respectively.
In addition, we put
$\frak{R}_{\max}=\frak{R}_{\max}(\g) :=
\max\{\frak{R}_i \,:\, 1\leq i\leq g\}$.\\
Moreover, for each $1 \le i \le g$, we define the {\em i\,th guide interval} to be the open interval
\begin{equation}\label{def-guide interval}
I_i =I_i(\mathcal{G}) :=(\frak{L}_i,\frak{R}_i).
\end{equation}

The following lemma is a slight extension of Lemma 3.1 in \cite{rnp}. We provide a proof of it for the sake of
completeness. Note that a sequence
$r_1,\ldots,r_n$ of real numbers is called {\em unimodal} if there
exists some $1 \le t \le n$ satisfying $r_1\le r_2\le \cdots \le r_t$
and $r_t\ge r_{t+1}\ge \cdots \ge r_n$.

\begin{lemma}\label{new-unimodular}
Suppose that  $\kappa \geq 3$ and $\lambda \geq 0$ are integers
with $\lambda \le \kappa-2$, and that
$\mathcal{G}=\Big((\gamma_i,\alpha_i,\beta_i)\Big)_{i=1}^{g+1}$
is a $(\kappa,\lambda)$-graphical sequence. Then the following
hold.
\begin{itemize}
\item[(i)] The inequality $\frak{R}_i\geq \frak{R}_1$ holds for
all $1 \le i \le g$, with equality holding if and only if
$(\gamma_i,\alpha_i,\beta_i) \in \{(1,\lambda, \kappa-\lambda-1),
(\kappa-\lambda-1,\lambda,1) \}$. \item[(ii)] For any $2 \le i
\leq g$, if $\beta_{i}\ge\gamma_{i}$ then $\frak{R}_{i-1} <
\frak{R}_{i}$. \item[(iii)] For any $2 \leq i \leq g-1$, if
$\beta_{i}\le\gamma_{i}$ then $\frak{R}_{i+1}< \frak{R}_{i}$.
\end{itemize}
In particular, by (ii) and (iii), it follows that  the sequence
$(\frak{R}_i)_{i=1}^{g}$ is unimodal.
\end{lemma}

\pf First note that by (G0) and (G1) in Definition~\ref{GS},
for each $1\leq i,j\leq g$, we
have
\begin{eqnarray}
\frak{R}_i-\frak{R}_j&=&
(\sqrt{\beta_j}-\sqrt{\gamma_j})^2-(\sqrt{\beta_i}-\sqrt{\gamma_i})^2,
\mbox{ and } \label{fr}\\
\beta_1&\geq &\gamma_i.\label{beta_1 >gamma_i}
\end{eqnarray}
Now, to see that (i) holds, note that by (G0),(G2) and
(\ref{beta_1 >gamma_i}), $\sqrt{\beta_1}-1\geq
\left|\sqrt{\beta_i}-\sqrt{\gamma_i} \right|$ holds. Hence
$\frak{R}_i\geq \frak{R}_1$ holds in view of (\ref{fr}) with
$j=1$. Moreover, equality holds if and only if
$(\gamma_i,\alpha_i,\beta_i)=(\gamma_1,\alpha_1,\beta_1)$ if
$\beta_i\geq \gamma_i$ and
$(\gamma_i,\alpha_i,\beta_i)=(\beta_1,\alpha_1,\gamma_1)$ if
$\gamma_i\geq \beta_i$.\\
\noindent To complete the proof of the lemma, note that (ii) and
(iii) follow from (\ref{fr}) and (G2), since $\beta_{i-1}\geq
\beta_i\geq \gamma_i\geq \gamma_{i-1}$ and $\gamma_{i+1}\geq
\gamma_i\geq \beta_i\geq \beta_{i+1}$ hold for (ii) and (iii),
respectively. \epf

\section{A Key Result}\label{central}

In this section we will state without proof
a key result (Theorem~\ref{b-i-small-thm}) that we will then
use to prove the main result of this paper (Theorem~\ref{b-i-thm}).
We will then give a sketch a proof of
this key result which we will prove
in Sections~\ref{no-thy-section} to \ref{distri},
inclusive.

For $w = (w_i)_{i=1}^n$ any sequence, we put
\begin{equation}\label{underline}
\underline{w} := \{w_i : 1 \le i \le n\}\,,
\end{equation}
i.e. the set consisting of all distinct terms in $W$.
To state Theorem~\ref{b-i-small-thm}, we
will require the following key definition:

\begin{defi}\label{RT-def}
Let $\kappa\geq 3$ and $\lambda\geq 0$ be integers with
$\lambda\leq \kappa-2$. A  {\em $(\kappa,\lambda)$-quadruple}
is a quadruple $(\g, \Delta; L,\ell)$ such that\\
(i) $\g:=\Big(\delta_i := (\gamma_i,\alpha_i,\beta_i)\Big)_{i=1}^{g+1}$ is a
$(\kappa,\lambda)$-graphical sequence (cf. Definition~\ref{GS}),\\
(ii) $\Delta = ( \delta_{i_p} )_{p=1}^\tau$
is a subsequence of $\g$ in which
$(1,\lambda,\kappa-\lambda-1) \in \underline{\Delta}$ (i.e., $i_1=1$) and
$(\gamma_{g+1},\alpha_{g+1},\beta_{g+1}) \not\in \underline{\Delta}$, and\\
(iii) $\ell: \{1,\dots,g+1\} \rightarrow \mathbb{N}$ and
$L: \{1,\dots,g+1\} \setminus \{i_1,\dots,i_{\tau}\} \rightarrow \mathbb{N}$
are functions with $\ell(g+1)=1$ and $L(i) = \ell(i)$ for all
$i \in \{1,\dots,g+1\} \setminus \{i_1,\dots,i_{\tau}\}$.\\
\end{defi}

\begin{theorem}\label{b-i-small-thm}
Let $\kappa\geq 3$ and $\lambda\geq 0$ be integers with
$\lambda\leq \kappa-2$, and let
$\g=\Big(\delta_i:=(\gamma_i,\alpha_i,\beta_i)\Big) _{i=1}^{g+1}$
be a $(\kappa,\lambda)$-graphical sequence. Suppose that $\Delta=(\delta_{i_p})_{p=1}^{\tau}$ is a subsequence of $\g$
with $(1,\lambda,\kappa-\lambda-1) \in \underline{\Delta}$ and
$(\gamma_{g+1},\alpha_{g+1},\beta_{g+1}) \not\in \underline{\Delta}$, and $L: \{1,\dots,g+1\}
\setminus \{i_1,\dots,i_{\tau}\} \rightarrow \mathbb{N}$ is a function. Suppose $\epsilon>0$ is a real
number, $C:=C(\kappa)>0$ is a constant, and
$\ell:\{1,\ldots, g+1\}\rightarrow \mathbb{N}$ is any function for which $(\g, \Delta; L,\ell)$ is a
$(\kappa,\lambda)$-quadruple and the associated $(\kappa,\lambda)$-tridiagonal
sequence $\T=\T(\g,\ell)$ satisfies\\
(i) Property (AC), \\
(ii) $D_{\T}\leq C\mbox{\em \HH}_{\T}$, and \\
(iii) $D_{\T}-(\mbox{\em\HH}_{\T}+\mbox{\em\TT}_{\T})>\epsilon
\mbox{\em\HH}_{\T}$,\\
where $\mbox{\em\HH}_{\T}$,
$\mbox{\em\TT}_{\T}$ and $D_{\T}$ are as defined in (\ref{def-h})--(\ref{D(T)}), respectively.\\
Then, there exist positive constants
$F:=F(\kappa,\g,\Delta,L)$ and
$H:=H(\kappa,\lambda,\epsilon,\g,\Delta,L)$
such that if $\ell (i_p)>F$ holds for all $1\leq p\leq \tau$, then $\mbox{\em \HH}_{\T} \leq H$ holds.
\end{theorem}

We will now use Theorem \ref{b-i-small-thm} to prove Theorem \ref{b-i-thm}, the main
theorem of this paper. To do this, we will make
use of the following result:

\begin{prop}\label{prop-main-thm}
Let $\kappa\geq 3$ and $\lambda\geq 0$ be integers with
$\lambda\leq \kappa-2$, and let
$\g=\Big(\delta_ i: = (\gamma_i,\alpha_i,\beta_i)\Big) _{i=1}^{g+1}$
be a $(\kappa,\lambda)$-graphical sequence.
Suppose $\epsilon>0$ is a real
number, $C:=C(\kappa)>0$ is a constant, and
$\ell:\{1,\ldots, g+1\}\rightarrow \mathbb{N}$ is any function with
$\ell(g+1)=1$, such that the associated $(\kappa,\lambda)$-tridiagonal
sequence $\T=\T(\g,\ell)$ satisfies\\
(i) Property (AC),\\
(ii) $D_{\T}\leq C\,\mbox{\em \HH}_{\T}$, and\\
(iii) $D_{\T}-(\mbox{\em \HH}_{\T}+\mbox{\em \TT}_{\T})>\epsilon\,
\mbox{\em \HH}_{\T}$.\\
Then there
exists a positive constant $H:=H(\kappa,\lambda,\epsilon,\g)$
such that $\mbox{\em \HH}_{\T} \leq H$ holds.
\end{prop}
\pf Suppose that $\kappa$, $\lambda$, $\g$, $\epsilon$, $C$, $\ell$
are as in the statement of the proposition.
First, we show that the following statement holds:
\begin{itemize}
\item[($\ddag$)] For each $i=0,\ldots,g-1$, there exists a
subsequence $\g_i$ of $\g$ with precisely $(i+1)$-terms satisfying
$(1,\lambda,\kappa-\lambda-1)\not\in \underline{\g_i}$ and
$(\gamma_{g+1},\alpha_{g+1},\beta_{g+1})\in \underline{\g_i}$ for
which there is a positive constant
$\mathbb{L}_i:=\mathbb{L}_i(\kappa,\lambda,\epsilon,\g,\g_i)$ such
that \[\ell(j)\leq \mathbb{L}_i\] holds for all
$(\gamma_j,\alpha_j,\beta_j)\in \underline{\g_i}$.
\end{itemize}
{\noindent{\em Proof of $(\ddag)$:}} We use induction on
$i$. In case $i=0$, ($\ddag$) holds for
the subsequence $\g_0 :=\Big( (\gamma_{g+1},\alpha_{g+1},\beta_{g+1})\Big)$
and constant $\mathbb{L}_0:=1$.

So, assume that ($\ddag$) holds for all $i=s$, with $0 \le s\leq
g-2$, i.e. there is a subsequence $\g_s=\Big((\gamma_{i_p},\alpha_{i_p},\beta_{i_p})\Big)_{p=1}^{s+1}$
of $\g$ with $(1,\lambda,\kappa-\lambda-1)\not\in
\underline{\g_s}$ and $(\gamma_{g+1},\alpha_{g+1},\beta_{g+1})\in
\underline{\g_s}$ for which there is a positive constant
$\mathbb{L}_s:=\mathbb{L}_s(\kappa,\lambda,\epsilon,\g,\g_s)>0$
such that $\ell(i_p) \le \mathbb{L}_s$ holds for all $1 \le p \le s+1$.

Let $\mathcal{L}(\{i_1,\dots,i_{s+1}\})$ denote the set consisting
of those functions $L:\{i_1,\ldots,i_{s+1}\}\rightarrow
\mathbb{N}$ satisfying $L(i_p)\leq \mathbb{L}_s$ for all $1\leq
p\leq s+1$. Note that the set $\mathcal{L}(\{i_1,\dots,i_{s+1}\})$
depends only on $\kappa,\lambda,\epsilon,\g$ and $\g_s$.
Let $\Delta_s$ denote the subsequence of $\g$
obtained by removing the terms in $\g_s$ from $\g$. Put
$m=m(\g,\g_s):=\min\{2\leq i\leq
g\,:\,(\gamma_i,\alpha_i,\beta_i)\in \underline{\Delta_s}\}$.

Define positive constants
$\widetilde{F}_s=\widetilde{F}_s(\kappa,\lambda,\epsilon,\g,\g_s)$
and $\widetilde{H}_s=\widetilde{H}_s(\kappa,\lambda,\epsilon,\g,\g_s)$ by
\begin{eqnarray*}
\widetilde{F}_s&:=&\max\{F(\kappa,\g,\Delta_s,L)\,:\,
L\in \mathcal{L}(\{i_1,\dots,i_{s+1}\})\},\\
\widetilde{H}_s&:=&\max\{H(\kappa,\lambda,\epsilon,\g,\Delta_s,L)\,:\,
L\in \mathcal{L}(\{i_1,\dots,i_{s+1}\})\},
\end{eqnarray*}
where $F(\kappa,\g,\Delta_s,L)$ and
$H(\kappa,\lambda,\epsilon,\g,\Delta_s,L)$
are the constants given by applying
Theorem~\ref{b-i-small-thm} to the $(\kappa,\lambda)$-quadruple
$(\g,\Delta_s;L,\ell)$.

Then, by Theorem~\ref{b-i-small-thm}, either (a)
$\ell(i)>\widetilde{F}_s$ holds for all
$(\gamma_i,\alpha_i,\beta_i)\in \underline{\Delta_s}$, in which
case we can let $\g_{s+1}$ be the sequence defined by adding the
term $(\gamma_m,\alpha_m,\beta_m)$ to the beginning of $\g_s$ and
put $\mathbb{L}_{s+1}
:=\max\{\mathbb{L}_s,C(\kappa)\,\widetilde{H}_s\}$, or (b) there
exists $(\gamma_n,\alpha_n,\beta_n)\in \underline{\Delta_s}$ such
that $\ell(n)\leq \widetilde{F}_s$ holds, in which case we can let
$\g_{s+1}$ be the sequence defined by inserting the term
$(\gamma_j,\alpha_j,\beta_j)$ with $j:=\max\{m,n\}$ into the
sequence $\g_s$ (according to its place in $\g$) and put
$\mathbb{L}_{s+1}:=\max\{\mathbb{L}_s,\,C(\kappa)\,\widetilde{F}_s\}$.
This completes the proof that statement ($\ddag$) holds. \epf

To complete the proof of the proposition, we apply ($\ddag$) for
$i=g-1$. In particular, for this choice of $i$,
$\g_{g-1}=\Big((\gamma_i,\alpha_i,\beta_i)\Big)_{i=2}^{g+1}$, and constant $\mathbb{L}_{g-1}$ depends only on
$\kappa,\lambda,\epsilon$ and $\g$, and hence the set
$\mathcal{L}(\{i_1,\dots,i_{g}\}) =\mathcal{L}(\{2,\dots,g+1\})$
of those functions $L:\{i_1,\ldots,i_{g}\}\rightarrow
\mathbb{N}$ satisfying $L(i_p)\leq \mathbb{L}_{g-1}$ for all
$1\leq p\leq g$ depends only on
$\kappa,\lambda,\epsilon$ and $\g$. Since
$(\g,\Delta;L,\ell)$ is a $(\kappa,\lambda)$-quadruple for the
subsequence $\Delta = \Big((\gamma_1,\alpha_1,\beta_1)=
(1,\lambda,\kappa-\lambda-1)\Big)$ of $\g$ and $L$ any function in
$\mathcal{L}(\{2,\dots,g+1\})$, it follows by applying
Theorem~\ref{b-i-small-thm} to $(\g,\Delta;L,\ell)$ that there exists a constant
\[
C=C(\kappa,\lambda,\epsilon,\g):=
\max\{F(\kappa,\g,\Delta,L),H(\kappa,\lambda,\epsilon,\g,\Delta,L)\,:\,
L\in \mathcal{L}(\{2,\dots,g+1\})\},
\]
where $F(\kappa,\g,\Delta,L),H(\kappa,\lambda,\epsilon,\g,\Delta,L)$
are the constants given by
Theorem~\ref{b-i-small-thm} so that
$$
\HH_{\T}\leq C$$ holds. This completes the proof of the
proposition.
\epf

In order to prove Theorem~\ref{b-i-thm}, we will also make use of the following
result from \cite{two-thm}, which generalizes results of
Bannai and Ito \cite{ban-ito-II, ban-ito-89} and Suzuki \cite{suzuki}:

\begin{theorem}(\cite[Theorem 1.2]{two-thm})\label{new-two-thm}{\ \\}
Suppose that $k\geq 3$ is a fixed integer. Then there exists a
positive number $\epsilon_0= \epsilon_0(k)$, depending only on
$k$, so that there are only finitely many distance-regular graphs
with valency $k$, head $\mbox{\em \HH}_{\Gamma}$, tail $\mbox{\em
\TT}_{\Gamma}$, and diameter $D_{\Gamma}$ that satisfy
\[D_{\Gamma}-(\mbox{\em \HH}_{\Gamma} + \mbox{\em \TT}_{\Gamma}) \leq  \epsilon_0 \,\mbox{\em \HH}_{\Gamma}.\]
\end{theorem}
\forme{Owing to this result, to prove Theorem~\ref{b-i-thm}, it is
enough to show that for given valency $k\geq 3$, there are only
finitely many distance-regular graphs $\Gamma$ with valency $k$,
head $\HH_{\Gamma}$, tail $\TT_{\Gamma}$ and diameter $D_{\Gamma}$
that satisfy
\[D_{\Gamma}-(\HH_{\Gamma}+\TT_{\Gamma})> \epsilon_0 \HH_{\Gamma}.\]}

\noindent{\em Proof of Theorem \ref{b-i-thm}: }
Let $k \geq 3$ be a fixed integer.
By Theorem~\ref{new-two-thm},
there exists a constant $\epsilon_0=\epsilon_0(k)>0$ (which
depends only on $k$) such that there
are only finitely many distance-regular graphs $\Gamma$ with
valency $k$, head $\HH_{\Gamma}$, tail $\TT_{\Gamma}$ and diameter
$D_{\Gamma}$ that satisfy
\[ D_{\Gamma}-(\HH_{\Gamma}+\TT_{\Gamma})\leq \epsilon_0 \HH_{\Gamma}.\]

Now, suppose that $\Gamma$ is any distance-regular graph with
valency $k$ that satisfies
\begin{equation}\label{middle-eq}
D_{\Gamma}-(\HH_{\Gamma}+\TT_{\Gamma})> \epsilon_0 \HH_{\Gamma}.
\end{equation}
Then, by Theorem~\ref{ivanov} and (\ref{middle-eq}),
the $(k,a_1)$-tridiagonal sequence
$\T_{\Gamma} = \T(\g_{\Gamma},\ell_{\Gamma})$
(cf. (\ref{induced tridiagonal seq}))
satisfies all of conditions (i)--(iii) in Proposition~\ref{prop-main-thm},
where $a_1$
is an intersection number of $\Gamma$.

Therefore, for any distance-regular graph $\Gamma$ with
valency $k$ that satisfies (\ref{middle-eq}), it follows that
$$
\HH_{\Gamma}\leq C(k):=\max\{H(k,a_1,\epsilon_0(k),\g)
: 0 \le a_1 \le k-2, \mbox{ $\g$ is a
$(k,a_1)$-graphical sequence} \}
$$
where $H(k, a_1,\epsilon_0(k),\g)$ is the constant
given by Proposition~\ref{prop-main-thm}
(note that in the formula for $C(k)$,
taking a maximum is appropriate since the number of
integers $a_1$ with $0\leq a_1 \leq k-2$ is finite,
and so is the number of $(k,a_1)$-graphical sequences).
Theorem~\ref{b-i-thm} now follows by applying
Corollary~\ref{ivanov-conclusion} with the constant $C(k)$. \epf

The strategy that we use to prove Theorem~\ref{b-i-small-thm}
(whose proof will be presented in Section~\ref{distri})
is quite involved, and so we will
now provide a brief overview of the proof before continuing.

Let $(\g,\Delta;L,\ell)$ be any $(\kappa,\lambda)$-quadruple
as in the statement of Theorem~\ref{b-i-small-thm},
and put $\g:=\Big((\gamma_i,\alpha_i,\beta_i) \Big)_{i=1}^{g+1}$
and $\T:=\T(\g,\ell)$.

By Lemma~\ref{new-unimodular} (i), for each $2\leq i\leq g$ satisfying
$(\gamma_i,\alpha_i,\beta_i)\in \underline{\g}\setminus \{(1,\lambda,\kappa-\lambda-1),(\kappa-\lambda-1,\lambda,1),
(\gamma_{g+1},\alpha_{g+1},\beta_{g+1})\}$,
there exists a closed interval
$\I=[\I_{\min},\I_{\max}]$ with $\I_{\min}<\I_{\max}$, which
we shall call a ``well-placed interval"
(see Section~\ref{wpi}), such that\\
\noindent(W1) $ \I \subseteq ( \frak{R}_1,\frak{R}_{\max})$;\\
\noindent(W2) If $\I \cap I_i \neq \emptyset$
then $ \I \subseteq I_i$ holds, $1 \leq i \leq g$;\\
\noindent(W3) $\I \subseteq I_i$ \\
all hold (cf. (\ref{R-L}), (\ref{def-guide interval})).

In the first step of the proof of
Theorem~\ref{b-i-small-thm}, we will
approximate the Christoffel numbers of
the eigenvalues of $\T$ inside a well-placed interval $\I$.
To do this, we define the quantities
\begin{eqnarray*}
\frak{c}=\frak{c}(\g, \I)&:=& \min \Big\{ \{ 2\leq i\leq g \,:\, \I_{\max} <\frak{L}_i\}\cup \{g+1\} \Big\};\\
\frak{d}=\frak{d}(\g, \I)&:=& \max \Big\{ \{ 2\leq i\leq g \,:\, \I_{\max}
<\frak{L}_i \}\cup \{\frak{c}\} \Big\}; \\
 \Gap(\I)=\Gap_{\g,\ell}(\I)&:= & \left\{
\begin{array}{ll}
\sum_{\frak{c}\leq j \leq \frak{d}} \ell(j) & \mbox{~~if~~}\frak{c}\leq g\\
0 & \mbox{~~if~~}\frak{c}=g+1,
\end{array}
\right.
\end{eqnarray*}
(cf. (\ref{c}), (\ref{d}), (\ref{def-gap})) and,
for any eigenvalue $\theta\in \I$ of $\T$, we
approximate the sum $\sum_{i=0}^{D_{\T}}\kappa_i u_i^2(\theta)$
(see Theorem~\ref{multi-thm})
by bounding the following three subsums (cf. (\ref{s(i)}), (\ref{D(T)}), (\ref{kappa_i}), (\ref{ui-form})):
\vskip0.01cm
\noindent(1) Head sum: $\sum_{i=0}^{s(\frak{a})-2}\kappa_i u_i^2(\theta)$;\\
\noindent(2) Gap sum: $\sum_{i=s(\frak{a})-1}^{s(\frak{b}+1)}\kappa_i u_i^2(\theta)$;\\
\noindent(3) Tail sum: $\sum_{s(\frak{b}+1)+1}^{D_{\T}}\kappa_i u_i^2(\theta)$.
\vskip0.01cm

We can use the theory of three-term recurrence relations,
to bound the Head sum and the Gap sum (see Theorem~\ref{0-s-sum} and Corollary~\ref{gap-sum}).
However, for the Tail sum, there may exist some real numbers near
to which we are unable to find good bounds for the Tail sum.
Let $\mathcal{B}$ denote the set of these real numbers (cf. (\ref{bad roots-def})).
In Theorem~\ref{bad-root-bd}, we show that
$\mathcal{B}$ is finite and depends only on $\g,\Delta$ and $L$.
In particular, for each $(\gamma_i,\alpha_i,\beta_i)\in
\underline{\g}\setminus
\{(1,\lambda,\kappa-\lambda-1),
(\kappa-\lambda-1,\lambda,1),
(\gamma_{g+1},\alpha_{g+1},\beta_{g+1})\}$,
there always exists a well-placed interval $\J\subseteq I_i$
such that $\J\cap \mathcal{B}=\emptyset$ (cf. Corollary~\ref{sub-WI}).
Note that such a well-placed interval $\J$ depends only on $\g,\Delta$ and $L$.
We strengthen the condition on the interval $\I$
by requiring that in addition to (W1)--(W3), it also
satisfies $\I\cap \mathcal{B}=\emptyset$.
Then for any such a well-placed interval, we can approximate the
Tail sum as long as we require that $\ell(i)>F$ holds
for all $(\gamma_i,\alpha_i,\beta_i)\in \underline{\Delta}$,
where $F$ is a positive constant depending only
on $\kappa,\g,\Delta$ and $L$ (cf. Theorem~\ref{s-d-sum}).

Now, by Condition (iii) of Theorem~\ref{b-i-small-thm}, we can find an element
$(\gamma_i,\alpha_i,\beta_i)\in \underline{\g}\setminus
\{(1,\lambda,\kappa-\lambda-1),
(\kappa-\lambda-1,\lambda,1),
(\gamma_{g+1},\alpha_{g+1},\beta_{g+1})\}$
satisfying $\ell(i)>\frac{\epsilon \HH_{\T}}{g}$, and
we can find a well-placed interval $\I\subseteq I_i$
such that $\I\cap \mathcal{B}=\emptyset$ and
$\Len(\I)> \frac{\epsilon \HH_{\T}}{g}$ both hold
(cf. (\ref{def-h}), (\ref{ell(I)})).

By the approximation given in Theorem~\ref{multi-thm} and Property
(AC), it follows that for any real number $\delta>0$, there exist
two positive constants
$C_1=C_1(\kappa,\lambda,\epsilon,\delta,\g,\Delta,L)$ and
$C_2=C_2(\kappa,\lambda,\delta)$ such that any two eigenvalues
$\theta, \eta\in \I$ of $\T$ which are conjugate algebraic numbers
must satisfy $|\theta-\eta|\leq \delta$ if $\HH_{\T}\geq C_1$ and
$\Gap(\I)\leq C_2 \HH_{\T}$ all hold (cf. Theorem
\ref{close-evs-gap(I)}). In Claim~\ref{main-prop-claim1}, we show,
by using interlacing, that the number of eigenvalues in $\I$ is at
least $C_3 \HH_{\T}$, where $C_3$ is a positive constant depending
only on $\I_{\max}-\I_{\min}, \epsilon$ and $\g$.

Now, we have to consider two cases: either $\Gap(\I)\leq C_2 \HH_{\T}$ or
$\Gap(\I)>C_2 \HH_{\T}$. In the first case, $\Gap(\I) \leq C_2 \HH_{\T}$,
we show by using Theorem~\ref{close-evs-gap(I)},
Claim~\ref{main-prop-claim1} and Theorem~\ref{limit-thm},
a result in number theory, that
\[ \lim_{\HH_{\T} \rightarrow  \infty} \frac{|\{ \eta\mbox{\,:\,eigenvalues of $\T$ that
have an algebraic conjugate in } \I \}|}{\HH_{\T}}= \infty\]
holds (cf. Proposition~\ref{main-prop}).
Since the number of eigenvalues of $\T$ is
exactly $D_{\T}+1$ (cf. (\ref{def-Eg}))
and $D_{\T}+1\leq (C(\kappa)+1)\HH_{\T}$ holds by condition (ii)
of Theorem~\ref{b-i-small-thm}, there exists a constant $H>0$
depending only on $\kappa, \lambda, \epsilon, \g,\Delta,L$
so that $\HH_{\T} \leq  H $ holds, as required.

In the second case, $\Gap(\I)>C_2 \HH_{\T}$, by
the unimodality of the sequence $(\frak{R}_i )_{i=1}^{g}$, we can find
another well-placed interval
$\I':=[\I'_{\min}, \I'_{\max}]$ which depends only on $\g,\Delta,L$ such that\\
\noindent(1) $\I'_{\min}> \I_{\max}$;\\
\noindent(2) $\I' \cap \mathcal{B}=\emptyset$;\\
\noindent(3) $\Len(\I')>\frac{\Gap(\I)}{g}$, \\
all hold (see Proposition~\ref{gap-ind}). So we can repeat the same process with $\I'$
instead of $\I$. Using the unimodality of the
sequence $(\frak{R}_i)_{i=1}^{g}$,
the condition $\I'_{\min}> \I_{\max}$
implies that $\frak{c}(\g,\I)\leq \frak{c}(\g,\I')$,
$\frak{d}(\g,\I')\leq \frak{d}(\g,\I)$ and $\Gap(\I')<\Gap(\I)$ all hold.
Hence, the second case can be repeated
at most $g$ times so that, finally, the first case must
be satisfied, from which Theorem~\ref{b-i-small-thm} again follows.

\section{Two Useful Results for Polynomials}\label{no-thy-section}

In this section, we prove two useful results concerning roots of
polynomials. The first one, Theorem \ref{ft-L-thm}, will be used
in Theorem~\ref{bad-root-bd} to show that the set $\mathcal{B}$
(as we introduced in Section \ref{central}) is finite. The second
result, Theorem~\ref{limit-thm}, analyzes the polynomials having
all roots in an interval. It will be used to bound the number of
eigenvalues of a distance-regular graph in the proof of
Proposition~\ref{main-prop}.

We denote the degree of any polynomial $p(x)$ by {\em
$\deg(p(x))$}. The polynomial $p(x)=0$ is called the {\em zero
polynomial} and, for technical reasons, we define the degree of
this polynomial to be $-1$ (cf. \cite[p.158]{hungerford}). Two
polynomials $p_1(x)$ and $p_2(x)$ are identical if their
difference $p_1(x)-p_2(x)$ is the zero polynomial.
 Let $\mathbb{R}$ and $\mathbb{C}$ be the fields of real and complex
numbers, respectively, and let $\mathbb{R}[x]$ denote the ring of
polynomials in one variable $x$ with real coefficients.

\begin{theorem}\label{ft-L-thm}
Let $q_1(x), q_2(x) \in \mathbb{R}[x]$ be two monic quadratic
polynomials which are not squares of linear polynomials, and let
$I \subseteq \mathbb{R}$ be the largest (infinite) interval on
which both $q_1(x)$ and $q_2(x)$ are non-negative. Suppose
$P_j(x)\in \mathbb{R}[x]~(1\leq j\leq 4)$ are such that $C:=\max
\{\deg(P_j(x))\,:\,1\le j\le 4\} \ge 0$. Put
\[
P(x):=P_1(x)+P_2(x)
\sqrt{q_1(x)}+P_3(x)\sqrt{q_2(x)}+P_4(x)\sqrt{q_1(x)q_2(x)}.
\]
Then the equation $P(x)=0$ has at most $4(C+2)$ roots in $I,$
unless $q_1(x)$ is identical to $q_2(x)$ and $P_2(x)+P_3(x)$,
$P_1(x)+q_1(x)P_4(x)$ are the zero polynomials, in which case
$P(x)=0$ for every $x \in I$.
\end{theorem}

\pf For each $i,j\in \{0,1\}$, we define
$$
P^{(ij)}(x) := P_1(x) +(-1)^i P_2(x) \, \sqrt{q_{1}(x)}
+ (-1)^{i+j}P_3(x) \, \sqrt{q_{2}(x)}+(-1)^jP_4(x)\,
\sqrt{q_{1}(x)q_{2}(x)}~,
$$
and put
\begin{equation}\label{product}
P^*(x) := P^{(00)}(x)\times P^{(01)}(x)\times P^{(10)}(x) \times P^{(11)}(x).
\end{equation}
Note that
$$P^{(00)}(x) P^{(01)}(x)= \Big(P_1(x)+P_2(x)\sqrt{q_1(x)}\Big)^2-q_2(x)
\Big(P_3(x)+P_4(x)\sqrt{q_1(x)}\Big)^2$$ has the form
$U(x)+V(x)\sqrt{q_1(x)}$ with $U(x), V(x) \in \mathbb{R}[x]$
satisfying  $\deg(U(x))\le 2C+4$ and $\deg(V(x)) \le 2C+2.$
Similarly,
$$P^{(10)}(x) P^{(11)}(x)= \Big(P_1(x)-P_2(x)\sqrt{q_1(x)}\Big)^2-q_2(x)
\Big(P_3(x)-P_4(x)\sqrt{q_1(x)}\Big)^2=U(x)-V(x)\sqrt{q_1(x)}.$$
Hence, by (\ref{product}), $P^*(x)=U(x)^2-V(x)^2 q_1(x)$ is a real
polynomial of degree at most $4C+8$. This proves the theorem in the
non-degenerate case when $P^*(x)$ is not the zero polynomial.

Assume now that $P^*(x)$ is the zero polynomial. We need to prove
that this happens only if $q_1(x)$ is identical to $q_2(x)$ and
$P_2(x)+P_3(x),$ $P_1(x)+q_1(x)P_4(x)$ are the zero polynomials.
We first prove the following.

\begin{claim}\label{q1 neq q2}
If $q_1(x) - q_2(x)$ is not the zero polynomial then $P^*(x)$ is
also not the zero polynomial.
\end{claim}
\noindent{\bf Proof of Claim \ref{q1 neq q2}} We first show that
$P_1^2(x)-q_1(x)q_2(x)P_4(x)^2$ is not the zero polynomial if at
least one of the polynomials $P_1(x), P_4(x)$ is not the zero
polynomial. Take a root $\gamma \in \mathbb{C}$ of $q_1(x)$ which
is not a root of $q_2(x)$. By the condition of the theorem,
$\gamma$ is the root of $q_1(x)q_2(x)$ of multiplicity $1$. Assume
that $P_1^2(x)-q_1(x)q_2(x)P_4(x)^2$ is the zero polynomial. Then
$\gamma$ is the root of $q_1(x)q_2(x)P_4(x)^2$ of odd multiplicity
but it is either not the root of $P_1(x)^2$ or it is its root of
even multiplicity, a contradiction. By the same argument,
$P_2(x)^2 q_1(x)-P_3(x)^2q_2(x)$ is not the zero polynomial if at
least one of the polynomials $P_2(x), P_3(x)$ is not the zero
polynomial. Since $C:=\max\{\deg(P_j(x))\,:\,1\le j\le 4 \}\ge 0$,
we always have either $P_1^2(x) \ne q_1(x)q_2(x)P_4(x)^2$ (if
$P_1(x)P_4(x)$ is not the zero polynomial) or $P_2(x)^2 q_1(x) \ne
P_3(x)^2q_2(x)$ (if $P_2(x)P_3(x)$ is not the zero polynomial) for
infinitely many $x \in I$.

Suppose $P^*(x)$ is the zero polynomial. Then one of the functions
$P^{(ij)}(x)$, where $i,j\in \{0,1\}$, must be zero identically on
$x \in I$. Hence
\begin{equation}\label{uuu1}
P_1(x)+(-1)^i P_2(x)\sqrt{q_1(x)}+ (-1)^{i+j}
P_3(x)\sqrt{q_2(x)}+(-1)^j P_4(x)\sqrt{q_1(x)q_2(x)}=0.
\end{equation}
Our aim is to show that this is only possible if all $P_j(x),$
$j=1,2,3,4,$ are the zero polynomials which is not the case by the
condition of the theorem.

We first claim that
\begin{equation}\label{uuu}
P_1(x)P_2(x)= q_2(x)P_3(x)P_4(x).
\end{equation}
Indeed, putting first two terms of (\ref{uuu1}) into the right
hand side and squaring we obtain
\begin{equation}\label{uuu2}
\Big(P_1(x)+(-1)^i P_2(x)\sqrt{q_1(x)}\Big)^2 = q_2(x)\Big((-1)^i
P_3(x)+ P_4(x)\sqrt{q_1(x)} \Big)^2.
\end{equation}
Since $q_1(x)$ is not the square of a linear polynomial, by the same
argument for roots multiplicity as above, the function
$S(x)+T(x)\sqrt{q_1(x)},$ where $S(x), T(x) \in \mathbb{R}[x]$, is
zero identically on $I$ if and only if $S(x)$ and $T(x)$ are the
zero polynomials. Therefore, collecting terms for $\sqrt{q_1(x)}$ in
(\ref{uuu2}) we obtain (\ref{uuu}).

Similarly, putting the first and the third term of (\ref{uuu1}) to
the right hand side, squaring and then using the same argument for
the ring $\mathbb{R}[x]+\mathbb{R}[x]\sqrt{q_2(x)},$ we deduce
that
\begin{equation}\label{uuu3}
P_1(x)P_3(x)= q_1(x)P_2(x)P_4(x).
\end{equation}

Suppose first that $P_1(x)$ is the zero polynomial. Then, by
(\ref{uuu}) and (\ref{uuu3}), $P_2(x), P_3(x)$ or $P_4(x)$ is zero
identically. If either $P_2(x)$ or $P_3(x)$ is the zero polynomial
then, by (\ref{uuu1}), all four $P_j(x)$ must be the zero
polynomials, a contradiction. If $P_4(x)$ is the zero polynomial
then $P_2(x)\sqrt{q_1(x)}+ (-1)^j P_3(x) \sqrt{q_2(x)}=0$. But this
yields $P_2(x)^2 q_1(x)=P_3(x)^2 q_2(x),$ a contradiction again. By
the same argument, if any of the polynomials $P_2(x), P_3(x),
P_4(x)$ is the zero polynomial, then by (\ref{uuu}) and (\ref{uuu3})
one more polynomial must be a zero polynomial. One then concludes as
above that all four polynomials are the zero polynomials.

Finally, if none of the polynomials $P_j(x)$ is the zero polynomial
then multiplying (\ref{uuu}) and (\ref{uuu3}) gives $P_1(x)^2
P_2(x)P_3(x)= q_1(x)q_2(x) P_2(x)P_3(x)P_4(x)^2$. Hence
$P_1(x)^2=q_1(x)q_2(x)P_4(x)^2$, which is a contradiction again.
\epf

Now, to complete the proof of the theorem, suppose that $q_1(x)$ is
identical to $q_2(x)$. Then
$P(x)=P_1(x)+q_1(x)P_4(x)+(P_2(x)+P_3(x))\sqrt{q_1(x)}$ for all
$x\in I$. If $P_1(x)+q_1(x)P_4(x)$ and $P_2(x)+P_3(x)$ are the zero
polynomials then $P(x)$ is zero identically. Otherwise,
$$P(x)(P_1(x)+q_1(x)P_4(x)-(P_2(x)+P_3(x))\sqrt{q_1(x)})
=(P_1(x)+q_1(x)P_4(x))^2- q_1(x)(P_2(x)+P_3(x))^2$$ is not the zero
polynomial. So $P$ has at most $2C+4$ roots in $I$, which is better
than required. \epf

In the remainder of this section, we will show the second useful
result,
Theorem~\ref{limit-thm}.\\
For any real number $\kappa\geq 2$, we denote by
$\mathcal{P}_{\kappa}$ the set of all irreducible monic
polynomials $p(x) \in \mathbb{Z}[x]$ such that all of the roots of
$p(x)$ are contained in the closed interval $[-\kappa,\kappa]$.
Note $\mathcal{P}_{\kappa} \subseteq \mathcal{P}_{\kappa'}$ if
$\kappa \le \kappa'$.

\begin{lemma}\label{no-thy-2}
Let $\kappa \geq 2$ be a real number and let $n$ be a positive
integer.
Then the following holds.\\
(i) The set consisting of all polynomials $p(x)\in
\mathcal{P}_{\kappa}$ of degree at most $n$ is finite.
\\
(ii) $\mathcal{P}_{\kappa}$ is an infinite set.
\end{lemma}
\pf (i) Obviously, any coefficient of each $p(x) \in
\mathcal{P}_{\kappa}$ of degree at most $n$ is in $[-(2\kappa)^n,
(2\kappa)^n],$ so $\mathcal{P}_{\kappa}$ contains at most $(2
(2\kappa)^n+1)^{n+1}$ of such polynomials. See also \cite[Lemma
7.1]{val5}.

\noindent (ii) Let $P_n(1,0,1)(x)$ be the characteristic
polynomial of the tridiagonal $(n \times n)$-matrix with zeroes on
the diagonal and ones on the subdiagonals and superdiagonals. Then
$P_n(1,0,1)(x)$ is a polynomial of degree $n$ and has $n$ distinct
roots, $2 \cos(\frac{i\pi}{n+1})$, $i =1, \ldots, n$
(\cite[p.11]{biggs}). Thus, if we factorize $P_n(1,0,1)(x)$ into
irreducible factors, say $q_1(x), \ldots, q_t(x)$, then $q_i(x)
\neq q_j(x)$ if $1 \leq i < j \leq t$ and $q_i(x) \in
\mathcal{P}_{2}$ for all $ 1 \leq i \leq t$. (ii) now follows
immediately from (i).\epf

In fact, an old result of R. M. Robinson \cite{robb} asserts that
if $J$ is an interval of length strictly greater than $4$ then
there are infinitely many irreducible monic polynomials whose
roots all lie in $J$. Moreover, none of them has a root of the
form $2\cos(\pi r)$ with $r \in \mathbb{Q}$ as those lying in
$\mathcal{P}_{2}$.

Now, for any real number $\zeta >0$, let ${\bf I_{\kappa, \zeta}}$
be the set of all closed intervals of length $\zeta$ which are
contained in the closed interval $[-\kappa, \kappa]$. For each
$p\in \mathcal{P}_{\kappa}$ and $I\in {\bf I_{\kappa,\zeta}}$, we
define
\begin{eqnarray}
\Upsilon_{\kappa}(p,I)&:=& \frac{\left| \{\theta \in I \,:\, p(\theta)=0\}\right|-1}{\deg(
p(x) ) }\,, \mbox{ and } \label{upsilonbase}\\
\Upsilon_{\kappa,\,\zeta}&:=&\sup \left\{ \Upsilon_{\kappa}(p,I)
\,:\, p\in \mathcal{P}_{\kappa},~I\in {\bf
I_{\kappa,\zeta}}\right\}. \label{def-Gamma}
\end{eqnarray}

\begin{remark}\label{upsilon>0}
Note that $\Upsilon_{\kappa,\,\zeta}$ is positive for all $\zeta
>0$ since by Lemma~\ref{no-thy-2} (ii) there exists a polynomial
$p(x) \in \mathcal{P}_{2}$ with degree $n > \frac{8
\kappa}{\zeta}$ and so, by the pigeon hole principle, there exists
an interval $I \in {\bf I_{\kappa,\zeta}}$ of length $\zeta$ such
that $p(x)$ has at least $\frac{n \zeta}{4 \kappa}$ roots in $I$.
Even so, we now show that the limit of $\Upsilon_{\kappa,\,\zeta}$
as $\zeta$ tends to $\infty$ is zero.
\end{remark}

\begin{theorem}\label{limit-thm}
Let $\kappa\geq 2$ be a real number. Then
\[
\lim_{\zeta \rightarrow 0 }\Upsilon_{\kappa,\,\zeta}=0.
\]
\end{theorem}

\pf Fix $\kappa\geq 2$ and $\zeta \in (0,1)$. Let $p(x) \in
\mathcal{P}_{\kappa}$ be of degree $n$, say, and let $I \in
\mathbf{I}_{\kappa,\,\zeta}$. Since $p(x)$ is irreducible in
$\mathbb{Z}[x]$, it has $n$ distinct roots $\alpha_1, \ldots,
\alpha_n \in [-\kappa, \kappa]$. Consider the discriminant
$\Delta(p)$ of $p$ given by
$$
\Delta(p) := \prod_{1 \leq i < j \leq n} (\alpha_i - \alpha_j)^2.
$$
Since $p(x)$ is a monic polynomial with integral coefficients, its
discriminant $\Delta(p)$ is an integer. Moreover,
$\Delta(p)$ is not zero as the roots of $p(x)$ are distinct and $\Delta(p) > 0$,
so $\Delta(p) \geq 1$.

Without loss of generality, assume that
$\{ \alpha_1, \ldots, \alpha_t \}$ is the set of
roots of $p(x)$ contained in $I$, for some $0 \le t \leq n$.
Let $\tau = \tau(p,I) := \frac{t}{n}$.

\begin{claim}\label{no-claim}
If $t \geq 2$ then $\tau^2 \leq -\frac{2 \ln{(2\kappa)}}{\ln
\zeta}$.
\end{claim}
{\bf Proof of Claim:} We have
\begin{eqnarray}
1&\leq & \prod_{1 \leq i<j \leq n}(\alpha_i-\alpha_j)^2 \nonumber \\
&= & \left( \prod_{1\leq i<j\leq t}(\alpha_i-\alpha_j)^2\right)
\left( \prod_{1 \leq i<j \leq n
~\mbox{ \small{and} $j>t$}}(\alpha_i-\alpha_j)^2\right)\nonumber \\
&\leq & \zeta^{\tau n (\tau n-1)} (2 \kappa )^{n (n-1)},\nonumber
\end{eqnarray}
since $t = \tau n \geq 2$, $|\alpha_i - \alpha_j| \leq \zeta$ for
$1 \leq i < j \leq t$ and $|\alpha_i - \alpha_j| \leq 2\kappa$ for
$1 \leq i < j \leq n$. Using $\tau n-1 \ge \tau (n-1)/2$ and
$0<\zeta < 1$ we find that $1 \leq \zeta^{\frac{\tau n \tau
(n-1)}{2}} (2 \kappa )^{n (n-1)}$, so $1 \le \zeta^{\tau^2/2}
2\kappa$. The claim follows by taking the logarithms of both sides
of the last inequality. \epf

Now, let $q(x) \in \mathcal{P}_{\kappa}$ and $I \in
\mathbf{I}_{\kappa,\,\zeta}$ be such that
$$
\Upsilon_{\kappa}(q,I) \ge \frac{1}{2}\Upsilon_{\kappa,\,\zeta}>0.
$$
Such a $q(x)$ exists, since, as remarked before the statement of
the theorem, $\Upsilon_{\kappa,\,\zeta}$ is positive. Since
$\Upsilon_{\kappa}(q,I)>0$, the polynomial $q(x)$ has at least 2
roots in $I$. Hence, by Claim \ref{no-claim} and
(\ref{upsilonbase}), we have
$$
\sqrt{-\frac{2 \ln{(2\kappa)}}{\ln \zeta}} \ge \frac{|\{ x\in I
\,:\, q(x) =0\}|}{\deg(q(x))} >\frac{|\{ x\in I \,:\, q(x)
=0\}|-1}{\deg(q(x))} =\Upsilon_{\kappa}(q,I) \ge
\frac{1}{2}\Upsilon_{\kappa,\,\zeta}>0,
$$
from which the theorem immediately follows.
\epf

\section{Preliminary Results for the Christoffel Numbers}\label{bad-roots-section}

In this section, we will prove some results which we
will use later in Section \ref{multi-section}
for the approximation of Christoffel numbers.

Suppose that
$\g:=\Big((\gamma_i,\alpha_i,\beta_i)\Big)_{i=1}^{g+1}$
is a $(\kappa,\lambda)$-graphical sequence, that
$(\g, \Delta;L,\ell)$ is a $(\kappa,\lambda)$-quadruple
as in Definition~\ref{RT-def}, and that $\Delta=(\delta_j)_{j=1}^{\tau}$.

Fix $i$ with $0 \leq i \leq \tau-1$.
Let $0\leq j_i <g$ be the integer for which
\begin{equation}\label{j(k)}
\delta_{\tau-i} =(\gamma_{g-j_i},\alpha_{g-j_i},\beta_{g-j_i})
\end{equation}
holds. We put $j_{-1}:=-1$, and note that
$j_i -j_{(i-1)}\geq 1$ necessarily holds.

Suppose $j_i-j_{(i-1)}\geq 2$. Then, for $n_i :=j_i-j_{(i-1)}-1$, we
define the sequence $z^{(i)}=(z^{(i)}_s)_{s=1}^{n_i}$ by putting
\[
z^{(i)}_s:=(\gamma_{g-j_{(i-1)}-s},\alpha_{g-j_{(i-1)}-s},
\beta_{g-j_{(i-1)}-s}),~~ s=1,\ldots,n_i.
\]
In addition,  for $N:=\sum_{z \in \underline{z^{(i)}}}L(z)$,
we let $w^{(i)}=(w^{(i)}_k)_{k=1}^N$ be the sequence
whose $k$th term $w_k^{(i)}$ is defined to be $z^{(i)}_j$ for
the necessarily unique $j$ for which
\begin{equation}\label{seq-w}
\sum_{s=1}^{j-1} L(z^{(i)}_s) < k \leq \sum_{s=1}^j L(z^{(i)}_s)
\end{equation}
holds.

Now, suppose that $\theta$ is a real number,
and that $v_0$ and $v_1$ are real numbers satisfying
$(v_0,v_1)\neq (0,0)$. In addition, let
$(v_j)_{j=0}^{N+1}$ be the sequence that
is defined by the recurrence relations
\begin{equation}\label{w-three-terms}
\tilde{\beta}_j v_{j-1} + (\tilde{\alpha}_j - \theta) v_j + \tilde{\gamma}_j v_{j+1} = 0 \ \ \ (j=1,2,
\ldots, N)~,
\end{equation}
where $(\tilde{\gamma}_j,\tilde{\alpha}_j, \tilde{\beta}_j)$ denotes the $j$th term $w_j^{(i)}$ of the
sequence $w^{(i)}$, and $N$ is as above if $j_i-j_{(i-1)}\geq 2$
and $N:=0$ else. Then, in view of (\ref{w-three-terms}),
for $j_i-j_{(i-1)}\geq 2$
there are polynomials $f^{(i)}_t(x)$, $g^{(i)}_t(x)$ in
$\mathbb{Q}[x]$ (of degree $s-1$ and $s-2$, respectively)
that, for $\theta \in [\frak{R}_{g-j_i},\kappa]$,
satisfy  $v_s = f^{(i)}_t(\theta) v_1 + g^{(i)}_t(\theta) v_0$
for each $s \geq 1$,
\begin{equation}\label{pi_1-f}
v_{N}=\left\{
\begin{array}{ll}
f_1^{(i)}(\theta) v_1+g_1^{(i)}(\theta)
v_0 & \mbox{~if~$j_i-j_{(i-1)}\geq 2$} \\
v_0 & \mbox{~if~$j_i-j_{(i-1)}=1$,}
\end{array}
\right.
\end{equation}
and
\begin{equation}\label{pi_1-f-1}
v_{N+1}=\left\{
\begin{array}{ll}
f_2^{(i)}(\theta) v_1+g_2^{(i)}(\theta) v_0 &\mbox{~if~$j_i-j_{(i-1)}\geq 2$} \\
v_1 & \mbox{~if~$j_i-j_{(i-1)}=1$.}
\end{array}
\right.
\end{equation}
In addition, in case $j_i-j_{(i-1)}=1$, we
let $f_t^{(i)}(x)$ and $g_t^{(i)}(x)~~(t=1,2)$ be
the polynomials in $\mathbb{Q}[x]$ for which
both $f_t^{(i)}(x)-t+1$ and $g_t^{(i)}(x)+t-2$ are
the zero polynomials for $t=1,2$. Note that
the degrees of the polynomials $f^{(i)}_t(x)$ and $g^{(i)}_t(x)$  are as follows:
\begin{eqnarray}
\deg(f_1^{(i)}(x))&=&\left\{
\begin{array}{ll}\label{degree-ftn}
-1+\sum_{z\in \underline{z^{(i)}}} L(z) &\mbox{~~if~~}j_i-j_{(i-1)}\geq 2\\
-1 &\mbox{~~if~~}j_i-j_{(i-1)}=1,
\end{array}
\right. ;\\
\deg(g_1^{(i)}(x))&=&\left\{
\begin{array}{ll}
\deg(f_1^{(i)}(x))-1 &\mbox{~~if~~}j_i-j_{(i-1)}\geq 2 \\
0 &\mbox{~~if~~}j_i-j_{(i-1)}=1
\end{array}
\right. ;\\
\deg(f_2^{(i)}(x))&=&\deg(f_1^{(i)}(x))+1 ;\\
\deg(g_2^{(i)}(x))&=&\deg(f_1^{(i)}(x)).\label{deg-f-g}
\end{eqnarray}
Note also that $f_t^{(i)}(x)$ and
$g_t^{(i)}(x)$ ($t=1,2$, $0\leq i\leq \tau -1$)
depend only on the triple $(\g,\Delta,L)$
(and not on the function $\ell$).

We now present the second key definition
of this section. For the $(\kappa,\lambda)$-graphical sequence
$\g=\Big((\gamma_i,\alpha_i,\beta_i)\Big)_{i=1}^{g+1}$, let
$x_{i}=x_i(\theta)$ and $y_{i}=y_i(\theta)$ (where
$|x_{i}|\,\geq\,|y_{i}|$) be the roots of the equation
\begin{equation}\label{def-x_i-y_i}
\gamma_{g-i}x^2+(\alpha_{g-i}-\theta)x+\beta_{g-i}=0~(0\le i < g).
\end{equation}

\begin{defi}\label{def-bad roots}
For any integers $\kappa \geq 3$ and $\lambda\geq 0$ with $\lambda
\leq \kappa-2$, let $(\g, \Delta; L,\ell)$ be a
$(\kappa,\lambda)$-quadruple with
$\g=\Big((\gamma_i,\alpha_i,\beta_i)\Big)_{i=1}^{g+1}$ and
$\Delta=(\delta_i)_{i=1}^{\tau}$. With reference to (\ref{R-L}),
(\ref{pi_1-f}), (\ref{pi_1-f-1}) and (\ref{def-x_i-y_i}), for $0
\leq i \leq \tau-1$ and
$\delta_{\tau-i}=(\gamma_{g-j_i},\alpha_{g-j_i},\beta_{g-j_i})\in
\underline{\Delta}$ satisfying $\beta_{g-j_i}\leq \gamma_{g-j_i}$,
and for any real numbers $v_0,v_1$ satisfying $(v_0,v_1)\neq
(0,0)$ we define the set
$\mathcal{B}_{i}=\mathcal{B}_{i}(\g,\Delta,L)$ by
\begin{equation}\label{B_i}
\mathcal{B}_{i}(\g,\Delta,L):=\left\{\theta \in
[\frak{R}_{g-j_i}(\g),\frak{R}_{\max}(\g)]\,:\, F_i(\theta)=0 \right\},
\end{equation}
where $F_i(x)$ is the polynomial in $\mathbb{R}[x]$ given by
\[
F_i(x):=\left\{
\begin{array}{ll}
\prod_{\xi \in \{x_{j_i},y_{j_i} \}} \left((x
-\alpha_{g+1})(f_1^{(i)}(x)\xi-f_2^{(i)}(x))+
\gamma_{g+1}(g_1^{(i)}(x)\xi-g_2^{(i)}(x))\right)
 &\mbox{if~$i=0$}\\
\prod_{(\xi,\chi) \in \{x_{j_i},y_{j_i}\}
\times \{x_{j_{(i-1)}},y_{j_{(i-1)}}\}}\left(
(f_1^{(i)}(x)\xi-f_2^{(i)}(x))\chi+
g_1^{(i)}(x)\xi-g_2^{(i)}(x)\right)
&\mbox{if~}i \ne 0.
\end{array}
\right.
\]
With reference to (\ref{B_i}), we also define the set
$\mathcal{B}=\mathcal{B}(\g,\Delta,L)$ by
\begin{equation}\label{bad roots-def}
\mathcal{B}(\g,\Delta,L):= \bigcup_{0\leq i\leq \tau-1
\mbox{~and~} \beta_{g-j_i}\leq \gamma_{g-j_i}}
 \mathcal{B}_i(\g,\Delta,L).
\end{equation}
\end{defi}

Note that since the polynomials $f_t^{(i)}(x)$ and $g_t^{(i)}(x)$
($t=1,2$, $0\leq i\leq \tau -1$) depend only on the triple
$(\g,\Delta,L)$, the polynomial $F_i(x)$ and the sets
$\mathcal{B}_i$ and $\mathcal{B}$ in the last definition all also
depend only on $(\g,\Delta,L)$. Note that, if
$\mathcal{R}_{\max}(\g)=\kappa$, then $F_i(\kappa)=0$ for all $i$,
$0\leq i\leq \tau-1$ (as the standard eigenvector for $\kappa$ is
the all-one vector and $x_i=1$ for all $0 \leq i \leq \tau-1$),
and hence in this case $\kappa \in \mathcal{B}$.

\begin{theorem}\label{bad-root-bd}
Let $\kappa \geq 3$ and $\lambda\geq 0$ be integers with $\lambda
\leq \kappa-2$. Suppose that $(\g, \Delta; L,\ell)$ is a
$(\kappa,\lambda)$-quadruple. Then there exists a constant
$C=C(\g,\Delta,L)>0$ such that
\[\left| \mathcal{B}\right| \leq C\]
holds, for $\mathcal{B}=\mathcal{B}(\g,\Delta,L)$
as defined in Definition~\ref{def-bad roots}.
\end{theorem}

\pf  Let $(\g,\Delta;L,\ell)$ be a $(\kappa,\lambda)$-quadruple,
put $\g=\Big((\gamma_i,\alpha_i,\beta_i)\Big)_{i=1}^{g+1}$ and let
$\T=\T(\g,\ell)$ be the associated $(\kappa,\lambda)$-tridiagonal
sequence. In addition, put $\Delta=(\delta_j)_{j=1}^{\tau}$ and
let $\delta_{\tau-i}=(\gamma_{g-j_i},\alpha_{g-j_i},\beta_{g-j_i})\in
\underline{\Delta}$ with $\beta_{g-j_i}\leq \gamma_{g-j_i}$ be as
defined in (\ref{j(k)}).

To prove the theorem, we will use Theorem~\ref{ft-L-thm} to
bound $|\mathcal{B}_i|$ by some constant depending only
on $\g$, $\Delta$ and $L$ for each $0 \le i \le \tau-1$.
To do this, we first define polynomials $q_s(x)$,
$s=1,2$, and $P_j(x)$,
$1 \le j \leq 4$ as in the statement of that theorem,
breaking this definition into cases depending on $i$:\\

(a) $i=0$: Let
$q_s(x)=(x-\alpha_{g-j_i})^2-4\beta_{g-j_i}\gamma_{g-j_i}~(s=1,2)$.
For each $\xi \in \{x_{j_i},y_{j_i}\}$, put
\begin{eqnarray*}
P_1(x)&:=&\frac{(x-\alpha_{g-j_i})((x-\alpha_{g+1})f_1^{(i)}
+\gamma_{g+1}g_1^{(i)})}{2\gamma_{g-j_i}}
-(x-\alpha_{g+1})f_2^{(i)}-\gamma_{g+1}g_2^{(i)},\\
P_2^{\xi}(x)&:=&(-1)^{\delta_{\xi,y_{j_i}}}
\left(\frac{(x-\alpha_{g+1})f_1^{(i)}
+\gamma_{g+1}g_1^{(i)}}{2\gamma_{g-j_i}}\right)\,,\\
P_3(x)&=&P_4(x)=0\,,
\end{eqnarray*}
where $\delta_{\xi,y_{j_i}}$ is the Kronecker
delta function, and let $P_2(x) = P_2^{\xi}(x)$. Then, for this
specific  choice of polynomials, the polynomial
$P(x) = P^{\xi}(x)$ in Theorem~\ref{ft-L-thm} becomes
$$
P(x)=
(x-\alpha_{g+1})(f_1^{(i)}(x)\xi -f_2^{(i)}(x))+
\gamma_{g+1}(g_1^{(i)}(x)\xi -g_2^{(i)}(x)),
$$
which is precisely the factor
that appears in the definition of
the polynomial $F_0(x)$ in
Definition~\ref{def-bad roots}.

Note that if $P_1(x)$ and $P_2(x)$ are the zero polynomials,
then $(x-\alpha_{g+1})f_s^{(i)}+\gamma_{g+1}g_s^{(i)}$ are
also the zero polynomials for $s=1,2$. This contradicts
(\ref{degree-ftn})--(\ref{deg-f-g}). Hence
$\max\{\deg(P_j(x))\,:\,1\leq j\leq 4\}\geq 0$.\\

\noindent(b) $i \ge 1$: Let
$q_1(x)=(x-\alpha_{g-j_i})^2-4\beta_{g-j_i}\gamma_{g-j_i}$
and
$q_2(x)=(x-\alpha_{g-j_{(i-1)}})^2-4\beta_{g-j_{(i-1)}}\gamma_{g-j_{(i-1)}}$.
Note that as $\beta_{g-j_{(i-1)}}\leq \beta_{g-j_i}\leq
\gamma_{g-j_i}\leq \gamma_{g-j_{(i-1)}}$ holds,
$q_1(y)\neq q_2(y)$ for some $y\in \mathbb{R}$.
For each $(\xi,\chi) \in \{x_{j_i},y_{j_i}\}
\times \{x_{j_{(i-1)}},y_{j_{(i-1)}}\}$, let
$P_s(x)\in \mathbb{Q}[x]$, $1\leq s\leq 4$, be
polynomials such that
$P(x):=(f_1^{(i)}(x)\xi -f_2^{(i)}(x))\chi
+ g_1^{(i)}(x)\xi - g_2^{(i)}(x)$, i.e.
the factor appearing in the
definition of the polynomial $F_i(x)$ in
Definition~\ref{def-bad roots}, $i \ge 1$. Note that if $P_s(x)=0~(1\leq
s\leq 4)$ are all the zero polynomials, then so are the polynomials
$f_t^{(i)}(x)$ and $g_t^{(i)}(x)$, $t=1,2$. Thus
$v_N(\theta)=v_{N+1}(\theta)=0$ hold for any real number $\theta$, which is
impossible as $(v_0,v_1)\neq (0,0)$. Hence
$\max\{\deg(P_s(x))\,:\,1\leq s \leq 4\}\geq 0$. \\

With these definitions in hand we can now apply
Theorem~\ref{ft-L-thm} simultaneously to cases (a) and (b).
(Clearly, $q_1(x)$ and $q_2(x)$ are not squares of linear polynomials.)
In particular, in view of (\ref{degree-ftn})--(\ref{deg-f-g}), \\
$$
\left|\mathcal{B}_i \right|\leq \left\{
\begin{array}{ll}
8\left(4+\deg(f_1^{(i)}(x))\right) &\mbox{~~if~~}i=0\\
16\left(4+\deg(f_1^{(i)}(x))\right) &\mbox{~~if~~}i\neq 0
\end{array}
\right.
$$
holds, from which the proof of the theorem now follows by taking
\[C(\g,\Delta,L):=16\,|\g|\,\left(
3+\sum_{(\gamma_{i},\alpha_{i},\beta_{i})\in
\underline{\g}\setminus \underline{\Delta}} L(i)\right).\] \epf

Now, for the $(\kappa,\lambda)$-quadruple
$(\g, \Delta;L,\ell)$, let $\T=\T(\g,\ell)$ be the
associated tridiagonal sequence. Let $\theta$ be a real number. Then, for each $0\leq i<g$ satisfying
$\theta>\frak{R}_{g-i}$, there exist complex numbers
$\nu_1^{(i)}(\theta)$ and $\nu_2^{(i)}(\theta)$ such that the
terms in the standard sequence $\Big(u_j=u_j(\theta)\Big)_{j=0}^{D_{\T}}$
satisfy
\begin{equation}\label{rec-tail}
u_{s(g-i+1)-j}(\theta)=\nu_1^{(i)}(\theta)x_{i}^j(\theta)+\nu_2^{(i)}(\theta)y_{i}^j(\theta)~~(0\le j\le \ell(g-i)+1),
\end{equation}
where $s(g-i+1)$, $x_i(\theta)$ and $y_i(\theta)$ are as defined in (\ref{s(i)}) and (\ref{def-x_i-y_i}).
Note that $x_i(\theta)-y_i(\theta) \neq 0$ holds, and that $(\nu_1^{(i)}(\theta),\,\nu_2^{(i)}(\theta))\neq (0,0)$
holds as $(u_0,u_1)\neq (0,0)$.
Taking $j=0,1$ in (\ref{rec-tail}) we obtain:
\begin{eqnarray}
\nu_1^{(i)}(\theta)  &=&  \left(
\frac{-y_i(\theta)}{x_i(\theta)-y_i(\theta)} \right)
u_{s(g-i+1)}(\theta) + \left( \frac{1}{x_i(\theta)-y_i(\theta)}
\right) u_{s(g-i+1)-1}(\theta) ;\label{nu1i}\\ \nu_2^{(i)}
(\theta) &=& \left( \frac{x_i(\theta)}{x_i(\theta)-y_i(\theta)}
\right)
u_{s(g-i+1)}(\theta)+\left(\frac{-1}{x_i(\theta)-y_i(\theta)}\right)
u_{s(g-i+1)-1}(\theta).\label{nu2i}
\end{eqnarray}

In particular, in view of (\ref{pi_1-f}), (\ref{pi_1-f-1}) and (\ref{rec-tail}), for each
$\delta_{\tau-i}=(\gamma_{g-j_i},\alpha_{g-j_i},\beta_{g-j_i})\in
\underline{\Delta}$, there exist polynomials
$f_t^{(i)}=f_t^{(i)}(x)$, $g_t^{(i)}=g_t^{(i)}(x)$
$(t=1,2)$ in $\mathbb{Q}[x]$ such that
\begin{equation}\label{u-pi_1-f}
u_{s(g-j_i+1)}(\theta)=\left\{
\begin{array}{ll}
f_1^{(i)}(\theta) u_{s(g-j_{(i-1)})-1}(\theta)+g_1^{(i)}(\theta) u_{s(g-j_{(i-1)})}(\theta)& \mbox{~if~$i\neq 0$} \\
f_1^{(i)}(\theta) u_{D-1}(\theta)+g_1^{(i)}(\theta) u_D(\theta)  &  \mbox{~if~$i=0$}
\end{array}
\right.
\end{equation}
and
\begin{equation}\label{u-pi_1-f-1}
u_{s(g-j_i+1)-1}(\theta)=\left\{
\begin{array}{ll}
f_2^{(i)}(\theta) u_{s(g-j_{(i-1)})-1}(\theta)+g_2^{(i)}(\theta) u_{s(g-j_{(i-1)})}(\theta) &\mbox{~if~$i\neq 0$} \\
f_2^{(i)}(\theta) u_{D-1}(\theta)+g_2^{(i)}(\theta) u_D(\theta) &  \mbox{~if~$i=0$}
\end{array}
\right.
\end{equation}
hold, where $f_t^{(i)}(x)-t+1$ and
$g_t^{(i)}(x)+t-2$ are the zero polynomials if $j_i-j_{(i-1)}=1$.\\

The last theorem of this section
will play an important role later on in obtaining an upper bound
for the Christoffel numbers of any eigenvalue
of $\T(\g,\ell)$ within some closed interval not intersecting $\mathcal{B}$.
For any non-empty closed real interval $I$, we define $I_{\min}$
and $I_{\max}$ to be the real numbers for which $I=[I_{\min},I_{\max}]$
holds.

\begin{theorem}\label{Lj}
Let $\kappa \geq 3$ and $\lambda\geq 0$ be integers with $\lambda
\leq \kappa-2$. Suppose that $(\g, \Delta; L,\ell)$ is a
$(\kappa,\lambda)$-quadruple and let
$\g=\Big((\gamma_i,\alpha_i,\beta_i)\Big)_{i=1}^{g+1}$. Suppose
that $I$ is a non-empty, closed subinterval of $(\frak{R}_1,\frak{R}_{\max})$
such that
 \begin{equation}\label{I-L_j}
I\cap \Big(\mathcal{B}\cup \{\frak{R}_i\,:\,1\leq i\leq g\}\Big)=\emptyset \mbox{  and  } 2\leq \frak{b}<g
\end{equation}
both hold, where $\frak{b}=\frak{b}(\g,I):=\max \{ 2\leq i\leq g
\,:\, I_{\max }< \frak{R}_i \}$. Then for each
$(\gamma_{g-i},\alpha_{g-i},\beta_{g-i})\in \underline{\Delta}$
with $\frak{b}+1\leq g-i\leq g$, there exist positive constants
$C_i=C_i(\kappa,\g,\Delta,L,I)\geq 1$ and $M_i=M_i(\kappa,
\g,\Delta,L,I)>1$ such that, if $\ell(g-j)>C_i$ hold for all $j<i$
with $(\gamma_{g-j},\alpha_{g-j},\beta_{g-j})\in
\underline{\Delta}$, then
\begin{equation}\label{eq-L_i}
\left| \frac{\nu_1^{(i)}(\theta) }{\nu_2^{(i)}(\theta)} \right| >
M_{i}\left( \frac{y_i(\theta)}{x_i(\theta)}  \right)^{C_{i}}
\end{equation}
holds for any real number $\theta\in I$, where $\frak{R}_i$,
$\underline{\Delta}$, $x_i(\theta),~y_i(\theta)$,
$\mathcal{B}=\mathcal{B}(\g,\Delta,L)$ and
$\nu_j^{(i)}(\theta)~(j=1,2)~$ are as defined in (\ref{R-L}),
(\ref{underline}), (\ref{def-x_i-y_i}), (\ref{bad roots-def}) and
(\ref{rec-tail}), respectively.
\end{theorem}

\pf Let $\mathcal{T}=\mathcal{T}(\g, \ell)$ be the
$(\kappa,\lambda)$-tridiagonal sequence associated to
$(\g,\Delta;L,\ell)$, let $D:=D_{\T}$ and let
$\Delta=(\delta_i)_{i=1}^{\tau}$. Note that for each $0\leq
j<g-\frak{b}$, $0<\frac{y_j(\theta)}{x_j(\theta)}<1$ holds for any
$\theta\in I$, and also that $\frac{y_j(\theta)}{x_j(\theta)}$ is a non-zero
continuous function (in $\theta$) on the closed interval $I$. Hence, there
exists a constant $0< P=P(\g,I)<1$ such that
\begin{equation}\label{ri}
\frac{y_j(\theta)}{x_j(\theta)} = \frac{\theta -\alpha_{g-j}-\sqrt{(\theta
-\alpha_{g-j})^2 -4 \beta_{g-j} \gamma_{g-j}}}{\theta
-\alpha_{g-j}+ \sqrt{(\theta -\alpha_{g-j})^2-4 \beta_{g-j}
\gamma_{g-j}}} \leq P  <1
\end{equation}
holds for any $0\leq j<g-\frak{b}$ and for any $\theta\in I$. Note
also that for each $0\leq j < g-\frak{b}$, $\beta_{g-j}\leq
\gamma_{g-j}$ holds by Lemma~\ref{new-unimodular}.

Now, for each
$\delta_{\tau-s}\in \underline{\Delta}$, let
\[\delta_{\tau-s}:=(\gamma_{g-j_s},\alpha_{g-j_s},\beta_{g-j_s})\]
for some $0\leq j_s<g$. We prove the theorem by
induction on $0\leq s\leq \frak{s}$, where
$$\frak{s}:=\max \{i\,:\,\delta_{\tau-i}\in
\underline{\Delta}\mbox{~~and~~}\frak{R}_{g-j_i}<I_{\min}\}.
$$

First, suppose $s=0$. Let
$\delta_{\tau}=(\gamma_{g-j_0},\alpha_{g-j_0},\beta_{g-j_0})\in
\underline{\Delta}$ for some $0\leq j_0 \leq g-\frak{b}-1$.\\
By (\ref{u_i-three-terms}),
$\gamma_{g+1 }u_{D-1}(\theta)+(\alpha_{g+1}-\theta)u_D(\theta)=0$
and $u_D(\theta)\neq 0$ for any $\theta\in I$. By
(\ref{nu1i})--(\ref{u-pi_1-f-1}), there exist polynomials
$f_t^{(0)}(x),\,g_t^{(0)}(x)\in
\mathbb{Q}[x]~(t=1,2)$ such that
\begin{eqnarray*}
\nu_1^{(j_0)}(\theta)&=& \left(\frac{(\theta
-\alpha_{g+1})(-f_1^{(0)}(\theta)y_{j_0}(\theta)+f_2^{(0)}(\theta))+\gamma_{g+1}(-g_1^{(0)}(\theta)y_{j_0}(\theta)+
g_2^{(0)}(\theta))}
{\gamma_{g+1}(x_{j_0}(\theta)-y_{j_0}(\theta))} \right) u_D(\theta) \mbox{~~~and }\\
\nu_2^{(j_0)}(\theta)&=&\left(\frac{(\theta-\alpha_{g+1})(f_1^{(0)}(\theta)x_{j_0}(\theta)-f_2^{(0)}(\theta))+
\gamma_{g+1}(g_1^{(0)}(\theta)
x_{j_0}(\theta)-g_2^{(0)}(\theta))}
{\gamma_{g+1}(x_{j_0}(\theta)-y_{j_0}(\theta))}\right) u_D(\theta)
\end{eqnarray*}
both hold for all $\theta\in I$. It follows by (\ref{I-L_j}) and
$u_D(\theta)\neq 0$ that $\nu_t^{(j_0)}(\theta)\neq 0~(t=1,2)$
for all $\theta\in I$. Since the function
$\frac{\nu_1^{(j_0)}(\theta)}{\nu_2^{(j_0)}(\theta)}$ is a
non-zero continuous function on the closed interval $I$, by
(\ref{ri}) there exist constants
$N_{j_0}:=N_{j_0}(\kappa,\mathcal{G},\Delta,L,I)>0$ and
$C_{j_0}:=C_{j_0}(\kappa,\mathcal{G},\Delta,L,I)\geq 1$ so that
\begin{equation}
\left| \frac{ \nu_1^{(j_0)}(\theta) }{ \nu_2^{(j_0)}(\theta) }
\right| \geq N_{j_0}
> P^{C_{j_0}-1}\geq \,\frac{1}{P}\left(\frac{y_{j_0}(\theta)}{x_{j_0}(\theta)} \right)^{C_{j_0}}~
\end{equation}
holds for any $\theta\in I$. Hence there exist constants
$C_{j_0}=C_{j_0}(\kappa,\g,\Delta,L,I)\geq 1$ and
$M_{j_0}=M_{j_0}(\kappa, \g,\Delta,L,I):=1/P$ such that
(\ref{eq-L_i}) holds for all $\theta \in I$. This completes the
proof of the base case.

Now, suppose $0<s\leq \frak{s}$,
and assume that the theorem holds for
all $\delta_{\tau-j}\in \underline{\Delta}$ with $0\leq j<s$.
Let $x:=x_{j_{(s-1)}}(\theta)$, $y:=y_{j_{(s-1)}}(\theta)$,
$\nu_t:=\nu_t^{(j_{(s-1)})}(\theta)~(t=1,2)$ and
$\ell:=\ell(g-j_{(s-1)})$. In view of (\ref{rec-tail}) and
(\ref{u-pi_1-f}) and (\ref{u-pi_1-f-1}), there exist polynomials
$f_t^{(s)}(x)$ and $g_t^{(s)}(x)$ $(t=1,2)$
such that, putting
$f_1:=f_1^{(s)}(\theta)$,
$f_2:=f_2^{(s)}(\theta)$,
$g_1:=g_1^{(s)}(\theta)$ and
$g_2:=g_2^{(s)}(\theta)$,
\begin{eqnarray}
u_{s(g-j_s+1)}&=& f_1 u_{s(g-j_{(s-1)})-1}+g_1  u_{s(g-j_{(s-1)})}
=(f_1 x +g_1)\nu_1 x  ^{\ell}+(f_1
y +g_1)\nu_2 y^{\ell}  ;\label{ft-l-eq-1}\\
u_{s(g-j_s+1)-1}& = & f_2 u_{s(g-j_{(s-1)})-1} + g_2
u_{s(g-j_{(s-1)})} = (f_2 x+g_2 )\nu_1 x ^{\ell}+(f_2 y +g_2
)\nu_2 y ^{\ell}\label{ft-l-eq-2}
\end{eqnarray}
both hold for all $\theta \in I$.
Let $x^{'}:=x_{j_s}(\theta)$ and
$y^{'}:=y_{j_s}(\theta)$, and define
$M_{j_s}=M_{j_s}(\kappa,\g,\Delta,L,I)$ by
\begin{equation}\label{M_i}
M_{j_s}(\kappa,\g,\Delta,L,I):=\max\left\{M_{j_{(s-1)}}, \,2
\left| \frac{(f_1y^{'}-f_2)y+g_1y^{'}-g_2} {(f_1y^{'}-f_2)x
+g_1y^{'}-g_2}\right|, \,2 \left|
\frac{(f_1x^{'}-f_2)y+g_1x^{'}-g_2} {(f_1x^{'}-f_2)x
+g_1x^{'}-g_2}\right| \,:\, \theta \in I\right \}.
\end{equation}
By the induction hypothesis, if $\ell(g-j_t)>C_{j_{(s-1)}}$ holds
for all $0\leq t<s-1$, then (\ref{eq-L_i}) holds for the case
$i=j_{(s-1)}$. Moreover, there exists an integer
$E_{j_s}:=E_{j_s}(\kappa,\g,\Delta,L,I)\geq 1$ so that
\begin{equation}\label{mathcalL}
\left| \frac{\nu_1 }{\nu_2 } \right|
\left(\frac{x}{y}\right)^{E_{j_s}}>M_{j_s}
\end{equation}
holds for any $\theta \in I $. Now take
$C^{'}_{j_s}(\kappa,\g,\Delta,L,I):=\max\{E_{j_s}, C_{j_{(s-1)}}\}$
and suppose that $\ell(g-j_t)>C^{'}_{j_s}$ holds for all
$0\leq t < s$. If $\nu_1^{(j_s)}(\eta)=0$ holds for some $\eta\in
I$, then $u_{s(g-j_s+1)-1}(\eta)=u_{s(g-j_s+1)}(\eta)
y_{j_s}(\eta)$ holds and so, as  $\eta\not\in \mathcal{B}$,
\begin{eqnarray*}
&& \bar{\nu}_1 \bar{x}^{\ell}\left( (\bar{f}_1 \bar{y}^{'}-\bar{f}_2)\bar{x}+\bar{g}_1\bar{y}^{'}-\bar{g}_2 \right)
=
\bar{\nu}_2\bar{y}^{\ell} \left( (-\bar{f}_1\bar{y}^{'}+\bar{f}_2)\bar{y}-\bar{g}_1\bar{y}^{'}+\bar{g}_2 \right)\,
\mbox{ and } \\
&& \left( (\bar{f}_1 \bar{y}^{'}-\bar{f}_2)\bar{x}+\bar{g}_1\bar{y}^{'}-\bar{g}_2 \right)\,
\left( (-\bar{f}_1\bar{y}^{'}+\bar{f}_2)\bar{y}-\bar{g}_1\bar{y}^{'}+\bar{g}_2 \right) \neq 0
\end{eqnarray*}
both hold, where $\bar{\nu}_t=\nu_t^{(j_{(s-1)})}(\eta)$,
$\bar{f}_t=f_t^{(s)}(\eta)$, $\bar{g}_t=g_t^{(s)}(\eta)~~(t=1,2)$,
$\bar{x}=x_{j_{(s-1)}}(\eta)$, $\bar{x}^{'}=x_{j_s}(\eta)$,
$\bar{y}=y_{j_{(s-1)}}(\eta)$, and $\bar{y}^{'}=y_{j_{(s-1)}}(\eta)$.
This contradicts (\ref{M_i}) and (\ref{mathcalL}) as
\[M_{j_s}< \left| \frac{\bar{\nu}_1}{\bar{\nu}_2}\right|\left(\frac{\bar{x}}{\bar{y}}\right)^{E_{j_s}}\leq  \left|
\frac{\bar{\nu}_1}{\bar{\nu}_2}\right|\left(\frac{\bar{x}}{\bar{y}}\right)^{\ell} =\left|
\frac{(\bar{f}_1\bar{y}^{'}-\bar{f}_2)\bar{y}+\bar{g}_1\bar{y}^{'}-\bar{g}_2}
{(\bar{f}_1 \bar{y}^{'}-\bar{f}_2)\bar{x}+\bar{g}_1 \bar{y}^{'}-\bar{g}_2}\right|.
\]
Similarly, it follows that
$\nu_t^{(j_s)}(\theta)\neq 0~(t=1,2)$ must hold for all $\theta \in I$.
Moreover, by (\ref{rec-tail}), (\ref{ft-l-eq-1}) and (\ref{M_i}), there exists a
constant $N_{j_s}:=N_{j_s}(\kappa,\mathcal{G},\Delta,L,I)>0$ such
that
\begin{eqnarray}
\left| \frac{\nu_1^{(j_s)}(\theta)}{\nu_2^{(j_s)}(\theta)} \right|
&=& \left|\frac{\nu_1 x^{\ell}\Big((-f_1 y^{'}+f_2)x-g_1y^{'}+g_2
\Big)+\nu_2y^{\ell} \Big((-f_1y^{'}+f_2)y-g_1y^{'}+g_2
\Big)}{\nu_1 x^{\ell } \Big((f_1x^{'}-f_2)x+g_1x^{'}-g_2 \Big)
+\nu_2 y^{\ell} \Big( (f_1x^{'}-f_2)y+g_1x^{'}-g_2 \Big) } \right| \nonumber\\
&\geq &\frac{M_{j_s}\Big|(f_1 y^{'}-f_2)x+g_1y^{'}-g_2 \Big|-
\Big| (f_1y^{'}-f_2)y+g_1y^{'}-g_2\Big|}{M_{j_s}\Big
|(f_1x^{'}-f_2)x+g_1x^{'}-g_2 \Big| + \Big|
(f_1x^{'}-f_2)y+g_1x^{'}-g_2\Big|}\nonumber \\
&\geq & \frac{\Big|(f_1y^{'}-f_2)y+g_1y^{'}-g_2 \Big|}
{M_{j_s}\Big |(f_1x^{'}-f_2)x+g_1x^{'}-g_2 \Big| + \Big|
(f_1x^{'}-f_2)y+g_1x^{'}-g_2 \Big|}\nonumber \\
 &\geq & N_{j_s} \label{c-1}
\end{eqnarray}
holds for any $\theta\in I$. This implies that there exists a
positive constant
$C^{''}_{j_s}=C^{''}_{j_s}(\kappa,\mathcal{G},\Delta,L,I)\geq 1$
such that $N_{j_s}> M_{j_s} P^{C^{''}_{j_s}}$ holds. Hence by
taking
$C_{j_s}=C_{j_s}(\kappa,\g,\Delta,L,I):=\max\{C^{'}_{j_s},C^{''}_{j_s}\}$,
it follows that
\[\left|\frac{\nu_1^{(j_s)}(\theta)}{\nu_2^{(j_s)}(\theta)} \right|\geq N_{j_s}> M_{j_s} P^{C_{j_s}}
\geq M_{j_s}
\left(\frac{y_{j_s}(\theta)}{x_{j_s}(\theta)}\right)^{C_{j_s}}\]
holds for all $\theta\in I$. By applying the
induction hypothesis, it follows that the desired result holds for each
$(\gamma_{g-j_s},\alpha_{g-j_s},\beta_{g-j_s})\in
\underline{\Delta}$ satisfying $0\leq j_s<g-\frak{b}$.
This completes the proof of the theorem. \epf

\section{Well-Placed Intervals}\label{wpi}

In this section, we define the concept of a
well-placed interval with respect to a graphical sequence, and
derive some simple properties of such intervals that will
be used later on. Note that our definition of a
well-placed interval is similar (but not identical) to the one
presented in \cite{rnp}.

Suppose that $\g=\Big( (\gamma_i,\alpha_i,\beta_i)\Big)_{i=1}^{g+1}$
is a $(\kappa,\lambda)$-graphical sequence,
where $\kappa\geq 3$ and $0\leq \lambda\leq \kappa-2$ are integers.


For any closed subinterval $I=[I_{\min}, I_{\max}]$ of
$(\frak{R}_1,\frak{R}_{\max})$ with positive length, define
integers $\frak{a}=\frak{a}(\g, I)$, $\frak{b}=\frak{b}(\g, I)$,
$\frak{c}=\frak{c}(\g, I)$ and $\frak{d}=\frak{d}(\g, I)$ (that
depend only on $\g$ and $I$) by
\begin{eqnarray}
\frak{a}(\g, I)&:=& \min \{ 2\leq i\leq g \,:\, I_{\max }< \frak{R}_i \}, \label{a} \\
\frak{b}(\g, I)&:=& \max \{ 2\leq i\leq g \,:\, I_{\max }< \frak{R}_i \}, \label{b} \\
\frak{c}(\g, I)&:=& \min \Big\{ \{ 2\leq i\leq g \,:\, I_{\max} <\frak{L}_i\}\cup \{g+1\} \Big\},\label{c}\\
\frak{d}(\g, I)&:=& \max \Big\{ \{ 2\leq i\leq g \,:\, I_{\max}
<\frak{L}_i \}\cup \{\frak{c}\} \Big\} . \label{d}
\end{eqnarray}
The interval $I$ is called
a {\em well-placed interval with respect to} $\mathcal{G}$ if
it satisfies the following conditions:
\begin{itemize}
\item[(W1)] $I$ is a closed subinterval of the open interval
$(\frak{R}_1,\frak{R}_{\max} )$ with positive length; \item[(W2)]
If $I\cap I_j\neq \emptyset$ then $I\subseteq I_j$ holds, $1\leq
j\leq g$; \item[(W3)] $ I \subseteq I_{\frak{a}}$, where
$\frak{a}:=\frak{a}(\g,I)$.
\end{itemize}
From now on, we will denote well-placed intervals using
calligraphic script (e.g. $\I$ instead of $I$) to help the reader
follow the text.

In the rest of the section, we will derive some properties of
well-placed intervals. We start with recording some simple
properties of the numbers $\frak{a},\frak{b},\frak{c}, \frak{d}$.

\begin{lemma}\label{h-m-t-partition}
Let $\kappa\geq 3$ and $\lambda\geq 0$ be integers with
$\lambda\leq \kappa-2$, and let
$\g=\Big((\gamma_i,\alpha_i,\beta_i)\Big) _{i=1}^{g+1}$
be a $(\kappa,\lambda)$-graphical sequence. Let $\I=[\I_{\min},
\I_{\max}]$ be a well-placed interval with respect to $\g$. For
the numbers $\frak{a},\frak{b},\frak{c}, \frak{d}$ as defined in
(\ref{a})--(\ref{d}), the following hold:\\
(i) $2\leq \frak{a}\leq \frak{b}\leq g$.\\
(ii) $\frak{c}\leq \frak{d}$. \\
(iii) If $\frak{c}\leq g$, then $2\leq \frak{a}< \frak{c}\leq
\frak{d}\leq \frak{b}\leq g$ holds.\\
(iv) $\{1\leq i\leq g \,:\, 1\leq i<\frak{a} \mbox{  or  }\frak{b}<i\leq g\}\subseteq \{1\leq i\leq g \,:\,
\frak{R}_i<\I_{\min}\}$.\\
(v) If $\frak{c}\leq g$, then $\{1\leq i\leq g \,:\, \frak{a}\leq
i<\frak{c} \mbox{ or }\frak{d}<i\leq \frak{b}\}
\subseteq \{1\leq i\leq g \,:\, \I\subseteq I_i\}$ holds. \\
(vi) If $\frak{c}= g+1$, then $\{1\leq i\leq g \,:\, \frak{a}\leq
i\leq \frak{b}\} = \{1\leq i\leq g \,:\, \I\subseteq I_i\}$ holds.
\end{lemma}

\pf (i)--(iii) are simple consequences of the
definitions of well-placed intervals and the
numbers $\frak{a}$, $\frak{b}$, $\frak{c}$ and $\frak{d}$.

\noindent (iv)--(vi) are direct consequences of the following
inequalities, which follow in view of the fact that the sequence
$(\frak{R}_i)_{i=1}^{g}$ is unimodal by Lemma~\ref{new-unimodular}:
$$
\max \{ \frak{R}_i \,:\, 1\leq i< \frak{a} \mbox{ or } \frak{b}< i
\leq g  \} <\I_{\min}<\I_{\max} < \min \{ \frak{R}_i \,:\,
\frak{a}\leq i \leq \frak{b} \}
$$
and
$$
\max\{\frak{L}_i \,:\, 1\leq i<\min\{\frak{c},~g+1\} \mbox{ or }
\min\{\frak{d},~g+1\} <i\leq g\} <\I_{\min}.
$$
\epf

We now present a result that ensures the existence of
well-placed intervals.

\begin{proposition}\label{existence of I}
Let $\kappa\geq 3$ and $\lambda\geq 0$ be integers with
$\lambda\leq \kappa-2$, and let
$\g=\Big((\gamma_i,\alpha_i,\beta_i)\Big) _{i=1}^{g+1}$
be a $(\kappa,\lambda)$-graphical sequence.
\begin{itemize}
\item[(i)] For each $(\gamma_i,\alpha_i,\beta_i) \in
\underline{\g}\setminus \{(1,\lambda,\kappa-\lambda-1),
(\kappa-\lambda-1,\lambda,1),(\gamma_{g+1},\alpha_{g+1},\beta_{g+1})\}$
and for any closed subinterval $I\subseteq
(\frak{R}_1,\frak{R}_i)$ with positive length, there exists a
well-placed interval $\J_i \subseteq I$ with respect to $\g$ (cf.
(\ref{R-L})). \item[(ii)] Let $\I$ be a well-placed interval with
respect to $\g$. Then any closed interval $\J\subseteq \I$ with
positive length is also a well-placed interval with respect to
$\g$. In particular, $\frak{a}(\g,\J)=\frak{a}(\g,\I)$, $\frak{b}(\g,\J)=\frak{b}(\g,\I)$,
$\frak{c}(\g,\J)=\frak{c}(\g,\I)$, $\frak{d}(\g,\J)=\frak{d}(\g,\I)$ must all hold (cf.
(\ref{a})--(\ref{d})).
\end{itemize}
\end{proposition}
\pf (i): Let $(\gamma_i,\alpha_i,\beta_i) \in
\underline{\g}\setminus \{(1,\lambda,\kappa-\lambda-1),
(\kappa-\lambda-1,\lambda,1),(\gamma_{g+1},\alpha_{g+1},\beta_{g+1})\}$
and suppose that $I=[I_{\min},I_{\max}]\subseteq
(\frak{R}_1,\frak{R}_{i})$ is a subinterval with positive length.
Define
$$
\frak{M}_i:=\max\{I_{\min},\,y \,:\, y\in \{\frak{R}_j,\,\frak{L}_j \,:\,
1\leq j\leq g\}\mbox{  and  } I_{\min}\leq y <I_{\max}\}.
$$
Then $I_{\min}\leq \frak{M}_i<I_{\max}$, and the closed interval
\[ \mathcal{J}_i:=\left[
\frac{I_{\max}+2\,\frak{M}_i}{3},\,\frac{2\,I_{\max}+
\frak{M}_i}{3} \right]\]
is a well-placed interval with respect to $\g$ satisfying
$\J_i\subseteq I$.\\
\noindent (ii): This follows immediately
from the definition of well-placed intervals. \epf

Now, suppose that $\ell: \{1,\dots,g+1\} \to \N$
is a function with $\ell(g+1)=1$.
For $\I$ a well-placed interval with respect
to $\g$, we define $\mathcal{C}=\mathcal{C}_{\g,\I}$, $\Len(\I)=\Len_{\g,\ell}(\I)$
and $\Gap(\I)=\Gap_{\g,\ell}(\I)$ as follows\,:

\begin{eqnarray}
\mathcal{C}_{\g,\I}&:=& \left\{
\begin{array}{ll}
\{1\leq i\leq g \,:\, \frak{a}\leq i<\frak{c} \mbox{~~or~~}\frak{d}<i\leq \frak{b}\} & \mbox{~~if~~}\frak{c}\leq g \\
\{1\leq i\leq g \,:\, \frak{a}\leq i\leq \frak{b}\} &
\mbox{~~if~~}\frak{c}=g+1
\end{array}
\right.\,,\nonumber \\
\Len_{\g,\ell}(\I)&:=& \sum_{j\in \mathcal{C}} \ell(j)~, \label{ell(I)}\\
 \Gap_{\g,\ell}(\I)&:= & \left\{
\begin{array}{ll}
\sum_{\frak{c}\leq j \leq \frak{d}} \ell(j) & \mbox{~~if~~}\frak{c}\leq g\\
0 & \mbox{~~if~~}\frak{c}=g+1
\end{array} \label{def-gap}
\right..
\end{eqnarray}

Using Proposition \ref{existence of I}, we now show that for any
$(\kappa,\lambda)$-graphical sequence $\g$, there is a certain
family of well-placed intervals with respect to  $\g$ each of
whose members avoid the set $\mathcal{B}(\g,\Delta,L)$ as defined
in Definition~\ref{def-bad roots}.

\begin{cor}\label{sub-WI}
Let $\kappa\geq 3$ and $\lambda\geq 0$ be integers with
$\lambda\leq \kappa-2$. Suppose that $(\g,\Delta;L,\ell)$ is a
$(\kappa,\lambda)$-quadruple and let
$\g=\Big(\delta_i:=(\gamma_i,\alpha_i,\beta_i)\Big) _{i=1}^{g+1}$.
Then for any closed subinterval $I\subseteq
(\frak{R}_1,\frak{R}_{i})$ with positive length, there exists a
well-placed interval $\J_i$ in $I$ such that $\J_i\cap
\mathcal{B}= \emptyset$ holds (cf. (\ref{R-L}), (\ref{underline})
and
(\ref{bad roots-def})).\\
In particular, ${\mbox{\em \Len}}(\J_i)\geq \ell(i)$ also holds.
\end{cor}

\pf Suppose that $(\gamma_i,\alpha_i,\beta_i)$ and $I$ are as in
the statement of the corollary. By Proposition~\ref{existence of
I} (i), there exists a well-placed interval $\I \subseteq I$ with
respect to $\g$. By Theorem~\ref{bad-root-bd}, the set
$\mathcal{B}$ is finite. Hence, by Proposition~\ref{existence of I}
(ii), we may take any closed subinterval $\J_i$ of $\I \setminus
\mathcal{B}$ with positive length to give the desired well-placed
interval. \epf

We conclude this section by showing that, in addition, well-placed
intervals satisfying certain other properties also exist.

\begin{prop}\label{gap-ind}
Let $\kappa\geq 3$ and $\lambda\geq 0$ be integers
with $\lambda\leq \kappa-2$. Let
$\g=\Big( (\gamma_i,\alpha_i,\beta_i)\Big)_{i=1}^{g+1}$
be a $(\kappa,\lambda)$-graphical sequence and
$\ell: \{1,\dots,g+1\} \to \N$
is a function with $\ell(g+1)=1$.
Suppose that $\I$ is a well-placed interval
with respect to $\g$ such that ${\mbox{\em \Gap}}(\I)\neq 0$
holds. Then there exists a well-placed interval $\mathcal{J}$
such that \\
(i) $\J_{\min}>\I_{\max}$,\\
(ii) ${\mbox{\em \Gap }}(\J)<{\mbox{\em\Gap}}(\I)$ and\\
(iii) ${\mbox{\em \Len}}(\J)>\frac{{\mbox{\em \Gap}}(\I)}{g}$\\
all hold, where ${\mbox{\em \Len}}(\J):={\mbox{\em \Len}}_{\g,\ell}(\J)$ and
${\mbox{\em \Gap }}(\I):={\mbox{\em \Gap }}_{\g,\ell}(\I)$ are as defined in (\ref{ell(I)}) and (\ref{def-gap}),
respectively.
\end{prop}

\pf As the sequence $(\frak{R}_i)_{i=1}^g$ is unimodal by
Lemma~\ref{new-unimodular} and since
$\Gap (\I)\neq 0$, there exists an integer $j$
with $\frak{c}\leq j\leq \frak{d}$ such
that $\ell(j)>\frac{\Gap(\I)}{g}$
and $\frak{R}_j>\I_{\max}$ both hold,
where $\frak{c}=\frak{c}(\g,\I)$ and
$\frak{d}=\frak{d}(\g,\I)$ are as defined in (\ref{c}) and (\ref{d}).
Hence, by Proposition~\ref{existence of I}, there exists such a
well-placed interval $\J \subseteq (\I_{\max},\frak{R}_j)\subseteq I_j$
as $\Len(\J)\geq \ell(j)>\frac{\Gap (\I)}{g}$. The result now follows.
\epf

\section{Christoffel Numbers}\label{multi-section}

In this section, we prove a result that will allow us to bound the
Christoffel numbers of the $(\kappa, \lambda)$-tridiagonal
sequence associated to a $(\kappa,\lambda)$-quadruple. We will
begin by stating the main theorem of this section,
whose proof will be split into several steps.
To state this result, we require some further definitions.

Let $\g=\Big((\gamma_i,\alpha_i,\beta_i)\Big)_{i=1}^{g+1}$
be a $(\kappa,\lambda)$-graphical sequence for
some integers $\kappa\geq 3$ and $\lambda\geq 0$
with $\lambda\leq \kappa-2$. Let $x$ be a real number.
For each $(\gamma_i,\alpha_i,\beta_i)\in
\underline{\g}\setminus
\{(\gamma_{g+1},\alpha_{g+1},\beta_{g+1})\}$, define
$\rho_i=\rho_i(x)$ and $\sigma_i=\sigma_i(x)$ to be the roots of the
(auxiliary) equation
\begin{equation}\label{def-rho}
\beta_i z^2+(\alpha_i-x)z+\gamma_i=0~,
\end{equation}
which, without loss of generality,
we assume to satisfy $|\rho_i|\geq |\sigma_i|$ for all $1\leq i \leq g$.

\begin{theorem}\label{multi-thm}
Let $\kappa \geq 3$ and $\lambda\geq 0$ be integers with $\lambda
\leq \kappa-2$. Suppose that $(\g,\Delta;L,\ell)$ is a
$(\kappa,\lambda)$-quadruple and let
$\g=\Big((\gamma_i,\alpha_i,\beta_i)\Big)_{i=1}^{g+1}$. Suppose
that $\I$ is a well-placed interval with respect to $\g$
satisfying $\I\cap\mathcal{B}(\g,\Delta,L)=\emptyset$, with
$\mathcal{B}(\g,\Delta,L)$ as defined in Definition~\ref{def-bad
roots}. Then there exist positive constants
$F:=F(\kappa,\g,\Delta,L,\I)~$, $C_1:=C_1(\kappa,\mathcal{G}, \I)$
and $C_2:=C_2(\kappa, \mathcal{G}, \Delta,L,\I)$ so that if
$\ell(i)>F$ holds for all $(\gamma_{i},\alpha_{i},\beta_{i})\in
\underline{\Delta}$ then, for any $\theta\in \I$, the following
holds :
\begin{eqnarray}
C_1\left(\frac{1}{9\kappa^4}\right)^{\mbox{\em
\Gap}(\I)}\,{\mbox{\em \Len}}(\I)\, \prod_{i=1}^{\frak{a}-1} \left(
\left(\frac{\beta_{i}}{\gamma_i} \right) \rho_i^2 \right)^{\ell(i)} &\leq& \sum_{i=0}^{D}\kappa_iu_i^2\nonumber\\
 &\leq & C_2 (9\kappa^4 )^{\mbox{\em \Gap}(\I)}\,{\mbox{\em \Len}}(\I)\,
 \prod_{i=1}^{\frak{a}-1}\left(\left(\frac{\beta_{i}}{\gamma_i}\right)\rho_i^2\right)^{\ell(i)}~,\nonumber
\end{eqnarray}
where $\kappa_i$ and $u_i:=u_i(\theta)$ are as defined in
(\ref{def-k_i}) and (\ref{vi-ui}) for the matrix
$L_1(\T(\g,\ell))$, and $D:=D_{\T(\g,\ell)}$,
$\frak{a}:=\frak{a}(\g,\I)$, ${\mbox{\em \Len}}(\I):={\mbox{\em
\Len}}_{\g,\ell}(\I)$, $\mbox{\em \Gap}(\I):=\mbox{\em
\Gap}_{\g,\ell}(\I)$ and $\rho_i:=\rho_i(\theta)$ are as defined
in (\ref{D(T)}), (\ref{a}), (\ref{ell(I)}), (\ref{def-gap}) and
(\ref{def-rho}), respectively.
\end{theorem}

To prove Theorem~\ref{multi-thm}, we will divide the sum
$\sum_{i=0}^{D}\kappa_iu_i^2$ into three parts: The {\em Head sum}
$\sum_{i=0}^{s(\frak{a})-2}\kappa_iu_i^2$, the {\em Gap sum}
$\sum_{i=s(\frak{a})-1}^{s(\frak{b}+1)}\kappa_iu_i^2$, and the {\em
Tail sum} $\sum_{i=s(\frak{b}+1)+1}^{D}\kappa_iu_i^2$. In
particular, in Section \ref{pre-multi} we will prove a preliminary
result concerning three-term recurrence relations and, for
completeness, recall some additional results on such recursions from
previous papers. We will then use these results in
Section~\ref{head-middle} to derive bounds for the Head and the Gap
sums (as well as to prove some results in Section~\ref{distri}).
Then, in Section~\ref{btail}, we will derive an upper bound for the
Tail sum which, together with the previous bounds, will be used to
prove Theorem~\ref{multi-thm}.

\subsection{Three-Term Recurrence Relations}\label{pre-multi}

Let $\kappa\geq 3$ and $\lambda\geq 0$ be integers with
$\lambda\leq \kappa-2$. Suppose that $\T=\T(\g,\ell)$ is a
$(\kappa,\lambda)$-tridiagonal sequence and let
$\g=\Big((\gamma_i,\alpha_i,\beta_i)\Big)_{i=1}^{g+1}$. Let $x$ be
a real number, and let $\rho_i:=\rho_i(x)$ and
$\sigma_i:=\sigma_i(x)$ be as defined in (\ref{def-rho}), noting
that without loss of generality we are assuming $|\rho_i|\geq
|\sigma_i|$ for all $1\leq i \leq g$. If $x \not \in
\{\frak{R}_i,\, \frak{L}_i \,:\, 1\leq i\leq g\}$, with
$\frak{R}_i$ and $\frak{L}_i$ as defined in (\ref{R-L}), then the
roots $\rho_i$ and $\sigma_i$ are distinct, and so, by standard
theory of recurrence relations, it follows that
\begin{equation}\label{ui-form}
u_{s(i)-1+j}=\omega_1^{(i)}\rho_i^j+\omega_2^{(i)}
\sigma_i^j~~(0\leq j\leq \ell(i)+1)
\end{equation}
holds for some complex numbers $\omega_1^{(i)}:=\omega_1^{(i)}(x)$
and $\omega_2^{(i)}:=\omega_2^{(i)}(x)$, where $u_i = u_i(x)$ are
the numbers associated to the matrix $L_1(\T)$ given by
(\ref{vi-ui}) and $s(i)$ is defined in (\ref{s(i)}). In this situation, note also that (1) if
$|x-\alpha_i|>2\sqrt{\beta_i\gamma_i}$ holds then the roots
$\rho_i$ and $\sigma_i$ are real numbers with
$|\rho_i|>\sqrt{\frac{\gamma_i}{\beta_i}}>|\sigma_i|$, and
$\omega_1^{(i)}, \omega_2^{(i)}$ are real, and (2) if $|x
-\alpha_i|<2\sqrt{\beta_i \gamma_i}$ holds then the roots $\rho_i$
and $\sigma_i$ are complex numbers with
$\sigma_i=\overline{\rho_i}$ and
$|\rho_i|=|\sigma_i|=\sqrt{\frac{\gamma_i}{\beta_i}}$, and
$\omega_1^{(i)}$, $\omega_2^{(i)}$ are complex numbers with
$\omega_2^{(i)}=\overline{\omega_1^{(i)}}$.

We now prove a result that is analogous with the result
\cite[Proposition~3.1]{bkm-length} that was proven to hold for
distance-regular graphs.

\begin{proposition}\label{new-ell-3.1}{\ \\}
Let $\kappa\geq 3$ and $\lambda\geq 0$ be integers with
$\lambda\leq \kappa-2$. Let
$\g=\Big((\gamma_i,\alpha_i,\beta_i)\Big)_{i=1}^{g+1}$ be a
$(\kappa,\lambda)$-graphical sequence, $\ell:\{1,\dots,g+1\} \to
\N$ be a function with $\ell(g+1)=1$, and $\T=\T(\g,\ell)$ be the
$(\kappa,\lambda)$-tridiagonal sequence associated to $\g$ and
$\ell$. Suppose that $\I$ is a well-placed interval with respect
to $\g$. Then, for all $\theta\in I$ the following hold (cf. (\ref{a}), (\ref{ui-form})):\\
(i) $0<\sigma_i(\theta)<\rho_i(\theta)<1$, for all $1\leq i\leq \frak{a}-1$.\\
(ii) $u_{s(i)-1}(\theta)>\prod_{j=1}^{i-1}\rho_j(\theta)^{\ell(j)}$,
for all $2\leq i\leq \frak{a}$.\\
(iii) $-\omega_1^{(i)}(\theta)<\omega_2^{(i)}(\theta)<0<\omega_1^{(i)}(\theta)$,
for all $1\leq i\leq \frak{a}$.
\end{proposition}
\pf Suppose $\theta \in \I$, and put $\rho_i:=\rho_i(\theta)$,
$\sigma_i:=\sigma_i(\theta)$, $(1\leq i\leq g)$,
$u_j:=u_j(\theta)$ $(1\leq j\leq D_{\T})$ and
$\omega_j^{(i)}:=\omega_j^{(i)}(\theta)~(j=1,2)$ as in
(\ref{ui-form}).

\noindent (i): Since $\theta>\frak{R}_i$ holds for all $1\leq
i<\frak{a}$ by (\ref{a}), $0<\sigma_i<\rho_i$ holds for all $1\leq
i< \frak{a}$. By (G0) in Definition \ref{GS} and by
Lemma~\ref{new-unimodular}, $0<\theta<\frak{R}_{\max} =\max
\{\kappa-(\sqrt{\beta_i}-\sqrt{\gamma_i})^2\,:\, 1\leq i\leq g\}
\leq \kappa$ and $\beta_i>\gamma_i~(1\leq i< \frak{a})$ both hold.
Hence
\begin{eqnarray*}
&& 2\beta_i-(\theta-\alpha_i)=(\kappa -\theta)+(\beta_i-\gamma_i)>0\mbox{~~~and~}\\
&& (2\beta_i-(\theta-\alpha_i))^2-((\theta-\alpha_i)^2-4\beta_i\gamma_i)=4\beta_i (\kappa-\theta)>0
\end{eqnarray*}
follow. Thus, (i) holds by (\ref{def-rho}) and
the fact that
$\rho_i=\frac{\theta-\alpha_i+\sqrt{(\theta-\alpha_i)^2-4\beta_i\gamma_i}}{2\beta_i}$
holds.\\

To prove that (ii) and (iii) hold, we will use the following
claim.

\begin{claim}\label{claim3.1}
(a) $\rho_{i+1}<\rho_i~~(1\leq i\leq \frak{a}-1)$.\\
(b) $u_{s(i)}>\rho_i\, u_{s(i)-1}~(1\leq i\leq \frak{a})$.\\
(c) $\omega_1^{(i)}>u_{s(i)-1}~(1\leq i\leq \frak{a})$.
\end{claim}
\noindent{\bf Proof of Claim \ref{claim3.1}}:
In view of Proposition~\ref{new-ell-3.1} (i)
and $\beta_j\rho_j^2+(\alpha_j-\theta)\rho_j+\gamma_j=0~(1\leq j\leq g)$, it
follows that
\begin{eqnarray*}
&&(\beta_i-1)\rho_i^2+(\alpha_i+1-\theta) \rho_i+\gamma_i
=\rho_i (1-\rho_i)>0\mbox{~~and} \\
&& \beta_i \rho_i^2+(\alpha_i-1-\theta) \rho_i+(\gamma_i+1)=1-\rho_i>0
\end{eqnarray*}
for all $1\leq i\leq \frak{a}-1$. Hence, by (G2) in Definition
\ref{GS}, statement (a) in the claim holds.

We now prove statements (b) and (c) by using induction on $i$.
Suppose $i=1$. By
\[
\beta_1 \left(\frac{\theta}{\kappa}\right)^2+(\lambda-\theta)
\left( \frac{\theta}{\kappa}\right) +1
=\left(1-\frac{\theta}{\kappa}\right)
\left(1+(1+\lambda)\left(\frac{\theta}{\kappa}\right) \right)>0~,
\]
$\rho_1u_0=\rho_1<\frac{\theta}{\kappa}=u_1$ hold. Thus, by
Proposition~\ref{new-ell-3.1} (i) and
$\rho_1u_0=\rho_1(\omega_1^{(1)}+\omega_2^{(1)})
<\omega_1^{(1)}\rho_1+\omega_2^{(1)}\sigma_1=u_1$, (b) and (c)
hold for $i=1$.

Now let $2\leq i< \frak{a}$, and suppose that (b) and (c) hold for
all $2\leq j\leq i$. By the induction hypothesis, it follows
$u_{s(i+1)}-\rho_i u_{s(i+1)-1}
=\omega_2^{(i)}\sigma_i^{\ell(i)}(\sigma_i-\rho_i)>0$ by (i) of
the proposition, and $u_{s(i+1)}>\rho_i u_{s(i+1)-1}>\rho_{i+1}
u_{s(i+1)-1}$ by statement (a) of the claim. Thus, by
(\ref{ui-form}), (b) and (c) hold for all $1\leq i\leq \frak{a}$,
which completes the proof of the claim. \epf

\noindent (ii): We prove this using induction on $i$. Suppose
$i=2$. Then by applying (b) and (c) of Claim~\ref{claim3.1} and
statement (i) (with $i=1$),
\[
u_{s(2)-1}-\rho_1^{\ell(1)}=
\rho_1^{\ell(1)}(\omega_1^{(1)}-1)+\omega_2^{(1)}\sigma_1^{\ell(1)}
>\rho_1(\omega_1^{(1)}-1)(\rho_1^{\ell(1)-1}-\sigma_1^{\ell(1)-1})>0.
\]
Therefore (ii) holds for $i=2$.

Now let $2\leq i<\frak{a}$, and suppose that (ii) holds for all
$2\leq j\leq i$. Using (i) and Claim~\ref{claim3.1} (c), it
follows that
\[
u_{s(i+1)-1}-u_{s(i)-1}\rho_i^{\ell(i)}
=\omega_2^{(i)}(\sigma_i^{\ell(i)}-\rho_i^{\ell(i)})>0~
\]
holds. Hence, by induction, $u_{s(i+1)-1}>\prod_{j=1}^{i}\rho_j^{\ell(j)}$.

\noindent (iii): Using (ii) and Claim~\ref{claim3.1} (c), it
follows
$0<u_{s(i)-1}=\omega_1^{(i)}+\omega_2^{(i)}<\omega_1^{(i)}$ for
all $1\leq i\leq \frak{a}$. Now, (iii) follows immediately.\epf

We now recall a result that was originally stated using
different terminology in \cite{two-thm} and \cite{ban-ito-II}.

\begin{lemma}(cf. \cite[Lemma 5.1]{two-thm},\,
\cite[Proposition~7]{ban-ito-II})\label{bounder}{\ \\} Let
$\kappa\geq 3$ and $\lambda\geq 0$ be integers with $\lambda\leq
\kappa-2$. Let
$\g=\Big((\gamma_i,\alpha_i,\beta_i)\Big)_{i=1}^{g+1}$ be a
$(\kappa,\lambda)$-graphical sequence, $\ell:\{1,\dots,g+1\} \to
\N$ be a function with $\ell(g+1)=1$, and $\T:=\T(\g,\ell)$ be the
$(\kappa,\lambda)$-tridiagonal sequence associated to $\g$ and
$\ell$ with diameter $D_{\T}$ (cf. (\ref{D(T)})). Let $\theta$ be
any real number with $|\theta|\leq \kappa$.
Then for each $i=1,\ldots, D_{\T}-1$,\\
(i)
$$
\frac{1}{3\kappa} \max\{ |u_i(\theta)|, |u_{i+1}(\theta)|\}
\leq \max\{ |u_{i-1}(\theta)|,|u_{i}(\theta)|\}
\leq 3\kappa\max\{ |u_i(\theta)|, |u_{i+1}(\theta)|\}
$$
and\\
(ii)
$$
\left(\frac{1}{9\kappa^4}\right)
\max\{\kappa_{i-1}u_{i-1}^2(\theta), \kappa_iu_i^2(\theta)\} \leq
\max\{\kappa_{i}u_{i}^2(\theta), \kappa_{i+1}u_{i+1}^2(\theta) \}
\leq 9\kappa^4 \max\{\kappa_{i-1}u_{i-1}^2(\theta),
\kappa_iu_i^2(\theta)\}\,
$$
hold, where $\kappa_i$ and $u_i(\theta)$ are as defined in
(\ref{def-k_i}) and (\ref{vi-ui}) for the matrix $L_1(\T)$.
\end{lemma}
\pf
(i): Since $|\theta|\leq \kappa$ and
$0<\beta_i,\,\gamma_i<\kappa~(1\leq i\leq g)$ hold,
it follows by (\ref{u_i-three-terms}) that
\[
|u_{i+1}(\theta)|=\left|
\left(\frac{\theta-\alpha_i}{\beta_i}\right)u_i(\theta)-
\left(\frac{\gamma_i}{\beta_i}\right)u_{i-1}(\theta)\right| \leq
2\,\kappa\,|u_i(\theta)|+\kappa\,|u_{i-1}(\theta)| \leq
3\kappa\,\max\{|u_{i-1}(\theta)|,\,|u_{i}(\theta)|\}
\]
and
\[
|u_{i-1}(\theta)|=\left|
\left(\frac{\theta-\alpha_i}{\gamma_i}\right)u_i(\theta)-
\left(\frac{\beta_i}{\gamma_i}\right)u_{i+1}(\theta)\right|
\leq 2\,\kappa\,|u_i(\theta)|+\kappa\,|u_{i+1}(\theta)|
\leq 3\kappa\,\max\{|u_i(\theta)|,\,|u_{i+1}(\theta)|\}
\]
all hold. Statement (i) now follows immediately.

\noindent (ii): Since $\frac{1}{\kappa}\kappa_{i+1} \leq
\kappa_i\leq \kappa\,\kappa_{i-1},~~i=1,\ldots, D_{\T}-1,$ holds
by (\ref{kappa_i}), statement (ii) follows immediately from (i).
\epf

For completeness, we now recall two results from \cite{two-thm}.

\begin{cor}\label{4.2-two-thm}(\cite[Corollary 4.2]{two-thm}){\ \\}
Suppose $N\geq 2$ is an integer,
and $\alpha\geq 0$, $\beta>0$, $\gamma>0$, $x_0$
and $x_1$ are real numbers satisfying $(x_0,x_1)\neq (0,0)$.
Let $\epsilon$ be a real number with $0<\epsilon <2\sqrt{\beta \gamma}$.
Then there exist positive real numbers
$C_s:=C_s(\beta,\gamma,\epsilon),~s=1,2,3,4$ such
that for every real number $\theta$
with $|\theta -\alpha|\leq 2\sqrt{\beta \gamma}-\epsilon$,
and for all real numbers $x_2,\ldots,x_N$
satisfying $\gamma x_{i-1}+(\alpha-\theta)x_i+\beta x_{i+1}=0~~(i=1,\ldots, N-1)$,
we have
\[C_1 \max\left\{x_0^2,\left(
\frac{\beta}{\gamma}\right) x_1^2\right\}\leq
\max\left\{\left(\frac{\beta}{\gamma}\right)^{i-1}x_{i-1}^2,\left(
\frac{\beta}{\gamma}\right)^i x_i^2\right\}\leq
C_2\max\left\{x_0^2,\left( \frac{\beta}{\gamma}\right)
x_1^2\right\}\] for $i=1,2,\ldots, N$, and
\[ C_3 N \max\left\{x_0^2,\left(
\frac{\beta}{\gamma}\right) x_1^2\right\}\leq \sum_{i=0}^{N}
\left(\frac{\beta}{\gamma}\right)^i x_i^2 \leq C_4 N \max
\left\{x_0^2,\left( \frac{\beta}{\gamma}\right) x_1^2\right\}.\]
\end{cor}

\begin{proposition}\label{4.3-two-thm}(\cite[Proposition 4.3]{two-thm}){\ \\}
Suppose $N\geq 2$ is an integer, and $\alpha\geq 0$, $\beta>0$,
$\gamma>0$, $x_0$ and $x_1$ are real numbers satisfying
$(x_0,x_1)\neq (0,0)$. Let $\kappa,\,\epsilon$ and $\epsilon'$ be
positive real numbers. Then there exist constants
$C_1=C_1(\kappa,\alpha,\beta,\gamma,\epsilon)>0$ and
$C_2=C_2(\beta,\gamma,\epsilon)>1$ such that, for every real
number $\theta$ with $|\theta-\alpha|\geq 2\sqrt{\beta
\gamma}+\epsilon$, $|\theta|\leq \kappa$, and
\[|x_1-x_0 \sigma|>\epsilon' \,\max \left\{|x_0|, \sqrt{\frac{\beta}{\gamma}}|x_1|\right\}\]
(with $\rho=\rho(\theta)$ and $\sigma=\sigma(\theta)$
the roots of $\beta x^2+(\alpha-\theta)x+\gamma=0$
with $|\rho|\geq |\sigma|$), and for all real numbers
$x_2,\ldots, x_N$ satisfying
$\gamma x_{i-1}+(\alpha-\theta)x_i+\beta x_{i+1}=0~~(i=1,\ldots, N-1)$, we have
\[\sum_{i=0}^{N}\left(\frac{\beta}{\gamma}\right)^i x_i^2\leq C_1
\left( \left(\frac{\beta}{\gamma}\right)\rho^2\right)^N
\max\left\{x_0^2,\left( \frac{\beta}{\gamma}\right)
x_1^2\right\}\] and, for all $n\leq N$,
\[x_n^2 \leq C_2\,
\rho^{2n}\,\max \left\{x_0^2,\left( \frac{\beta}{\gamma}\right)
x_1^2\right\}.\]
\end{proposition}

\subsection{Bounding Head and Gap Sums}\label{head-middle}

In this subsection, we obtain bounds for Head sum
$\sum_{i=0}^{s(\frak{a})-2} \kappa_i u_i^2 $ and Gap sum
$\sum_{i=s(\frak{a})-1}^{s(\frak{b}+1)} \kappa_i u_i^2 $. In more
detail, we prove the following:

\begin{theorem}\label{0-s-sum}
Let $\kappa\geq 3$ and $\lambda\geq 0$ be integers with
$\lambda\leq \kappa-2$. Let
$\g=\Big((\gamma_i,\alpha_i,\beta_i)\Big)_{i=1}^{g+1}$ be a
$(\kappa,\lambda)$-graphical sequence, $\ell:\{1,\dots,g+1\} \to
\N$ be a function with $\ell(g+1)=1$, and $\T:=\T(\g,\ell)$ be the
$(\kappa,\lambda)$-tridiagonal sequence associated to $\g$ and
$\ell$. Suppose that $\I$ is a well-placed interval with respect
to $\g$. Then there exist positive constants
$C_i:=C_i(\kappa,\mathcal{G},\I)~(1\leq i\leq 11)$ such that for
any element $\theta$ in $\I$, the following all hold:
\begin{enumerate}
\item[(i)] $\sum_{i=0}^{s(\frak{a})-2} \kappa_i u_i^2 \leq C_1
\prod_{i=1}^{\frak{a}-1}\left ( \left( \frac{\beta_i}{\gamma_i}
\right) \rho_i^2 \right
    )^{ \ell(i)
    }$.
\item[(ii)] $\prod_{i=1}^{\frak{a}-1}\rho_i^{2 \ell(i)} < \max
\left\{u_{s(\frak{a})-1}^2, \left(
    \frac{\beta_{\frak{a}}}{\gamma_{\frak{a}}}\right)u_{s(\frak{a})}^2 \right\} \leq C_2
    \prod_{i=1}^{\frak{a}-1} \rho_i^{ 2 \ell(i)}$.

\item[(iii)] Let $\widehat{\frak{c}}:=\min\{\frak{c},\frak{b}+1\}$. Then
\[ C_3 \left(\sum_{i=\frak{a}}^{\widehat{\frak{c}}-1}\ell(i)\right)
    \prod_{j=1}^{\frak{a}-1}\left(\left(\frac{\beta_j}{\gamma_j}\right)\rho_j^2\right
)^{\ell(j)}
 \leq  \sum_{i=s(\frak{a})-1}^{s(\widehat{\frak{c}})}\kappa_iu_i^2  \leq C_4
\left(\sum_{i=\frak{a}}^{\widehat{\frak{c}}-1}\ell(i)\right)\prod_{j=1}^{\frak{a}-1}\left
(\left(\frac{\beta_j}{\gamma_j}\right) \rho_j^2\right )^{\ell(j)}.
\]
\item[(iv)] Let $\widehat{\frak{c}}:=\min\{\frak{c},\frak{b}+1\}$.
Then \[C_5 \prod_{i=1}^{\frak{a}-1}\rho_i^{2 \ell(i)} \leq
\max\left\{ u_{s(\widehat{\frak{c}})-1}^2,
\left(\frac{\beta_{\widehat{\frak{c}}-1}}{\gamma_{\widehat{\frak{c}}-1}}\right)
u_{s(\widehat{\frak{c}})}^2\right\}
\prod_{i=\frak{a}}^{\widehat{\frak{c}}-1}\left(\frac{\beta_i}{\gamma_i}\right)^{\ell(i)}
\leq C_6 \prod_{i=1}^{\frak{a}-1}\rho_i^{2 \ell(i)}.\]

\item[(v)]
If $\frak{c}\leq g$, then
$\sum_{i=s(\frak{c})}^{s(\frak{d}+1)}\kappa_iu_i^2< C_7
\left(9\kappa^4\right)^{\mbox{\em \Gap }(\I)}\prod_{i=1}^{\frak{a}-1}\left(\left(\frac{\beta_i}{\gamma_i}\right)
\rho_i^2\right )^{\ell(i)}.$

\item[(vi)] If $\frak{c}\leq g$ and $\frak{d}<\frak{b}$ both hold,
then
\begin{eqnarray}
&& C_8 \left(\frac{1}{9\kappa^4}\right)^{\mbox{\em \Gap} (\I)
}~\left(\sum_{i=\frak{d}+1}^{\frak{b}}\ell(i)\right)
\prod_{j=1}^{\frak{a}-1}\left(\left(\frac{\beta_j}{\gamma_j}\right)\rho_j^2\right
)^{\ell(j)}\nonumber \\ & \leq &
\sum_{i=s(\frak{d}+1)-1}^{s(\frak{b}+1)}\kappa_i u_i^2 \leq C_9 (9
\kappa ^4)^{\mbox{\em \Gap}(\I) }~
\left(\sum_{i=\frak{d}+1}^{\frak{b}}\ell(i)\right)\prod_{j=1}^{\frak{a}-1}
\left(\left(\frac{\beta_j}{\gamma_j}\right)\rho_j^2\right)^{\ell(j)}.\nonumber
\end{eqnarray}

\item[(vii)] If $\frak{c}\leq g$, then
\begin{eqnarray}
C_{10} \left(\frac{1}{9\kappa^4}\right)^{\mbox{\em \Gap}(\I)}~\prod_{i=1}^{\frak{a}-1}\left(
    \left(\frac{\beta_i}{\gamma_i}\right)\rho_i^2\right )^{\ell(i)} &\leq & \kappa_{s(\frak{b}+1)}
    \max\left\{u_{s(\frak{b}+1)-1}^2,
    \left(\frac{\beta_{\frak{b}}}{\gamma_{\frak{b}}}\right)u_{s(\frak{b}+1)}^2 \right\}\nonumber\\ &\leq & C_{11}
    (9\kappa^4)^{\mbox{\em \Gap}(\I)}~\prod_{i=1}^{\frak{a}-1}\left(
    \left(\frac{\beta_i}{\gamma_i}\right)\rho_i^2
    \right)^{\ell(i)}~,\nonumber
\end{eqnarray}
\end{enumerate}
where $\kappa_i$ and $u_i:=u_i(\theta)$ are as defined in
(\ref{def-k_i}) and (\ref{vi-ui}) relative to the matrix
$L_1(\T)$, and $s(i)$, $\frak{a}$, $\frak{b}$, $\frak{c}$,
$\frak{d}$, $\mbox{\em \Gap}(\I):=\mbox{\em \Gap}_{\g,\ell}(\I)$
and $\rho_i:=\rho_i(\theta)$ are as defined in (\ref{s(i)}),
(\ref{a})--(\ref{d}), (\ref{def-gap}) and (\ref{def-rho}),
respectively.
\end{theorem}

\pf Suppose that $\mathcal{G}$, $\ell$, $\T$, $\I$ and $\theta$
are as in the statement of the theorem.

\noindent{(i) and (ii):} In order to apply
Proposition~\ref{4.3-two-thm}, we first prove that there are positive constants
$\epsilon_1:=\epsilon_1(\mathcal{G},\I)$ and
$\epsilon_2:=\epsilon_2(\mathcal{G},\I)$ such that, for
all $1 \le i \le \frak{a}-1$,\\
(a) $|\theta -\alpha_i|\ge 2\sqrt{\beta_i\gamma_i}+\epsilon_1$,  and\\
(b) $|u_{s(i)}-\sigma_iu_{s(i)-1}|> \epsilon_2
\max\left\{|u_{s(i)-1}|,\sqrt{\frac{\beta_i}{\gamma_i}}|u_{s(i)}| \right\}$\\
both hold.

For statement (a), we can take
$\epsilon_1=\epsilon_1(\mathcal{G},\I):=
\min\{\I_{\min}-(\alpha_i+2\sqrt{\beta_i\gamma_i})\,:\, 1\le i\le
\frak{a}-1\}$, in view of (\ref{a}) and (W2).

By (\ref{a}) and (\ref{def-rho}), inequalities $M_1 \geq
\rho_{i}>0$ and
$\frac{\rho_i-\sigma_i}{\sqrt{\beta_i/\gamma_i}}\ge M_2>0 $ all
hold for any $\theta\in \I$ and for any $1 \le i \le \frak{a}-1$,
where
\begin{eqnarray*}
M_1=M_1(\mathcal{G},\I)&:=&
\max\left \{\frac{(\I_{\max}-\alpha_i)+ \sqrt{(\I_{\max}-\alpha_i)^2-4\beta_i\gamma_i}}{2\beta_i}\,:\, 1\le i\le
\frak{a}-1
\right \}~\mbox{~and~}\\
M_2=M_2(\mathcal{G},\I)&:=& \min\left\{\frac{\sqrt{\gamma_i}
 \sqrt{(\I_{\min}-\alpha_i)^2-4\beta_i\gamma_i}}{\beta_i\sqrt{\beta_i}} \,:\,
1\le i\le \frak{a}-1 \right\}.
\end{eqnarray*}
By (i) and (iii) of Proposition~\ref{new-ell-3.1},
\[
\frac{|u_{s(i)}-\sigma_{i}u_{s(i)-1}|}
{\max\left\{|u_{s(i)-1}|,\sqrt{\frac{\beta_{i}}{\gamma_{i}}}|u_{s(i)}|
\right\}}
=\left|\frac{\omega_1^{(i)}}{\omega_1^{(i)}\rho_{i}+\omega_2^{(i)}\sigma_{i}}\right|
~\frac{|\rho_{i}-\sigma_{i}|}{\sqrt{\frac{\beta_{i}}{\gamma_{i}}}}
> \frac{1}{|\rho_{i}|}\frac{|\rho_{i}-\sigma_{i}|}
{\sqrt{\frac{\beta_{i}}{\gamma_{i}}}}\ge \frac{M_2}{M_1} > 0
\]
holds, and hence (b) holds for
$\epsilon_2=\epsilon_2(\mathcal{G},\I)=\frac{M_2}{M_1}$.\\
Now by Proposition~\ref{4.3-two-thm}, there exist constants
$M_3=M_3(\kappa,\mathcal{G},\I)>1$ and
$M_4=M_4(\kappa,\mathcal{G},\I)>0$ such that for all $1 \le i\le
\frak{a}-1$,
\begin{equation}\label{4.3-i}
\max\left\{u_{s(i+1)-1}^2,
\left(\frac{\beta_{i+1}}{\gamma_{i+1}}\right)u_{s(i+1)}^2
\right\}\le M_3\,\rho_{i}^{2
\ell(i)}\max\left\{u_{s(i)-1}^2,\left(\frac{\beta_{i}}{\gamma_{i}}\right)
u_{s(i)}^2\right\}
\end{equation}
and
\begin{equation}\label{4.3-ii}
\sum_{j=0}^{\ell(i)+1}\kappa_{s(i)-1+j}u_{s(i)-1+j}^2\le  M_4\,
\kappa_{s(i)-1}\left(\left(\frac{\beta_{i}}{\gamma_{i}}\right)\rho_{i}^2\right)^{\ell(i)}
\max\left\{u_{s(i)-1}^2,
\left(\frac{\beta_{i}}{\gamma_{i}}\right)u_{s(i)}^2 \right\}
\end{equation}
both hold. By applying (\ref{4.3-i}) inductively and
also using (\ref{4.3-ii}), it follows
that, for each $1\leq i\leq \frak{a}-1$,
\begin{eqnarray*}
\max\left\{u_{s(i+1)-1}^2,
\left(\frac{\beta_{i+1}}{\gamma_{i+1}}\right)u_{s(i+1)}^2
\right\}&\le & M_3^{i}\prod_{j=1}^{i}\rho_{j}^{2 \ell(j)
}\max\left\{u_0^2,
\left(\frac{\beta_{1}}{\gamma_{1}}\right)u_1^2\right\}< \kappa\, M_3^{i}\prod_{j=1}^{i}\rho_{j}^{2 \ell(j) }\,
\mbox{ and }\\
\sum_{j=0}^{\ell(i)+1}\kappa_{s(i)-1+j}u_{s(i)-1+j}^2&\le& \kappa^2\, M_4 \,M_3^{i-1}\,\prod_{j=1}^{i}
\left(\left(\frac{\beta_{j}}{\gamma_{j}}\right)
\rho_{j}^2\right)^{ \ell(j) }\le  \kappa^2\, \,M_4\,M_3^{\frak{a}-2} \,\prod_{j=1}^{\frak{a}-1}
\left(\left(\frac{\beta_j}{\gamma_j}\right)\rho_j^2\right)^{
\ell(j) }~
\end{eqnarray*}
all hold. Statements (i) and (ii) now follow by taking
$$
C_1:=(\frak{a}-1) \,\kappa^2\,M_4\,M_3^{\frak{a}-2}~,~~C_2:=\kappa \,M_3^{\frak{a}-1}
$$
and noting that $u^2_{s(\frak{a})-1}>\prod_{j=1}^{\frak{a}-1}\rho_{j}^{2\ell(j)}$
holds by Proposition~\ref{new-ell-3.1} (ii).\\

\noindent{(iii):} By Lemma~\ref{h-m-t-partition}, $\I\subseteq
I_i$ holds for each $\frak{a}\leq i\leq \widehat{\frak{c}}-1$. Let
$ \epsilon=\epsilon(\mathcal{G},\I)
:=\min\{|\alpha_i+2\sqrt{\beta_i\gamma_i}- \I_{\max}|,
|\I_{\min}-(\alpha_i -2\sqrt{\beta_i\gamma_i})| \,:\, \frak{a}\le
i\le \widehat{\frak{c}}-1\} $. Then $|\theta -\alpha_i|\le
2\sqrt{\beta_i\gamma_i}-\epsilon$ and $0<\epsilon
<2\sqrt{\beta_i\gamma_i}$ both hold for all $\frak{a} \le i \le
\widehat{\frak{c}}-1$. Hence by Corollary~\ref{4.2-two-thm}, there
exist constants $M_j:=M_j(\mathcal{G},\I)>0~(5\leq j \leq 8)$ such
that, for any $\frak{a} \le i \le \widehat{\frak{c}}-1$,
\begin{equation}\label{4.2-i}
M_5\max\left\{u_{s(i)-1}^2,\left(\frac{\beta_{i}}{\gamma_{i}}\right)u_{s(i)}^2
\right\}\le \left(\frac{\beta_{i}}{\gamma_{i}}\right)^{\ell(i)}
\max\left\{u_{s(i+1)-1}^2,
\left(\frac{\beta_{i}}{\gamma_{i}}\right)u_{s(i+1)}^2 \right \}
\le M_6
\max\left\{u_{s(i)-1}^2,\left(\frac{\beta_{i}}{\gamma_{i}}\right)u_{s(i)}^2
\right\}
\end{equation}
and
\begin{equation}\label{4.2-ii}
M_7~
\ell(i)\max\left\{u_{s(i)-1}^2,\left(\frac{\beta_{i}}{\gamma_{i}}\right)u_{s(i)}^2
\right\} \le
\sum_{j=0}^{\ell(i)+1}\left(\frac{\beta_{i}}{\gamma_{i}}\right)^ju_{s(i)-1+j}^2
\le M_8~
\ell(i)\max\left\{u_{s(i)-1}^2,\left(\frac{\beta_{i}}{\gamma_{i}}\right)u_{s(i)}^2\right\}
\end{equation}
both hold. Note that for each $1\leq i\leq g$, the following hold for all $0 \le j \le \ell(i)+1$:
\begin{equation}\label{val-ineq}
\frac{1}{\kappa} \left(\frac{\beta_i}{\gamma_i}\right)^j \prod_{m=1}^{i-1}\left(\frac{\beta_m}{\gamma_m}\right)^{\ell(m)}
\leq \kappa_{s(i)-1+j}\leq \kappa
\left(\frac{\beta_i}{\gamma_i}\right)^j \prod_{m=1}^{i-1}\left(\frac{\beta_m}{\gamma_m}\right)^{\ell(m)}.
\end{equation}

Hence, by applying (\ref{val-ineq}) to (\ref{4.2-ii}) and by using
(\ref{4.2-i}) and statement (ii) of the theorem, it follows that,
for each $\frak{a}\leq i\leq \widehat{\frak{c}}-1$, there exists a
constant $M_9:=M_9(\kappa,\mathcal{G},\I)>0$ such that
\begin{eqnarray*}\label{mid-upper}
\sum_{j=0}^{\ell(i)+1}\kappa_{s(i)-1+j}u_{s(i)-1+j}^2 &\le &
\kappa~ M_8 \, \ell(i)\,
\prod_{j=1}^{i-1}\left(\frac{\beta_{j}}{\gamma_{j}}\right)^{\ell(j)}~\max\left\{u_{s(i)-1}^2,
\left(\frac{\beta_{i}}{\gamma_{i}}\right)u_{s(i)}^2\right\}\nonumber
\\ &\le & \kappa M_8 M_6^{i-\frak{a}}\,\ell(i)\, \prod_{j=1}^{\frak{a}-1}\left(\frac{\beta_j}{\gamma_j}\right)^{\ell(j)}
\max\left\{u_{s(\frak{a})-1}^2,
\left(\frac{\beta_{\frak{a}}}{\gamma_{\frak{a}}}\right)u_{s(\frak{a})}^2 \right\}\\
&\le & \kappa \max\{1,M_6^{\widehat{\frak{c}}}\}M_8 M_9
\,\ell(i)\,
\prod_{j=1}^{\frak{a}-1}\left(\left(\frac{\beta_j}{\gamma_j}\right)\rho_j^2\right)^{\ell(j)}
\nonumber
\end{eqnarray*}
holds and, similarly,
\[ \sum_{j=0}^{\ell(i)+1}\kappa_{s(i)-1+j}u_{s(i)-1+j}^2 \geq  \frac{M_7}{\kappa}\,\min\{1, M_5^{\widehat{\frak{c}}} \}\,
\ell(i)\,
\prod_{j=1}^{\frak{a}-1}\left(\left(\frac{\beta_j}{\gamma_j}\right)\rho_j^2\right)^{\ell(j)}.
\]
Hence (iii) follows by taking
$$
C_3(\kappa,\mathcal{G},\I):=\frac{M_7}{3 \kappa} \min \{1, M_5^{\widehat{\frak{c}}}
\},~C_4(\kappa,\g,\I):=\kappa \max\{1,M_6^{\widehat{\frak{c}}} \} M_8 M_9,
$$
in light of that fact that each element $\kappa_{s(i)-1+j}u_{s(i)-1+j}^2$
appears in the sum
$\sum_{i=\frak{a}}^{\widehat{\frak{c}}-1}
\sum_{m=0}^{\ell(i)+1}\kappa_{s(i)-1+m}u_{s(i)-1+m}^2$ at most three times.\\

\noindent{(iv):} This follows from (\ref{4.2-i}) and statement (ii) of the theorem. \\

\noindent{(v):} By Lemma~\ref{bounder} (ii), statement (iv) of the
theorem, (\ref{def-gap}) and (\ref{val-ineq}), there exists a
constant $C_7=C_7(\kappa,\mathcal{G},\I)>0$ so that
\begin{eqnarray*}
\sum_{i=s(\frak{c})}^{s(\frak{d}+1)}\kappa_iu_i^2&<& 2 \left( 9
\kappa^4\right)^{s(\frak{d}+1)-s(\frak{c})+1}\max\{\kappa_{s(\frak{c})-1}u_{s(\frak{c})-1}^2,
\kappa_{s(\frak{c})}u_{s(\frak{c})}^2\}\\ &\leq & C_7 \left( 9
\kappa^4\right)^{\Gap(\I)}
\prod_{i=1}^{\frak{a}-1}\left(\left(\frac{\beta_i}{\gamma_i}\right)\rho_i
^2\right)^{\ell(i)}
\end{eqnarray*}
holds. (v) follows immediately.\\

\noindent{(vi) and (vii):} Using the same proof as for statement
(iii), it can be seen that if $\frak{d}<\frak{b}$ then
(\ref{4.2-i}) and (\ref{4.2-ii}) both hold for all $\frak{d}+1\leq
i\leq \frak{b}$.\\
By (\ref{4.2-i})--(\ref{val-ineq}), Lemma~\ref{bounder} (ii) and
statement (iv) of the theorem, there exist constants
$M_j=M_j(\kappa, \g, \I)>0~(10\leq j\leq 15)$ such that
\begin{eqnarray*}
\kappa_{s(\frak{b}+1)}\max\left\{u_{s(\frak{b}+1)-1}^2,
\left(\frac{\beta_{\frak{b}}}{\gamma_{\frak{b}}}\right)u_{s(\frak{b}+1)}^2
\right\}&\leq &M_{10}\,
\kappa_{s(\frak{d}+1)}\max\left\{u_{s(\frak{d}+1)-1}^2,
\left(\frac{\beta_{\frak{d}}}{\gamma_{\frak{d}}}\right)u_{s(\frak{d}+1)}^2\right\}\\
&\leq & M_{11}\, \kappa_{s(\frak{c})} \,
(9\kappa^4)^{\Gap(\I)}\max\left\{u_{s(\frak{c})-1}^2,
\left(\frac{\beta_{\frak{c}-1}}{\gamma_{\frak{c}-1}}\right)u_{s(\frak{c})}^2\right\}\\
&\leq & M_{12}\,  (9\kappa^4)^{\Gap(\I)}
~\prod_{j=1}^{\frak{a}-1}\left(\left(\frac{\beta_i}{\gamma_i}\right)\rho_i^2\right
)^{\ell_i}
\end{eqnarray*}
holds, and moreover, if $\frak{d}<\frak{b}$ then for each
$\frak{d}+1\leq i\leq \frak{b}$,
\begin{eqnarray*}
\sum_{j=0}^{\ell(i)+1}\kappa_{s(i)-1+j}u_{s(i)-1+j}^2 &\leq &
M_{13}\, \ell(i)\,  (9\kappa^4)^{\Gap(\I)}
\prod_{j=1}^{\frak{c}-1}\left(\frac{\beta_j}{\gamma_j}\right)^{\ell(j)}\max\left\{u_{s(\frak{c})-1}^2,
\left(\frac{\beta_{\frak{c}}}{\gamma_{\frak{c}}}\right)u_{s(\frak{c})}^2\right\}\\
&\leq & M_{14} \, \ell(i) \, (9\kappa^4)^{\Gap(\I)}
\prod_{j=1}^{\frak{a}-1}\left(\left(\frac{\beta_j}{\gamma_j}\right)
\rho_j^2 \right)^{\ell(j)}
\end{eqnarray*}
and
\[\sum_{j=0}^{\ell(i)+1}\kappa_{s(i)-1+j}u_{s(i)-1+j}^2\geq M_{15} \, \ell(i)
\left( \frac{1}{9\kappa^4} \right)^{\Gap(\I)}
~\prod_{j=1}^{\frak{a}-1}\left(\left(\frac{\beta_j}{\gamma_j}\right)
\rho_j^2 \right)^{\ell(j)}\] all hold. By taking
$C_8(\kappa,\g,\I):=\frac{M_{15}}{3}$, $C_9(\kappa,\g,\I):=M_{14}$
and $C_{11}(\kappa,\mathcal{G},\I):=M_{12}$, it can be seen that
the inequalities in (vi) and (vii) involving these constants all
hold. It can also be seen in a similar fashion that there exists a
constant $C_{10}=C_{10}(\kappa,\g,\I)>0$ such that the left-hand
inequality in (vii) holds. \epf

By using the previous theorem, we now obtain bounds for Gap sum
$\sum_{i=s(\frak{a})-1}^{s(\frak{b}+1)}\kappa_iu_i^2$.

\begin{cor}\label{gap-sum}
Let $\kappa\geq 3$ and $\lambda\geq 0$ be integers with
$\lambda\leq \kappa-2$. Let
$\g=\Big((\gamma_i,\alpha_i,\beta_i)\Big)_{i=1}^{g+1}$ be a
$(\kappa,\lambda)$-graphical sequence, $\ell:\{1,\dots,g+1\} \to
\N$ be a function with $\ell(g+1)=1$, and $\T:=\T(\g,\ell)$ be the
$(\kappa,\lambda)$-tridiagonal sequence associated to $\g$ and
$\ell$. Suppose that $\I$ is a well-placed interval with respect
to $\g$. Then there exist positive constants
$C:=C(\kappa,\mathcal{G},\I)$ and $C':=C'(\kappa,\mathcal{G},\I)$
such that for any element $\theta$ in $\I$,
\begin{eqnarray}
C\left(\frac{1}{9\kappa^4}\right)^{\mbox{\em
\Gap}(\I)}\,{\mbox{\em \Len}}(\I)\, \prod_{i=1}^{\frak{a}-1}
\left(
\left(\frac{\beta_{i}}{\gamma_i} \right) \rho_i^2 \right)^{\ell(i)} &\leq&
\sum_{i=s(\frak{a})-1}^{s(\frak{b}+1)}\kappa_iu_i^2\nonumber\\
 &\leq & C' (9\kappa^4 )^{\mbox{\em \Gap}(\I)}\,{\mbox{\em \Len}}(\I)\,
 \prod_{i=1}^{\frak{a}-1}\left(\left(\frac{\beta_{i}}{\gamma_i}\right)\rho_i^2\right)^{\ell(i)}~,\nonumber
\end{eqnarray}
where $\kappa_i$ and $u_i:=u_i(\theta)$ are as defined in
(\ref{def-k_i}) and (\ref{vi-ui}) for the matrix
$L_1(\T(\g,\ell))$, and $s(i)$, $\frak{a}:=\frak{a}(\g,\I)$, $\frak{b}:=\frak{b}(\g,\I)$, ${\mbox{\em
\Len}}(\I):={\mbox{\em \Len}}_{\g,\ell}(\I)$, $\mbox{\em
\Gap}(\I):=\mbox{\em \Gap}_{\g,\ell}(\I)$ and
$\rho_i:=\rho_i(\theta)$ are as defined in (\ref{s(i)}), (\ref{a}), (\ref{b}), (\ref{ell(I)})-(\ref{def-rho}),
respectively.
\end{cor}
\pf Constants $C_i~(i=3,4,7,8,9)$ in this proof are the constants
in Theorem \ref{0-s-sum}. Note that $\Len(\I)\geq 1$ and $\left(9
\kappa^4 \right)^{\Gap(\I)}\geq 1$. We break the proof into three
cases:

\noindent (1) $\frak{c}=g+1$: By (\ref{ell(I)}) and
(\ref{def-gap}), $\Len(\I)=\sum_{i=\frak{a}}^{\frak{b}}\ell(i)\geq
1$ and $\Gap(\I)=0$. By applying Theorem~\ref{0-s-sum} (iii) with
$\widehat{\frak{c}}=\frak{b}+1$, Corollary \ref{gap-sum} holds for
$C:=C_3$ and $C':=C_4$.

\noindent (2) $\frak{c}\leq g$ and $\frak{d}=\frak{b}$: Then
$\Len(\I)=\sum_{i=\frak{a}}^{\frak{c}-1}\ell(i)\geq 1$, and by
applying Theorem~\ref{0-s-sum} (iii) and (v) for
$\widehat{\frak{c}}=\frak{c}$, the result follows for $C:=C_3$ and
$C':=C_4+C_7$ as

\begin{eqnarray*}
&& C_3\left(\frac{1}{9\kappa^4}\right)^{ \Gap(\I)}\,{ \Len }(\I)\,
\prod_{i=1}^{\frak{a}-1} \left( \left(\frac{\beta_{i}}{\gamma_i}
\right) \rho_i^2 \right)^{\ell(i)} \leq  C_3 {\Len}(\I)\,
\prod_{i=1}^{\frak{a}-1} \left( \left(\frac{\beta_{i}}{\gamma_i}
\right) \rho_i^2 \right)^{\ell(i)}\nonumber \\ &\leq&
\sum_{i=s(\frak{a})-1}^{s(\frak{c})}\kappa_iu_i^2\nonumber\\
 &\leq & C_4 {\Len }(\I)\,
 \prod_{i=1}^{\frak{a}-1}\left(\left(\frac{\beta_{i}}{\gamma_i}\right)\rho_i^2\right)^{\ell(i)}
  \leq C_4 (9\kappa^4 )^{\Gap (\I)}\,{\Len }(\I)\,
 \prod_{i=1}^{\frak{a}-1}\left(\left(\frac{\beta_{i}}{\gamma_i}\right)\rho_i^2\right)^{\ell(i)}
\end{eqnarray*}
and
\begin{eqnarray*}
\sum_{i=s(\frak{c})}^{s(\frak{b}+1)}\kappa_iu_i^2 < C_7 (9\kappa^4
)^{\Gap (\I)}\,
 \prod_{i=1}^{\frak{a}-1}\left(\left(\frac{\beta_{i}}{\gamma_i}\right)\rho_i^2\right)^{\ell(i)}
\leq  C_7 (9\kappa^4 )^{\Gap (\I)}\,{ \Len }(\I)\,
 \prod_{i=1}^{\frak{a}-1}\left(\left(\frac{\beta_{i}}{\gamma_i}\right)\rho_i^2\right)^{\ell(i)}
\end{eqnarray*}
all hold.

\noindent (3) $\frak{c}<g$ and $\frak{d}<\frak{b}$: In this case,
$\Len(\I)=\sum_{i=\frak{a}}^{\frak{c}-1}\ell(i)+\sum_{i=\frak{d}+1}^{\frak{b}}\ell(i)$
and by Theorem~\ref{0-s-sum} (iii), (v) and (vi), the following
all hold:
\begin{eqnarray*}
&& \min\{C_3,C_8\}\left(\frac{1}{9\kappa^4}\right)^{
\Gap(\I)}\sum_{i=\frak{a}}^{\frak{c}-1}\ell(i)\,\prod_{i=1}^{\frak{a}-1}
\left( \left(\frac{\beta_{i}}{\gamma_i} \right) \rho_i^2
\right)^{\ell(i)} \leq  C_3
\sum_{i=\frak{a}}^{\frak{c}-1}\ell(i)\, \prod_{i=1}^{\frak{a}-1}
\left( \left(\frac{\beta_{i}}{\gamma_i} \right) \rho_i^2
\right)^{\ell(i)}\nonumber \\ &\leq&
\sum_{i=s(\frak{a})-1}^{s(\frak{c})}\kappa_iu_i^2\nonumber\\
 &\leq & C_4 \sum_{i=\frak{a}}^{\frak{c}-1}\ell(i)\,
 \prod_{i=1}^{\frak{a}-1}\left(\left(\frac{\beta_{i}}{\gamma_i}\right)\rho_i^2\right)^{\ell(i)}
  \leq \max\{C_4,C_7,C_9\} (9\kappa^4 )^{\Gap (\I)}\,{\Len }(\I)\,
 \prod_{i=1}^{\frak{a}-1}\left(\left(\frac{\beta_{i}}{\gamma_i}\right)\rho_i^2\right)^{\ell(i)},
\end{eqnarray*}
\begin{eqnarray*}
\sum_{i=s(\frak{c})}^{s(\frak{d+1})}\kappa_iu_i^2\leq
C_7(9\kappa^4 )^{\Gap (\I)}\,
 \prod_{i=1}^{\frak{a}-1}\left(\left(\frac{\beta_{i}}{\gamma_i}\right)\rho_i^2\right)^{\ell(i)}
  \leq \max\{C_4,C_7,C_9\} (9\kappa^4 )^{\Gap (\I)}\,{\Len }(\I)\,
 \prod_{i=1}^{\frak{a}-1}\left(\left(\frac{\beta_{i}}{\gamma_i}\right)\rho_i^2\right)^{\ell(i)}\\
\end{eqnarray*}
and
\begin{eqnarray*}
&& \min\{C_3,C_8\}\left(\frac{1}{9\kappa^4}\right)^{
\Gap(\I)}\sum_{i=\frak{d}+1}^{\frak{b}}\ell(i)\,\prod_{i=1}^{\frak{a}-1}
\left( \left(\frac{\beta_{i}}{\gamma_i} \right) \rho_i^2
\right)^{\ell(i)}\leq
\sum_{i=s(\frak{d}+1)-1}^{s(\frak{b}+1)}\kappa_iu_i^2\nonumber\\
 &\leq & C_9 (9\kappa^4)^{
\Gap(\I)}\sum_{i=\frak{d}+1}^{\frak{b}}\ell(i)\,\prod_{i=1}^{\frak{a}-1}
\left( \left(\frac{\beta_{i}}{\gamma_i} \right) \rho_i^2
\right)^{\ell(i)}\leq \max\{C_4,C_7,C_9\} (9\kappa^4)^{
\Gap(\I)}\Len(\I)\,\prod_{i=1}^{\frak{a}-1} \left(
\left(\frac{\beta_{i}}{\gamma_i} \right) \rho_i^2
\right)^{\ell(i)}.
\end{eqnarray*}
Hence, the result now follows by taking
\[C:=\frac{\min\{C_3,C_8\}}{2},~~C':=3\,\max\{C_4,C_7,C_9\},\]
in light of the fact $s(\frak{c})\leq s(\frak{d}+1)-1$. The
corollary now follows. \epf

\subsection{Bounding Tail Sum}\label{btail}

In this section, we obtain an upper bound for the Tail sum
$\sum_{i=s(\frak{b}+1)+1}^{D}\kappa_iu_i^2$. Namely:\\

\begin{theorem}\label{s-d-sum}
Let $\kappa \geq 3$ and $\lambda\geq 0$ be integers with $\lambda
\leq \kappa-2$. Suppose that $(\g,\Delta;L,\ell)$ is a
$(\kappa,\lambda)$-quadruple and let $\g=\Big(
(\gamma_i,\alpha_i,\beta_i)\Big)_{i=1}^{g+1}$. Suppose that $\I$
is a well-placed interval with respect to $\g$ satisfying $\I\cap
\mathcal{B}=\emptyset$ and $\frak{b}<g$ (cf. (\ref{bad roots-def})
and (\ref{b})). Then there exist positive constants
$F:=F(\kappa,\g,\Delta,L,\I)~$ and
$C:=C(\kappa,\mathcal{G},\Delta,L,\I)$ so that if $\ell(i)>F$
holds for all $(\gamma_{i},\alpha_{i},\beta_{i})\in
\underline{\Delta}$ with $\frak{b}<i \leq g$ then, for any $\theta
\in \I $,
\[\sum_{i=s(\frak{b}+1)+1}^{D}\kappa_iu_i^2\leq C (9\kappa^4)^{\mbox{\em \Gap}(\I)}\prod_{i=1}^{\frak{a}-1}
\left(\left( \frac{\beta_i}{\gamma_i}\right) \rho_i^2
\right)^{\ell(i)} \] holds, where $\kappa_i$ and
$u_i:=u_i(\theta)$ are as defined in (\ref{def-k_i}) and
(\ref{vi-ui}) for the matrix $L_1(\T(\g,\ell))$, and $s(i)$,
$D:=D_{\T(\g,\ell)}$, $\frak{a}:=\frak{a}(\g,\I)$, $\mbox{\em
\Gap}(\I):=\mbox{\em \Gap}_{\g,\ell}(\I)$ and
$\rho_i:=\rho_i(\theta)$ are as defined in (\ref{s(i)}),
(\ref{D(T)}), (\ref{a}), (\ref{def-gap}) and (\ref{def-rho}),
respectively.
\end{theorem}

\pf Suppose $(\g,\Delta;L,\ell)$ and $\I$ are as in the statement
of the theorem. By Theorem~\ref{Lj}, for each $0\leq i\leq
g-\frak{b}-1$ satisfying
$(\gamma_{g-i},\alpha_{g-i},\beta_{g-i})\in \underline{\Delta}$
there exist constants $C_i=C_i(\kappa,\g,\Delta,L,\I)\geq 1$ and
$M_i=M_i(\kappa,\g,\Delta,L,\I)>1$ such that if $\ell(g-j)>C_i $
holds for all $(\gamma_{g-j},\alpha_{g-j},\beta_{g-j})\in
\underline{\Delta}$ with $j<i$, then (\ref{eq-L_i}) holds for all
$\theta\in \I $. Now put
\begin{eqnarray*}
F&=&F(\kappa,\g,\Delta,L,\I):=\max\{C_i(\kappa,\g,\Delta,L,\I)\,:\,0\leq i\leq g-\frak{b}-1\mbox{~and~}
(\gamma_{g-i},\alpha_{g-i},\beta_{g-i})\in \underline{\Delta}
\};\\
M&=&M(\kappa,\g,\Delta,L,\I):=\min\{M_i(\kappa,\g,\Delta,L,\I)\,:\,0\leq i\leq g-\frak{b}-1\mbox{~and~}
(\gamma_{g-i},\alpha_{g-i},\beta_{g-i})\in \underline{\Delta}
\}.
\end{eqnarray*}
\begin{itemize}
\item[($\dag$)] Suppose that if $\{(\gamma_{g-i},\alpha_{g-i},\beta_{g-i})\in
\underline{\Delta}\,:\,0\leq i\leq g-\frak{b}-1\}\neq \emptyset$
then $\ell(g-i)>F$ holds for all $(\gamma_{g-i},\alpha_{g-i},\beta_{g-i})\in \underline{\Delta}$ with
$0\leq i\leq g-\frak{b}-1$.
\end{itemize}
Let $\theta \in \I$. We will use the following:
\begin{claim} \label{tail-claim}
There exist constants $C_1=C_1(\g,\I)>0$
and $C_m=C_m(\kappa,\g,\Delta,L,\I)>0~(m=2,3)$
such that, for all $0\leq i\leq g-\frak{b}-1$, the following hold:\\
(a)
\[
|u_{s(g-i+1)-j}|\leq C_1 \max\{|u_{s(g-i+1)-1}|,|u_{s(g-i+1)}|
\}~x_i^j~~(0\leq j\leq \ell(g-i)+1).
\]
(b)
\[ \max \{|u_{s(g-i)-1}|, |u_{s(g-i)}|\} > C_2 \max \{ |u_{s(g-i+1)-1}|, |u_{s(g-i+1)}|\}~x_i^{\ell(g-i)}.\]
(c)
\[\prod_{j=0}^{g-\frak{b}-1}{x_j}^{2 \ell(g-j)}\max\{u_{D-1}^2, u_{D}^2 \}< C_3
\max\{u_{s(\frak{b}+1)-1}^2, u_{s(\frak{b}+1)}^2 \},\] where $x_j$
is defined in (\ref{def-x_i-y_i}).
\end{claim}

\noindent{\bf Proof of Claim \ref{tail-claim}}: Let $0\leq i\leq
g-\frak{b}-1$. By (\ref{def-x_i-y_i}) and
Lemma~\ref{h-m-t-partition} (iv), $x_i > y_i>0$. Let
$\nu_j^{(i)}=\nu_j^{(i)}(\theta)~~(j=1,2)$ be as defined in
(\ref{rec-tail}). \vskip0.01cm
\noindent (a): First suppose that
$\nu_1^{(i)}\nu_2^{(i)}>0$ holds. Then for all $0\leq j\leq
\ell(g-i)+1$, (a) follows since
\[|u_{s(g-i+1)-j}|=|\nu_1^{(i)}|\, x_i^j+|\nu_2^{(i)}|\, y_i^j<(|\nu_1^{(i)}|+|\nu_2^{(i)}|) \, x_i^j
=|u_{s(g-i+1)}|\, x_i^j.\]
Now suppose $\nu_1^{(i)}\nu_2^{(i)}<0$. By (\ref{nu1i}) and (\ref{nu2i}),
\begin{eqnarray}
\max\{|\nu_1^{(i)}|, |\nu_2^{(i)}|\} &\leq & 2
\max \left \{ \frac{1}{x_i-y_i},  \frac{x_i}{x_i-y_i}\right\} \max\{|u_{s(g-i+1)-1}|, |u_{s(g-i+1)}|\}\nonumber\\
&\leq & C_1 \max\{ |u_{s(g-i+1)-1}|,
|u_{s(g-i+1)}|\}\label{e-ineq}
\end{eqnarray}
holds, where
$$
C_1=C_1(\g,\I):=2
\max\left\{\frac{\I_{\max}-\alpha_{m}}{\sqrt{(\I_{\min}-\alpha_{m})^2-4\beta_{m}
\gamma_{m}}},
\frac{\gamma_{m}}{\sqrt{(\I_{\min}-\alpha_{m})^2-4\beta_{m}
\gamma_{m}}} \,:\, 0\leq m \leq  g-\frak{b}-1 \right\}.
$$
Since $|u_{s(g-i+1)-j}|\leq \max\{|\nu_1^{(i)}|,
|\nu_2^{(i)}|\}\,x_i^j$ holds by $\nu_1^{(i)}\nu_2^{(i)}<0$ and $x_i > y_i>0$,
(a) follows by (\ref{e-ineq}).\\

\noindent{(b):} Suppose $(\gamma_{g-i},\alpha_{g-i},\beta_{g-i})\in
\underline{\g} \setminus \underline{\Delta}$. Then $\ell(g-i)=
L(g-i)$ and, by Lemma~\ref{bounder} (i) and
$0<x_i<\frac{\I_{\max}-\alpha_{g-i}}{\gamma_{g-i}}$, it follows that
\begin{eqnarray*}
\max \{ |u_{s(g-i)-1}|, |u_{s(g-i)}| \} &\geq & \left(\frac{1}{3
\kappa\,x_i}\right) ^{\ell(g-i)} \max \{|u_{s(g-i+1)-1}|,
|u_{s(g-i+1)}|\}\, x_i^{\ell(g-i)}\\ &\geq & \left(
\frac{\gamma_{g-i}}{3 \kappa
(\I_{\max}-\alpha_{g-i})}\right)^{L(g-i)} \max \{|u_{s(g-i+1)-1}|,
|u_{s(g-i+1)}|\}\,x_i^{\ell(g-i)},
\end{eqnarray*}
and thus (b) follows by taking $C_2=C_2(\kappa,\g,\Delta,L,\I)$,
where
 \[C_2:=\min\left\{\left( \frac{\gamma_{g-m}}{3 \kappa
(\I_{\max}-\alpha_{g-m})}\right)^{L(g-m)} \,:\, 0\leq m\leq
g-\frak{b}-1 \mbox{ and
}(\gamma_{g-m},\alpha_{g-m},\beta_{g-m})\in \underline{\g}
\setminus \underline{\Delta}\right\}.\] Now suppose
$(\gamma_{g-i},\alpha_{g-i},\beta_{g-i})\in \underline{\Delta}$.
By Theorem~\ref{Lj} with ($\dag$),
\begin{eqnarray}
\max \{ |u_{s(g-i+1)-1}|, |u_{s(g-i+1)}| \} &\leq & \max\{1,x_i\}(|\nu_1^{(i)}|+|\nu_2^{(i)}|)\nonumber \\
&<& \left(1+ \frac{1}{M}
\left(\frac{x_i}{y_i}\right)^{F} \right)(1+x_i)\, |\nu_1^{(i)}|\nonumber\\ &<
&\left(1+\frac{1}{M}\left(\frac{(\I_{\max}-\alpha_{g-i})^2}{\beta_{g-i}\gamma_{g-i}}\right)^{F} \right)
\left(1+\frac{\I_{\max}-\alpha_{g-i}}{\gamma_{g-i}}\right)
\,|\nu_1^{(i)}|\label{tail-nu1-bd}
\end{eqnarray}
and
\begin{eqnarray}
\max \{ |u_{s(g-i)-1}|, |u_{s(g-i)}| \} &\geq &
|\nu_1^{(i)}x_i^{\ell(g-i)}|-|\nu_2^{(i)}y_i^{\ell(g-i)}|\nonumber \\
&>& |\nu_1^{(i)}|\,x_i^{\ell(g-i)} \left(1-\frac{1}{M}\left(\frac{x_i}{y_i}\right)^{F}
\left(\frac{y_i}{x_i}\right)^{\ell(g-i)}\right)\nonumber\\
&> &
\left(1-\frac{1}{M}\right)|\nu_1^{(i)}|\,x_i^{\ell(g-i)}\label{tail-nu1-bd-1}
\end{eqnarray}
all hold. By (\ref{tail-nu1-bd}) and (\ref{tail-nu1-bd-1}),
statement (b) now follows by taking
$C_2=C_2(\kappa,\mathcal{G},\Delta,L,\I)$, where
$$
C_2:=\frac{1-\frac{1}{M}} {\max\left\{
\left(1+\frac{1}{M}\left(\frac{(\I_{\max}-\alpha_{g-m})^2}{\beta_{g-m}\gamma_{g-m}}\right)^{F}
\right) \left(1+\frac{\I_{\max}-\alpha_{g-m}}{\gamma_{g-m}}\right)
\,:\, 0\leq m\leq g-\frak{b}-1 \mbox{ and
}(\gamma_{g-m},\alpha_{g-m},\beta_{g-m})\in
\underline{\Delta}\right\}}.
$$

\noindent{(c):} This follows by applying (b) inductively on $i$
for $0\leq i\leq g-\frak{b}-1$.\epf

Let $0\leq i\leq g-\frak{b}-1$. By (a) and (c) of the claim,
$\left(\frac{\gamma_{g-i}}{\beta_{g-i}}\right)x_i^2>1$ and
Theorem~\ref{0-s-sum} (vii), there exist constants
$M_j=M_j(\kappa,\g,\I)>0~(j=1,2)$ and
$M_j=M_j(\kappa,\g,\Delta,L,\I)>0~(j=3,4)$ such that
\begin{eqnarray}
&& \sum_{j=0}^{\ell(g-i)-1}\kappa_{s(g-i+1)-j}u_{s(g-i+1)-j}^2
\nonumber\\&<& \kappa \kappa_D
\prod_{j=0}^{i-1}\left(\frac{\gamma_{g-j}}{\beta_{g-j}}\right)^{\ell(g-j)}~\sum_{m=0}^{\ell(g-i)-1}
\left(\frac{\gamma_{g-i}}{\beta_{g-i}}\right)^m u_{s(g-i+1)-m}^2
\nonumber\\ &\leq & M_1 \kappa_D
\prod_{j=0}^{i-1}\left(\frac{\gamma_{g-j}}{\beta_{g-j}}\right)^{\ell(g-j)}
\max\{u_{s(g-i+1)-1}^2,  u_{s(g-i+1)}^2 \}
\sum_{m=0}^{\ell(g-i)-1}
\left(\frac{\gamma_{g-i}}{\beta_{g-i}}\right)^m
x_i^{2m}\nonumber\\ &\leq & M_2 \kappa_D \max\{u_{D-1}^2,
u_D^2\}\prod_{j=0}^{i}\left(\left(\frac{\gamma_{g-j}}{\beta_{g-j}}\right)x_j^2\right)^{\ell(g-j)}\nonumber
\\ &\leq & M_2 \kappa_D \max\{u_{D-1}^2, u_D^2
\}\prod_{j=0}^{g-\frak{b}-1}\left(\left(\frac{\gamma_{g-j}}{\beta_{g-j}}\right)x_j^2\right)^{\ell(g-j)}\nonumber
\\ &\leq & M_3 \kappa_{s(\frak{b}+1)} \max\{u_{s(\frak{b}+1)-1}^2,u_{s(\frak{b}+1)}^2\}\nonumber\\
&\leq & M_4 (9\kappa^4)^{\Gap(\I)}
\prod_{j=1}^{\frak{a}-1}\left(\left(\frac{\beta_j}{\gamma_j}\right)\rho_j^2\right)^{\ell(j)}\label{tail-sum-element}
\end{eqnarray}
holds. From (\ref{tail-sum-element}) and
\[
\sum_{j=s(\frak{b}+1)+1}^{D}\kappa_ju_j^2=\sum_{i=0}^{g-\frak{b}-1}~\sum_{j=0}^{\ell(g-i)-1}
\kappa_{s(g-i+1)-j}u_{s(g-i+1)-j}^2~,
\]
Theorem~\ref{s-d-sum} now follows
by taking $C(\kappa,\g,\Delta,L,\I):=(g-\frak{b})M_4$.
\epf

With these results in hand, we can now prove
the main theorem of this section:

{\noindent{\em Proof of Theorem \ref{multi-thm}:}} Theorem
\ref{multi-thm} follows immediately by Theorem \ref{0-s-sum} (i),
Corollary \ref{gap-sum} and Theorem \ref{s-d-sum}. \epf

\section{Distribution of Eigenvalues and Proof of Theorem \ref{b-i-small-thm}}\label{distri}

In this section we prove Theorem~\ref{b-i-small-thm} and thus
complete the proof of the Bannai-Ito conjecture. To do this we will
first prove two results concerning the distribution of the
eigenvalues of a graphical sequence in a well-placed interval with
respect to this sequence, using the results from the last four
sections.

\begin{theorem}\label{close-evs-gap(I)}
Let $\kappa\geq 3$ and $\lambda\geq 0$ be integers with
$\lambda\leq \kappa-2$, and let $\g=\Big(\delta_i:=(\gamma_i,\alpha_i,\beta_i)\Big)_{i=1}^{g+1}$ be a
$(\kappa,\lambda)$-graphical sequence. Suppose that $\Delta=(\delta_{i_p})_{p=1}^{\tau}$ is a subsequence of $\g$ with
$(1,\lambda,\kappa-\lambda-1) \in \underline{\Delta}$ and
$(\gamma_{g+1},\alpha_{g+1},\beta_{g+1}) \not\in \underline{\Delta}$, $L:\{1,\ldots,g+1\}\setminus
\{i_1,\ldots,i_{\tau}\} \rightarrow \mathbb{N}$ is a function, and $\I$ is a well-placed interval with respect to $\g$ satisfying $\I\cap \mathcal{B}(\g,\Delta,L)=\emptyset$ (cf. (\ref{bad roots-def})). Suppose that $\epsilon>0$ is a real number, $C:=C(\kappa)>0$ is a constant, and $\ell: \{1,\ldots,g+1\}\rightarrow \mathbb{N}$ is
any function for which $(\g, \Delta; L,\ell)$ is a $(\kappa,\lambda)$-quadruple and the associated $(\kappa, \lambda)$-tridiagonal sequence
$\mathcal{T}=\mathcal{T}(\g,\ell)$ satisfies\\
(i) Property (AC),\\
(ii) $D_{\T}\leq C \mbox{\em \HH}_{\T}$, and\\
(iii) $\mbox{\em \Len}(\I)\geq \epsilon \mbox{\em \HH}_{\T}$,\\
where $\mbox{\em
\HH}_{\T}$, $D_{\T}$ and $\mbox{\em \Len}(\I):=\mbox{\em \Len}_{\g,\ell}(\I)$ are as defined in (\ref{def-h}), (\ref{D(T)}) and (\ref{ell(I)}), respectively.\\
Then for any real number $\delta>0$, there exist positive
constants $F:=F(\kappa,\g,\Delta,L,\I)$,
$C_1:=C_1(\kappa,\lambda,\epsilon,\delta,\g,\Delta,L,\I)$ and
$C_2:=C_2(\kappa,\lambda,\delta)$ such that if $\ell(i_p)>F$ holds
for all $1\leq p\leq \tau$ and if there exist two conjugate
algebraic numbers $\theta $ and $\eta$ in $\mathcal{E}_{\T}\cap
\I$ satisfying $|\theta-\eta|>\delta$ then
\[\mbox{either~~~}\mbox{\em \HH}_{\T}<C_1\mbox{~~or~~}\mbox{\em \Gap}(\I)> C_2\mbox{\em \HH}_{\T}\]
holds, where $\mathcal{E}_{\mathcal{T}}$ and $\mbox{\em \Gap}(\I):=\mbox{\em \Gap}_{\g,\ell}(\I)$ are as defined in (\ref{def-Eg}) and (\ref{def-gap}), respectively.
\end{theorem}

\pf Suppose that $\kappa$, $\lambda$, $\epsilon$, $C$, $\g$,
$\Delta$, $L$, $\I$, $\ell$ and $\T$ are as in the statement of
the theorem, and put $\HH:=\HH_{\T}$ and $D:=D_{\T}$. Let $\delta$
be any positive real number, and let $\theta$ and $\eta$ be two
conjugate algebraic numbers in $\mathcal{E}_{\mathcal{T}}\cap \I$
satisfying $|\theta -\eta|>\delta$. Without loss of generality, we
assume $\eta-\theta
> \delta$. \\

By applying Theorem~\ref{multi-thm} and the conditions
$\epsilon\, \HH\leq \Len(\I)< D\leq C\,\HH$ given by (ii) and (iii) in the statement of the theorem, it follows that there exist
positive constants $F:=F(\kappa,\g,\Delta,L,\I)$, $M_1:=M_1(\kappa,\g,\I)$ and
$M_2:=M_2(\kappa,\g,\Delta,L,\I)$ so that if $\ell(i_p)>F$ holds for all $1\leq p\leq \tau$ then
\begin{equation}\label{no-thy-1}
\epsilon\,\HH \, M_1\left(\frac{1}{9\kappa^4}\right)^{\Gap (
\I)}\prod_{i=1}^{\frak{a}-1}\left(\left(\frac{\beta_i}{\gamma_i}\right)
\rho_i^2(x)\right)^{\ell(i)} \le
\sum_{i=0}^{D}\kappa_iu_i^2(x) \le \HH\,M_2 C (9 \kappa^4
)^{\Gap (\I)}
\prod_{i=1}^{\frak{a}-1}\left(\left(\frac{\beta_i}{\gamma_i}\right)\rho_i^2(x)\right)^{\ell(i)}
\end{equation}
holds for any $x\in \I$, where $\kappa_i$ and $u_i:=u_i(x)$ are as defined in (\ref{def-k_i}) and
(\ref{vi-ui}) for the matrix $L_1(\T)$, and $\frak{a}$, $\rho_i(x)$
are as defined in (\ref{a}) and (\ref{def-rho}), respectively.\\
By Proposition \ref{new-ell-3.1} (i) and $\eta>\theta$, it follows that
\begin{equation}\label{rhoi(theta)<rhoi(theta')}
0<\rho_i(\theta)<\rho_i(\eta)<1~~(i=1,\ldots, \frak{a}-1),
\end{equation}
and moreover, by (\ref{rhoi(theta)<rhoi(theta')}) and $\eta-\theta>\delta$,
\begin{equation}\label{rho1(theta)<rho1(eta)}
\rho_1(\eta)>\rho_1(\theta)+\frac{\delta}{2(\kappa-\lambda-1)}>
\left(1+\frac{\delta}{2(\kappa-\lambda-1)}\right)\rho_1(\theta).
\end{equation}
By applying
(\ref{rhoi(theta)<rhoi(theta')}) and (\ref{rho1(theta)<rho1(eta)}) to (\ref{no-thy-1}), it follows that
\begin{eqnarray}
\sum_{i=0}^{D}\kappa_iu_i^2(\eta)&\geq & \epsilon\,\HH \, M_1\left(\frac{1}{9\kappa^4}\right)^{\Gap (
\I)}\prod_{i=1}^{\frak{a}-1}\left(\left(\frac{\beta_i}{\gamma_i}\right)
\rho_i^2(\eta)\right)^{\ell(i)} \nonumber \\
&>& \epsilon\,\HH \, M_1\left(\frac{1}{9\kappa^4}\right)^{\Gap (
\I)}  \left(1+\frac{\delta}{2(\kappa-\lambda-1)}\right)^{2\HH}  \prod_{i=1}^{\frak{a}-1}\left(\left(\frac{\beta_i}{\gamma_i}\right)
\rho_i^2(\theta)\right)^{\ell(i)} \nonumber \\
&\geq & \frac{\epsilon M_1}{M_2 C}\left(\frac{1}{9\kappa^4}\right)^{2\Gap (\I)}
\left(1+\frac{\delta}{2(\kappa-\lambda-1)}\right)^{2\HH}\sum_{i=0}^{D}\kappa_iu_i^2(\theta).\label{ln-no-thy-0}
\end{eqnarray}
Since $\theta $ and $\eta$ are algebraic conjugates,
$\sum_{i=0}^{D}\kappa_iu_i^2(\eta)=\sum_{i=0}^{D}\kappa_iu_i^2(\theta)>0$
holds by Property (AC). Hence, by (\ref{ln-no-thy-0}),
\begin{equation}\label{ln-no-thy-1}
\ln
\left(1+\frac{\delta}{2(\kappa-\lambda-1)}\right)^2<\frac{\ln\left(\frac{M_2
C}{\epsilon M_1}\right)}{\HH}+\frac{\Gap (\I)}{\HH}\ln (9\kappa
^4)^2.
\end{equation}
Now, put
\[
C_1:=\frac{\ln \left(\frac{M_2 C}{\epsilon
M_1}\right)}{\ln \left(1+\frac{\delta}{2(\kappa-\lambda-1)}
\right)}\mbox{~~~and~~~}C_2:=\frac{\ln
\left(1+\frac{\delta}{2(\kappa-\lambda-1)}\right)}{2 \ln (9\kappa
^4)}.
\]

If $0< \frac{M_2 C}{\epsilon M_1}\leq 1$ then $\Gap(\I)>C_2 \HH$ holds as $\frac{\Gap
(\I)}{\HH}>\frac{\ln
\left(1+\frac{\delta}{2(\kappa-\lambda-1)}\right)}{\ln (9\kappa
^4)}>\frac{\ln
\left(1+\frac{\delta}{2(\kappa-\lambda-1)}\right)}{2 \ln (9\kappa
^4)}$ by (\ref{ln-no-thy-1}). Moreover, if $\frac{M_2 C}{\epsilon M_1}> 1$ and $\Gap(\I)\leq C_2 \HH$, then $\HH<\frac{\ln \left(\frac{M_2 C}{\epsilon
M_1}\right)}{\ln \left(1+\frac{\delta}{2(\kappa-\lambda-1)}
\right)}$ holds by (\ref{ln-no-thy-1}). Therefore Theorem \ref{close-evs-gap(I)} now follows for this choice of $C_1$ and $C_2$. \epf

\begin{prop}\label{main-prop}
Let $\kappa\geq 3$ and $\lambda\geq 0$ be integers with
$\lambda\leq \kappa-2$, and let $\g=\Big(\delta_i:=(\gamma_i,\alpha_i,\beta_i)\Big)_{i=1}^{g+1}$ be a $(\kappa,\lambda)$-graphical sequence. Suppose that $\Delta=(\delta_{i_p})_{p=1}^{\tau}$ is a subsequence of $\g$ with $(1,\lambda,\kappa-\lambda-1) \in \underline{\Delta}$ and
$(\gamma_{g+1},\alpha_{g+1},\beta_{g+1}) \not\in \underline{\Delta}$, $L:\{1,\ldots,g+1\}\setminus \{i_1,\ldots,i_{\tau}\} \rightarrow \mathbb{N}$ is a function, and $\I$ is a well-placed interval with respect to $\g$ satisfying $\I\cap \mathcal{B}(\g,\Delta,L)=\emptyset$ (cf. (\ref{bad roots-def})). Suppose that $\epsilon>0$ is a real number, $C:=C(\kappa)>0$ is a constant, and $\ell: \{1,\ldots,g+1\}\rightarrow \mathbb{N}$ is
any function for which $(\g, \Delta; L,\ell)$ is a $(\kappa,\lambda)$-quadruple and the associated $(\kappa, \lambda)$-tridiagonal sequence
$\mathcal{T}=\mathcal{T}(\g,\ell)$ satisfies\\
(i) Property (AC),\\
(ii) $D_{\T}\leq C \mbox{\em \HH}_{\T}$, and\\
(iii) $\mbox{\em \Len}(\I)\geq \epsilon \mbox{\em \HH}_{\T}$,\\
where $\mbox{\em\HH}_{\T}$, $D_{\T}$ and $\mbox{\em \Len}(\I):=\mbox{\em \Len}_{\g,\ell}(\I)$ are as defined in (\ref{def-h}), (\ref{D(T)}) and (\ref{ell(I)}), respectively.\\ Then for any real number $\mu >0$, there exist positive constants $F:=F(\kappa,\g,\Delta,L,\I)$,
$G:=G(\kappa,\lambda,\epsilon,\mu,\g,\I)$ and
$H:=H(\kappa,\lambda,\epsilon,\mu,\g,\Delta,L,\I)$ such that if $\ell(i_p)>F$ holds for all $1\leq p\leq \tau$, and $\mbox{\em \HH}_{\T} \geq H$ and $\mbox{\em \Gap}(\I) \leq G\,
\mbox{\em \HH}_{\T}$ also hold, then the number of eigenvalues of
$\mathcal{T}$ that have an algebraic conjugate in $\I$ is at least
$\mu \mbox{\em \HH}$, where $\mbox{\em \Gap}(\I):=\mbox{\em \Gap}_{\g,\ell}(\I)$ is as defined
in (\ref{def-gap}).
\end{prop}

\pf Suppose that $\kappa$, $\lambda$, $\epsilon$, $C$, $\g$, $\Delta$, $L$, $\I$, $\ell$ and $\T$ are as in the statement of the proposition, and put $\HH:=\HH_{\T}$ and $D:=D_{\T}$. Let $\mu$ be any positive real number.\\
In view of Theorem~\ref{limit-thm}, there exists a constant
$M_1:=M_1(\kappa,\epsilon,\mu,\g,\I)>0$ such that for any positive
real number $\zeta$ satisfying $\zeta< M_1$,
\begin{equation}\label{upsilon-0}
\Upsilon_{\kappa,\zeta}\leq \frac{1}{2+\frac{48 \pi \kappa \mu
g}{\epsilon |\I|}}
\end{equation}
holds (cf. (\ref{def-Gamma})). Put
\begin{equation}\label{zeta_0-def}
\zeta_0=\zeta_0(\kappa,\epsilon,\mu,\g,\I):=\min\left\{|\I|,\frac{M_1}{2}
\right\}~\mbox{  and  }\Upsilon:=\Upsilon_{\kappa,\zeta_0}.
\end{equation}
Then by (\ref{upsilon-0}) and (\ref{zeta_0-def}),
\begin{equation}\label{zeta_0-eq}
\Upsilon \leq \frac{1}{2+ \frac{48 \pi \kappa \mu g}{\epsilon
|\I|}}<\frac{1}{2}.
\end{equation}
By Lemma \ref{no-thy-2} (i) and Remark \ref{upsilon>0}, there exists a constant
$M_2:=M_2(\kappa,\epsilon,\mu,\g,\I)>0$ such that
$$\left|
\left\{p(x)\in \mathcal{P}_{\kappa} \,:\, \deg(p(x))\leq
\frac{1}{\Upsilon}\right \}\right| \leq M_2$$ holds, and therefore
\begin{equation}\label{main-prop-eq-1}
\left|  \left\{x \in  \I\cap \mathcal{E}_{\T} \,:\, \deg(x)\leq
\frac{1}{\Upsilon}\right \}\right| \leq \frac{M_2}{\Upsilon},
\end{equation}
where $\deg(x)$ is the degree of the minimal polynomial of
an algebraic number $x$ (cf. (\ref{def-Eg})).\\

Now, let $F:=F(\kappa,\g,\Delta,L,\I)$,
$C_1:=C_1(\kappa,\lambda,\epsilon,\mu,\g,\Delta,L,\I)$ and $C_2:=C_2(\kappa,\lambda,\epsilon,\mu,\g,\I)$ be the positive constants given by Theorem
\ref{close-evs-gap(I)} by taking $\delta:=\frac{\zeta_0}{2}$, and put
\[H:=\max\left\{\frac{6 \kappa\pi g}{\epsilon |\I|},\,\frac{24 \kappa \pi g M_2}{\epsilon \Upsilon |\I|},\,C_1
\right\}\mbox{~and~}G:=C_2.\] We now show that for this choice of
$F$, $H$ and $G$, the proposition holds. To this end, let $\theta$
be any element in $\mathcal{E}_{\T}\cap \I$ satisfying
$\deg(\theta)>\frac{1}{\Upsilon}$, and let $p_{\,\theta}(x)\in
\mathcal{P}_{\kappa}$ be a minimal polynomial of $\theta$. Then by
Theorem \ref{close-evs-gap(I)}, all roots of $p_{\theta}(x)$ must
lie in the closed interval $[\theta-\frac{\zeta_0}{2},
\theta+\frac{\zeta_0}{2}]$. Hence, by (\ref{upsilonbase}),
(\ref{def-Gamma}) and
$\deg(\theta)=\deg(p_{\theta})>\frac{1}{\Upsilon}$,
\begin{equation}\label{main-prop-eq-2}
\left|  \{x\in \I \,:\, p_{\,\theta}(x)=0\}\right|\leq \left|
\left\{x\in \left[\theta-\frac{\zeta_0}{2},
\theta+\frac{\zeta_0}{2}\right] \,:\, p_{\,\theta}(x)=0\right
\}\right| \leq \Upsilon\, \deg(p_{\theta})+1 < 2\,\Upsilon\,
\deg(p_{\,\theta}).
\end{equation}
Now, we prove the following claim.
\begin{claim}\label{main-prop-claim1}
The number of eigenvalues of $\T$ in $\I$ is at least $\left(
\frac{\epsilon |\I|}{12 \kappa \pi g}\right) \,\mbox{\em \HH}$.
\end{claim}
{\noindent{\bf Proof of Claim \ref{main-prop-claim1} :}} As $\HH\geq
H$ and $|\I|<\kappa$ (by (W1)), $\HH \geq \frac{6\kappa \pi
g}{\epsilon |\I|}>\frac{1}{\epsilon}$ holds. Hence, as $\Len(\I)\geq
\epsilon \HH>1$ (by statement (iii) of the proposition), there
exists $m\in \{2,\ldots,g\}$ so that $\ell(m)> \frac{\epsilon
\HH}{g}$ and $\I \subseteq I_m$ hold, where $I_m$ is the $m$\,th
guide interval. Put
$(\gamma,\alpha,\beta):=(\gamma_{m},\alpha_{m},\beta_{m})$,
$\ell:=\ell(m)$ and $e:=\left| \{j\in \{1,\ldots, \ell\} \,:\,
\alpha+2\sqrt{\beta\gamma}\cos\left(\frac{j\pi}{\ell+1}\right)\in \I
\}\right|$. Note that $\I \subseteq I_m$ and
$\alpha+2\sqrt{\beta\gamma}\cos\left(\frac{j\pi}{\ell+1}\right)\in
I_{m}$, for all $1\leq j\leq \ell$. Since
\begin{eqnarray*}
\left(\alpha+2\sqrt{\beta\gamma}\cos\left(\frac{(j-1)\pi}{\ell+1}\right)\right)
-\left(\alpha+2\sqrt{\beta\gamma}\cos\left(\frac{j\pi}{\ell+1}\right)\right)
\le \frac{2\pi \sqrt{\beta\gamma}}{\ell+1}\le \frac{\kappa
\pi}{\ell+1}
\end{eqnarray*}
holds for all $2\le j\le \ell$, it follows by $\ell >\frac{\epsilon
\HH}{g}$ and $\HH\geq \frac{6 \kappa \pi g}{\epsilon |\I|}$ that
$e\ge \lfloor \frac{(\ell+1)|\I|}{\kappa \pi} \rfloor > \lfloor
\frac{\epsilon \HH |\I|}{\kappa \pi g}\rfloor\geq \frac{\epsilon \HH
|\I|}{2\kappa \pi g } \geq 3$. Hence by Lemma \ref{two-thm-3.2}
(ii), there exists an eigenvalue $\theta \in
\mathcal{E}_{\mathcal{T}} \cap \I$ and, moreover,
\[\left|  \mathcal{E}_{\T}\cap \I
\right| \geq \left\lfloor\frac{e}{3}\right\rfloor > \frac{e}{6}
>\frac{\epsilon \HH  |\I|}{12 \kappa \pi g}\]
holds. Claim \ref{main-prop-claim1} now follows immediately. \epf

By applying Claim \ref{main-prop-claim1}, (\ref{main-prop-eq-1}) and
$\HH\geq \frac{24 \kappa \pi g M_2}{\epsilon \Upsilon |\I|}$ (by
$\HH\geq H$), it follows that
\begin{equation}\label{I-cap-high-degree}
\left| \{x \in \mathcal{E}_{\mathcal{T}}\cap \I \,:\,
\deg(x)>\frac{1}{\Upsilon} \}\right| \geq \frac{\epsilon |\I|
\HH}{12 \kappa \pi g}-\left| \{x \in \mathcal{E}_{\mathcal{T}} \cap
\I \,:\, \deg(x) \leq \frac{1}{\Upsilon}\}\right| \geq
\frac{\epsilon |\I| \HH}{24 \kappa \pi g}.
\end{equation}

Now, for each integer $i>\frac{1}{\Upsilon}$, let $\Delta_i$ be
the set of those elements in $\mathcal{E}_{\T}\cap \left([-\kappa,
\kappa]\setminus \I
 \right)$ of degree $i$ that have an algebraic conjugate which is contained in
 $\I$, and let $\Theta_i$ be the set of those elements in $\mathcal{E}_{\T}\cap \I$ that have degree $i$. Then by
(\ref{zeta_0-eq}) and (\ref{main-prop-eq-2}), each element in
$\Theta_i$ has an algebraic conjugate in $[-\kappa,
\kappa]\setminus \I$. This implies that $\Delta_i$ is a non-empty
set if and only if $\Theta_i$ is a non-empty set. Hence, for each integer
$i>\frac{1}{\Upsilon}$ satisfying $\Theta_i\neq \emptyset$, the number of elements in the set
\[\Lambda_i:=\{(\theta,\eta)\in \Theta_i\times \Delta_i \,:\, \mbox{$\theta$ and
$\eta$ are conjugate algebraic numbers}\}\]
is bounded above and below as follows:
\begin{equation}\label{Lambda_i}
(1-2\,\Upsilon)\,i\,\left|\Theta_i \right| < \left| \Lambda_i
\right| < 2\,i\,\Upsilon\, \left|\Delta_i \right|.
\end{equation}
Hence, by (\ref{zeta_0-eq}), (\ref{I-cap-high-degree}) and (\ref{Lambda_i}),
the inequality
\begin{equation}\label{main-prop-key-eq}
\frac{(1-2\,\Upsilon)\epsilon |\I|\,\HH}{24 \kappa \pi g}\leq
(1-2\,\Upsilon)\sum_{i>\frac{1}{\Upsilon},~\Theta_i\neq
\emptyset}\left|\Theta_i\right| <
2\,\Upsilon\sum_{i>\frac{1}{\Upsilon},~\Delta_i\neq \emptyset}\left|
\Delta_i \right|
\end{equation}
holds, and therefore by (\ref{zeta_0-eq}) and (\ref{main-prop-key-eq}), it follows that
\begin{equation}\label{last-main-prop}
\sum_{i>\frac{1}{\Upsilon},\,\Delta_i\neq \emptyset}\left| \Delta_i
\right| >\frac{(1-2\,\Upsilon)\epsilon |\I| \,\HH}{48 \kappa \pi
g\,\Upsilon} \geq \mu\, \HH
\end{equation}
holds. Since the number of eigenvalues of $\T$ which have an
algebraic conjugate in $\I$ is at least
$\sum_{i>\frac{1}{\Upsilon},\Delta_i\neq \emptyset}\left|
\Delta_i\right|$, the proposition now follows immediately by
(\ref{last-main-prop}). \epf

\pfmainthm Suppose that $\kappa$, $\lambda$, $\epsilon$, $C$, $\g$, $\Delta$, $L$, $\ell$ and $\T$ are as in the statement of the theorem, and put $\HH:=\HH_{\mathcal{T}}$, $\TT:=\TT_{\mathcal{T}}$, $D:=D_{\mathcal{T}}$, $\Len:=\Len_{\g,\ell}$ and $\Gap:=\Gap_{\g,\ell}$ (cf. (\ref{ell(I)}), (\ref{def-gap})). \\
By statement (iii) of the theorem and Lemma \ref{new-unimodular} (i), there exists an integer $s_0\in
\{2,\ldots,g \}$ such that $\frak{R}_{s_0}>\frak{R}_1$ and
$\ell(s_0)>\left(\frac{\epsilon}{g}\right)\HH$ (cf. (\ref{R-L})). On the other hand, by Corollary \ref{sub-WI}, there exists a
well-placed interval $\mathcal{J}_{0}$ in the $s_0$\,th guide
interval $I_{s_0}=(\frak{L}_{s_0},\frak{R}_{s_0})$ (relative to $\g$) such that
$\J_{0}\cap \mathcal{B}=\emptyset$ and
$\Len(\mathcal{J}_{0})>\left(\frac{\epsilon}{g}\right) \HH$ both hold as $\Len(\J_{0})\geq \ell(s_0)$ (cf. (\ref{def-guide interval}), (\ref{bad roots-def})). It follows by
Proposition \ref{main-prop} for $(\epsilon,\mu):=\left(
\frac{\epsilon}{g}, C(\kappa)+2 \right)$ that there exist
positive constants $F_{0}:=F_{0}(\kappa,\g,\Delta,L)$,
$G_{0}:=G_{0}(\kappa,\lambda,\epsilon,\g)$ and
$H_{0}:=H_{0}(\kappa,\lambda,\epsilon,\g,\Delta,L)$ such that
if $\ell(i)>F_{0}$ holds for all
$(\gamma_i,\alpha_i,\beta_i)\in \underline{\Delta}$ then
\[\mbox{either~~}\HH < H_{0}  \mbox{~~or~~}\Gap (\mathcal{J}_{0})>G_{0}\HH \mbox{~~holds},\]
as $\T$ has exactly $D+1$ distinct eigenvalues (cf. (\ref{def-Eg})) and $D\leq C\,\HH$ holds by statement (ii) of the theorem.\\

Now, if $\HH< H_{0}$, then the theorem
follows by taking $H:=H_{0}$ and $F:=F_{0}$.\\

Otherwise, $\HH \geq H_{0}$ and
$\Gap(\mathcal{J}_{0})>G_{0}\HH$ both hold, so by Corollary \ref{sub-WI} and Proposition \ref{gap-ind} for $\I:=\J_0$, there exists an integer $s_1$,
$\frak{c}(\g, \mathcal{J}_{0})\leq s_1 \leq \frak{d}(\g, \mathcal{J}_{0})$, and a well-placed interval $\mathcal{J}_{1}$ in the $s_1$\,th guide interval $I_{s_1}$ such that
\[\Len(\mathcal{J}_{1})>\frac{\Gap(\mathcal{J}_{0})}{g}>\left(\frac{G_{0}}{g}\right)\HH.\]
By applying Proposition \ref{main-prop} again for $(\epsilon,
\mu):=(\frac{G_{0}}{g}, C(\kappa)+2)$, there exist positive
constants $F_{1}:=F_{1}(\kappa,\g,\Delta,L)$,
$G_{1}:=G_{1}(\kappa,\lambda,\epsilon,\g)$ and
$H_{1}:=H_{1}(\kappa,\lambda,\epsilon,\g, \Delta,L)$ such
that if $\ell(i)>F_{1}$ holds for all
$(\gamma_i,\alpha_i,\beta_i)\in \underline{\Delta}$ then
\[\mbox{either~~}\HH < H_{1}  \mbox{~~or~~}\Gap(\mathcal{J}_{1})>G_{1}\HH \mbox{~~holds.} \]
Since $(\frak{R}_i)_{i=1}^{g}$ is a finite unimodal
sequence by Lemma \ref{new-unimodular}, it follows by iteratively repeating this argument (if necessary) that there exist an integer $m$, $1\leq m\leq g$, and positive constants $F_{j}:=F_{j}(\kappa,\g,
\Delta,L)$ and $H_{j}:=H_{j}(\kappa,\lambda,\epsilon,\g,
\Delta,L)~(0\leq j\leq m)$ given by Proposition \ref{main-prop} such that if
$\ell(i)>\max\{F_{j}\,:\,0\leq j\leq m\}$ holds for all
$(\gamma_i,\alpha_i,\beta_i)\in \underline{\Delta}$ then
\[
\max\{H_{j} \,:\, 0\leq j \leq m-1\}\leq \HH < H_{m}
\]
holds. Theorem \ref{b-i-small-thm} now follows by taking \[H:=\max\{H_{j} \,:\, 0\leq j \leq m\}\mbox{~~and~~}F:=\max\{F_{j}\,:\, 0\leq j\leq m\,\}.\] \epf

\section{Distance-Regular Graphs of Order $(s,t)$}\label{10}

In this section, we shall use our main result to show that, for
fixed integer $t >1$, there are only finitely many
distance-regular graphs of order $(s,t)$ whose smallest eigenvalue
is different from $-t-1$. We begin by recalling the relevant
definitions and some previous results.

Let $\Gamma$ be a distance-regular graph. For
any vertex $x$, the {\it local graph} of a vertex $x$ is
the subgraph of $\Gamma$ induced by $\Gamma_1(x)$.
For an integer $s\geq 1$, a {\em clique of size $s$} (or, $s$-clique)
is a set of $s$ vertices which are pairwise adjacent.
Following H. Suzuki (see \cite{suzukinotes}), we
say that a distance-regular graph $\Gamma$
is of {\em order $(s,t)$} for some positive integers $s, t$, if
the local graph of any vertex is the disjoint
union of $t+1$ cliques of size $s$.
In particular, a non-complete
distance-regular graph with valency $k \geq 3$ and $c_2 =1$ is
of order $(s,t)$ with $s = a_1 +1$ and $t = \frac{k}{a_1 +1}$.

Note that the
Hamming graph $H(n,q)$ is a distance-regular graph of order
$(n-1, q-1)$. Hence, for fixed positive integer $t$, there
are infinitely many distance-regular graphs of order $(s,t)$ where
$s$ is a positive integer. In addition,
B.~Mohar and J.~Shawe-Taylor \cite{t=1}
(see also \cite[Theorem 4.2.16]{bcn}) showed
that any distance-regular graph of order
$(s,1)$ with $s>1$ is isomorphic to the line graph of a Moore graph
or the point graph of some generalized $2D$-gon of order $(s,1)$,
where $D \in \{3,4,6\}$. Since the point graph of a generalized
$2D$-gon of order $(s, 1)$ is exactly the same as the flag graph of
a regular generalized $D$-gon of order $(s,s)$,
there are infinitely many distance-regular graphs
of order $(s,1)$ with $s> 1$.

The following proposition is well-known; we
include its proof for completeness.

\begin{prop}\label{st-prop}
For $s,t$ positive integers, let
$\Gamma$ be a distance-regular graph of order $(s,t)$
with diameter $D \geq 2$. Then the
smallest eigenvalue $\theta_D$ of $\Gamma$ satisfies $\theta_D \geq -t -1$.
Moreover, if $s > t$, then $\theta_D = -t-1$ holds.
\end{prop}

\pf Let $\mathcal C$ be the set of $(s+1)$-cliques in $\Gamma$.
Let $M$ be the vertex-clique of size $s+1$ incidence matrix, that
is, $M$ is the $(|V(\Gamma)| \times |\mathcal{C}|)$-matrix such
that the $(x,C)$-entry of $M$ is $1$ if $x \in C$ and 0 otherwise.
Then $MM^T = A + (t+1)I$, where $M^T$ is the transpose of $M$. As
$MM^T$ is positive semidefinite, it follows that all the
eigenvalues of $\Gamma$ are at least $-t-1$. Note that
$|\mathcal C|(s+1) = |V(\Gamma)|(t+1)$ so that if $s > t$ then
$|\mathcal{C}| < |V(\Gamma)|$, and, as the rank of $M$ is at most
$|\mathcal{C}|$, it follows that $A + (t+1)I$ is singular.
This shows that $A$ has $-t-1$ as its smallest eigenvalue. \epf

\begin{cor}\label{corst}
Let $t \geq 1$ be an integer. Then there are only finitely many
distance-regular graphs of order $(s,t)$ with $s \geq 1$ and $st
\neq 1$ which have smallest eigenvalue not equal to $-t-1$.
\end{cor}

\pf Let $t \geq 1$. If $\Gamma$ is a distance-regular graph of
order $(s, t)$ such that its smallest eigenvalue is different from
$-t-1$, then $s \leq t$ holds by Proposition~\ref{st-prop}.
As the valency of $\Gamma$ equals $s(t+1) \leq t(t+1)$,
the corollary follows by Theorem~\ref{b-i-thm} as long as
$s(t+1) \neq 2$. \epf

\begin{remark}
Not much is known concerning
distance-regular graphs of order $(s,t)$ with $t \geq 2$. The
distance-regular graphs of order $(1,2)$ and $(2,2)$
were classified by N.~L.~Biggs, A.~G.~Boshier and J.~Shawe-Taylor
\cite{$k=3$} and by A.~Hiraki, K.~Nomura and H.~Suzuki \cite{hiraki-suzuki},
respectively. In \cite{yamazaki}, N.~Yamazaki
presented some strong results
concerning distance-regular graphs of order $(s,2)$ with $s >2$.
However, it is not known whether there are infinitely
distance-regular graphs of order $(s,2)$ with $s \geq 2$ and $c_2 =1$.
\end{remark}

\section{Concluding Remarks}\label{11}

In Section~\ref{intro}, we mentioned that Sims' conjecture on
permutation groups could be used to prove that there are only
finitely many finite, connected distance-transitive graphs of
fixed valency greater than two. We conclude by recalling and
discussing a combinatorial version of Sims' conjecture that is
related to the Bannai-Ito conjecture.

To state this conjecture, we first
recall the definition of association schemes
(as defined by E.~Bannai and T.~Ito
\cite{banito}). An {\em association scheme}
$(X,\mathcal{R})$ is a finite set $X$
together with a collection $\mathcal{R}=\{R_0,R_1,\ldots, R_r\}$
of non-empty binary relations on
$X$ satisfying the following conditions: \\
(i) $\mathcal{R}$ is a partition of $X\times X$;\\
(ii) $R_0=\{(x,x)\,:\,x\in X\}$;\\
(iii) for each $R_i\in \mathcal{R}$, there exists $i'$
such that $R_{i'}=\{(y,x)\,:\,(x,y)\in R_i\}$;\\
(iv) for any $0 \leq i,j,h\leq r$ and for any $(x,y)\in R_h$, the
number $\left|\{z\in X\,:\,(x,z)\in R_i \mbox{~and~}(z,y)\in R_j\}
\right|$ is a constant $p^h_{ij}$ which depends only on $i,j,h$ not
on the choice of $(x,y)$.

Note that an association scheme in this sense is also called a
{\em homogeneous coherent configuration} (see \cite{higman}).
Also, an association scheme $(X,\mathcal{R})$ is called
{\em primitive} if any non-trivial relation $R_i~(i\neq 0)$
induces a directed connected graph on the vertex set $X$.

Let $(X,\mathcal{R})$ be a primitive association scheme. Then each
non-trivial relation $R_i\in \mathcal{R}~(i\neq 0)$ induces a
directed, connected, regular graph of valency $k_i:=p^0_{i i'}$.
L.~Pyber \cite[p.207]{pyber} and M.~Hirasaka \cite[p.105]{hirasaka} attribute
the following conjecture to L.~Babai.

\begin{conj}
[Babai's Conjecture]\label{babai}{\ \\}
There exists an integral function $f$ such that for any primitive
association scheme $(X,\{R_0,R_1,\ldots, R_r\})$,
\[k_{\max} \leq f(k_{\min})\] holds, where
$k_{\max}:=\max \{k_i\,:\,1\leq i\leq r\}$
and $k_{\min}:=\min \{k_i\,:\,1\leq i\leq r\}$.
\end{conj}

For a primitive permutation group $G$ on a finite set $\Omega$, the
orbits $R_i$ of the induced action of $G$ on $\Omega\times \Omega$
determine a primitive association scheme, denoted by $AS(G)$. Sims'
conjecture follows from Conjecture~\ref{babai} by considering the
association scheme $AS(G)$ for a primitive permutation group $G$.
Note also
that the cyclotomic schemes (for a definition see \cite[p.106]{hirasaka})
provide examples of primitive association schemes with
fixed smallest non-trivial valency and
an unbounded number of classes. Therefore, in Conjecture~\ref{babai}
we cannot expect to provide a bound for $r$ in terms of $k_{\min}$.

The main theorem of this paper, Theorem~\ref{b-i-thm}, implies that
Conjecture~\ref{babai} is true for primitive distance-regular
graphs with diameter $D$ as the sequence
$\left(k_i\right)_{1\leq i\leq D}$
is unimodal by \cite[Proposition 5.1.1 (i)]{bcn} and
$k_i \geq \sqrt{k}$ holds for all $i \geq 1$ by \cite[Proposition 5.6.1]{bcn}.

One could also ask whether there exists an integral function $f$
such that for any primitive commutative association scheme
$(X,\{R_0, R_1, ..., R_r\})$ with multiplicities
$m_i~(i=0,1,\ldots,r)$ with $m_0=1$, \[m_{\max}\leq f(m_{\min})\]
holds, where $m_{\max}:=\max \{m_i\,:\,1\leq i\leq r\}$ and
$m_{\min}:=\min \{m_i\,:\,1\leq i\leq r\}$. Such a function is not
known to exist even for the class of $Q$-polynomial association
schemes (for a definition see \cite[p.58]{bcn}), although the dual
statement of Theorem~\ref{b-i-thm} has been shown to be true by
W.~J.~Martin and J.~S.~Williford \cite{dual-BI-conj}. In
particular, they showed that for any $m_1 >2$, there are only
finitely many $Q$-polynomial association schemes with the property
that the first idempotent in a $Q$-polynomial ordering has rank
$m_1$.

\begin{center}
{\bf Acknowledgements}
\end{center}
The first author was supported by the Korea Research Foundation
Grant funded by the Korean Government (MOEHRD, Basic Research
Promotion Fund) KRF-2008-359-C00002. The third author was
partially supported by a grant of the Korea Research Foundation
funded by the Korean Government (MOEHRD) under grant number
KRF-2008-314-C00007 and he was also partially supported by
the Basic Science Research Program through the National Research
Foundation of Korea (NRF) funded by the Ministry of Education,
Science and Technology (grant number 2009-0089826). 

\end{document}